# PRECISE ASYMPTOTICS OF SMALL EIGENVALUES OF REVERSIBLE DIFFUSIONS IN THE METASTABLE REGIME

By Michael Eckhoff

*Universität Zürich*

We investigate the close connection between metastability of the reversible diffusion process $X$ defined by the stochastic differential equation

$$dX_t = -\nabla F(X_t)\,dt + \sqrt{2\varepsilon}\,dW_t, \qquad \varepsilon > 0,$$

and the spectrum near zero of its generator $-L_\varepsilon \equiv \varepsilon\Delta - \nabla F \cdot \nabla$, where $F : \mathbb{R}^d \to \mathbb{R}$ and $W$ denotes Brownian motion on $\mathbb{R}^d$. For generic $F$ to each local minimum of $F$ there corresponds a metastable state. We prove that the distribution of its rescaled relaxation time converges to the exponential distribution as $\varepsilon \downarrow 0$ with optimal and uniform error estimates. Each metastable state can be viewed as an eigenstate of $L_\varepsilon$ with eigenvalue which converges to zero exponentially fast in $1/\varepsilon$. Modulo errors of exponentially small order in $1/\varepsilon$ this eigenvalue is given as the inverse of the expected metastable relaxation time. The eigenstate is highly concentrated in the basin of attraction of the corresponding trap.

## Contents

**1. Introduction.** We address in this work the problem of characterizing—in terms of potential theoretic quantities—the low-lying spectrum of the following second-order, elliptic differential operator:

$$(1.1) \qquad L_\varepsilon \equiv -\varepsilon e^{F/\varepsilon} \nabla \cdot e^{-F/\varepsilon} \nabla = -\varepsilon \Delta + \nabla F \cdot \nabla, \qquad \varepsilon > 0,$$

on $L^2(\mathbb{R}^d, e^{-F/\varepsilon} dx)$, where the precise conditions on $F : \mathbb{R}^d \to \mathbb{R}$ are given in Assumption 1.2. Our main motivation is to derive precise uniform control in the limit $\varepsilon \downarrow 0$ of the distribution of *metastable transition times* $\tau(x)$ of the diffusion process $(X_t^x)$ on $\mathbb{R}^d$ generated by $-L_\varepsilon$, that is, the solution to the stochastic differential equation

$$(1.2) \qquad dX_t^x = -\nabla F(X_t^x) \, dt + \sqrt{2\varepsilon} \, dW_t, \qquad X_0^x = x.$$

Here $(W_t)$ denotes Brownian motion on $\mathbb{R}^d$ starting in zero. By definition $\tau(x)$ is the first time of a transition from the basin of attraction corresponding to a given local attractor $x$ of $\nabla F$, that is, a local minimum of $F$, to small vicinities of the more stable local attractors. The precise definition of $\tau(x)$ is given in (1.16).

We continue the work started in [4] and generalize the analysis of [3] from the discrete to the continuous state space setting. To each local attractor $x$ there corresponds a simple eigenvalue $\lambda_x$ of $L_\varepsilon$ which is exponentially small in $1/\varepsilon$. Modulo this type of error this eigenvalue equals the inverse of the expectation of $\tau(x)$. With the same precision an eigenfunction corresponding to $\lambda_x$ is constant in the basin of attraction of $x$ and exponentially small in "deeper" basins which correspond to attractors $y$ satisfying $\lambda_y < \lambda_x$. The results obtained in [4] then yield in terms of $F$ the leading-order asymptotic behavior of these eigenvalues. Moreover, below some threshold of order $\varepsilon^N$ no other eigenvalues occur. The control of the low-lying part of the spectrum implies that the rescaled (by its expectation) distribution of a metastable transition time converges—again modulo in $1/\varepsilon$ exponentially small errors—to the exponential distribution with parameter 1.

Metastability in random dynamical systems is an intensively studied phenomenon. A Markov process in the metastable regime, roughly speaking, exhibits quasi-invariant sets of the state space, which may be viewed as metastable states, in which the process is captured for long time periods. For systems with discrete state space in this regime in [22] and [23] as well as [30] and [7] the authors study different aspects in this area. Concerning systems with continuous state space in this regime, we refer the reader to [20, 33, 40, 41], where the authors develop a large deviation technique for diffusion processes to study spectral and dynamical properties. From the point of view of asymptotic expansions in the small parameter $\varepsilon > 0$ we mention [5, 6, 13, 14, 15, 19, 21, 30, 31, 32]. In most of these papers the authors consider the process up to the time of exit from a single domain of attraction associated to the unperturbed dynamical system. For the investigation of the spectrum and its connection to metastability it is necessary to consider the process as it continues from one domain to another. In [9, 10, 11, 15, 21, 23] and in [5, 6, 32, 34], where in the latter two articles the full description of the low-lying spectrum is accomplished, the authors investigate properties of the spectrum of the generator of the dynamical system that are connected to metastability. Unfortunately, these approaches encounter the following shortcomings. Generally speaking, rigorous asymptotic expansions, though giving sharp error estimates, suffer from strong regularity assumptions. On the other hand, $L^2$-methods as applied in [34] and [9, 10, 11] as well as large deviation theory lead to rough error estimates. In [3] and in [18] we establish the characterization of the low-lying spectrum in the context of Markov chains in the metastable regime. A key idea of [3, 4] and [17, 18] for irreversible chains is to analyze metastability from the dynamical or from the spectral point of view by potential theoretic methods, which particularly leads to a clear description of the spectrum in terms of the geometry of $F$. In addition to the work in [5, 6] and in [32] we are able to establish the same precise relation of the small eigenvalues to the geometric properties of $F$. Our approach also considerably improves the range of applicability as well as the quality of the error estimates. In [17, 18] this aspect is particularly emphasized. Here we concentrate on the main new technical complications which do not exist in systems with finite, discrete state space.

The technical tool to connect spectral to potential theory already appears in [40] or in [35], relying on work of [44], and was rediscovered in [3]. Reference [40] contains a description of the spectrum in terms of the underlying Markovian process while in [44] and [35] the analytical counterpart is used to investigate criticality of elliptic operators. This characterization is far more transparent for processes with discrete state space as is demonstrated in [3] and in [17, 18]. The fact that a point in discrete space can be visited by the process with strictly positive probability, that is, has strictly positive capacity, might be seen as a main reason for this difference. In continuous state spaces small balls are the equivalent of points in discrete spaces. This choice entails the disadvantage that a function a priori may change its sign on a small ball. Using level sets of functions instead of balls, one quickly runs into technical complications whose solutions go beyond the questions we are addressing. The approach presented in the previously mentioned references naturally requires to work in the $L^\infty$-context. We hence shall follow the strategy to first establish rather strong pointwise $L^\infty$-estimates. Compared to [3] and [34] a second complication arises from the fact that the state space is noncompact. A treatment of the analogous problem concerning irreversible, infinite-state Markov chains can be found in [18]. There is a well-established $L^2$-theory of weighted estimates of solutions of second-order elliptic differential equations as developed in [1] or in [26, 27] involving a small parameter from which pointwise bounds on solutions can be obtained. The development of weighted estimates will serve to gain control of the growth of eigenfunctions at infinity. We would like to mention that the methods introduced in [4] and in [18] suffice to prove the same kind of estimates for which $L^2$-weighted estimates are not available, even if the process is irreversible though technically simpler.

The a priori input enables us to relate small eigenvalues of $L_\varepsilon$ to the capacity matrix introduced in [28]. For generic $F$ the analysis of this matrix then is a straightforward generalization of that in [3] and [18] for Markov chains. Let us mention that this matrix representation also can be used to treat the degenerate situation, where there exist attractors of $\nabla F$ of equal strength with respect to each other. It turns out that to each small eigenvalue there corresponds a quasi-invariant set and a time scale, which roughly speaking equals the expected time the process generated by $L_\varepsilon$ is captured this set. These time scales are defined in terms of capacities and the invariant measure of the process. As is shown in [2] and [18] in the discrete state space setting they determine the long-time behavior of the process in a precise manner. These kind of results were extended in [4] to the diffusion process generated by $L_\varepsilon$. They will serve as a crucial tool to investigate small eigenvalues of $L_\varepsilon$.

We now recall the main potential theoretic background. A set $\Gamma$ with locally $\mathcal{C}^{2,\alpha}$ boundary for some $\alpha > 0$ henceforth will be referred to as a *regular* set. Fix disjoint, nonempty closed regular sets $A, B \subset \mathbb{R}^d$ such that $\Gamma \equiv \mathbb{R}^d \setminus A \setminus B$ is connected (usually in the sequel $A$ and $B$ are balls). The $\lambda$-capacity of the capacitor $(A, B)$ is given by

$$(1.3) \qquad \mathrm{cap}_A^\lambda(B) \equiv \varepsilon \int_{\partial A} e^{-F/\varepsilon} \partial_n h_{A,B}^\lambda \, d\sigma - \lambda \int_A e^{-F/\varepsilon} \, dy,$$

where locally there is $\alpha > 0$ such that $F : \mathbb{R}^d \to \mathbb{R}$ is $\mathcal{C}^{1,\alpha}$ for some $\alpha > 0$, $\sigma$ always denotes the Euclidean surface measure on the set the integration is taken over, $n$ is the unit normal at this surface pointing towards $A \cup B$ and the normal derivative is taken from outside $A$ and $B$. Here $h_{A,B}^\lambda$ denotes the electrostatic equilibrium potential of the capacitor, that is, the weak solution $h \in W^{1,2}(\Gamma, e^{-F/\varepsilon} \, dx)$ of the Dirichlet problem

$$(1.4) \quad (L_\varepsilon - \lambda) h(x) = g(x), \qquad x \in \Gamma, \qquad h - f \in W_0^{1,2}(\Gamma, e^{-F/\varepsilon} \, dx),$$

where $\Gamma \equiv \mathbb{R}^d \backslash A \backslash B$, $g \equiv 0$, $f \equiv \mathbb{1}_A$ and where $W_0^{1,2}(\Gamma, e^{-F/\varepsilon} \, dx)$ denotes the closure of $\mathcal{C}_0^1(\Gamma)$ in $W^{1,2}(\Gamma, e^{-F/\varepsilon} \, dx)$, the space of weakly differentiable functions with first partial derivatives in $L^2(\Gamma, e^{-F/\varepsilon} \, dx)$. Under Assumption 1.2 standard regularity theory will show that (1.4) is uniquely solvable and that the solution is $\mathcal{C}^{2,\alpha}$ up to the boundary. Functions $h$ satisfying (1.4) for some $f$ and $g \equiv 0$ we sometimes refer to as (weakly) $(L_\varepsilon - \lambda)$-*harmonic* functions (with respect to the measure $e^{-F/\varepsilon} \, dx$). In the commonly used terminology of partial differential equations they are called weakly $(\Lambda_\varepsilon - e^{-F/\varepsilon} \lambda)$-harmonic (with respect to Lebesgue measure), where we introduce the formally symmetric, locally elliptic, second-order differential operator in divergence form

$$(1.5) \qquad \Lambda_\varepsilon \equiv -\varepsilon \nabla \cdot e^{-F/\varepsilon} \nabla = e^{-F/\varepsilon} L_\varepsilon.$$

In particular, the well-developed regularity theory for divergence-type operators is available. The *communication height* between sets $A$ and $B$ is defined by

$$(1.6) \qquad \hat{F}(A, B) \equiv \inf_{\substack{c \colon [0,1] \to \mathbb{R}^d \\ c(0) \in A, c(1) \in B}} \max F(c([0,1])),$$

where the infimum is taken over all continuous curves. If $A \equiv \{x\}$ is a singleton, for convenience we write $\hat{F}(x, B) \equiv \hat{F}(\{x\}, B)$ instead. Furthermore, for a finite set of points $I \cup x$ such that $B_{I \cup x}$ is a disjoint union of open balls, where

$$(1.7) \qquad B_J \equiv \bigcup_{y \in J} B(y, \varepsilon/4), \qquad J \subset \mathbb{R}^d,$$

we introduce

$$(1.8) \qquad A_{x,I} \equiv \{y \in \mathbb{R}^d | \hat{F}(y, x) < \hat{F}(y, I \backslash x)\}.$$

In analogy to [18] we define the time scales

$$(1.9) \qquad T_{x,I} \equiv \frac{\int_{A_{x,I}} e^{-F/\varepsilon} \, dy}{\mathrm{cap}^0_{B_x}(B_{I \backslash x})}.$$

We recall from Theorem 3.1 in [4] the classical Eyring formula for the capacity.

THEOREM 1.1. *Fix regular, disjoint, nonempty sets $A$ and $B$. Assume that there is only one solution of $F(z^*) = \hat{F}(A, B) > \max F(A \cup B) + R\varepsilon \log(1/\varepsilon)$ such that $z^*$ is a critical point of $F$. If in addition to the condition $F \in \mathcal{C}^{1,\alpha}$ for some $\alpha > 0$ the Hessian at $z^*$ of $F$ exists and is nondegenerate, then for some $R$,*

$$(1.10) \quad \mathrm{cap}^0_A(B) = (1 + \mathcal{O}(1)\varepsilon \log(1/\varepsilon)) \frac{(2\pi)^{d/2-1}|\lambda^*|}{\sqrt{|\det \mathrm{Hess}\, F(z^*)|}} \varepsilon^{d/2} e^{-\hat{F}(A,B)/\varepsilon},$$

where $\lambda^*$ is the unique, negative eigenvalue of the Hessian at $z^*$. The modulus of the Landau symbol is dominated by a constant $C \equiv C(d, F)$.

Let $\mathcal{M}$ denote the set of *local minima* of $F$. For $x \in \mathcal{M}$ and $I \subset \mathcal{M}\setminus x$ nonempty with nondegenerate Hessian at $x$ and $z^*$ as in Theorem 1.1, we obtain from (1.10) that the time scale introduced above satisfies

$$(1.11) \quad T_{x,I} = (1 + \mathcal{O}(1)\varepsilon \log(1/\varepsilon))2\pi\sqrt{\frac{|\det \operatorname{Hess} F(z^*)|}{|\lambda^*|\det \operatorname{Hess} F(x)}}e^{(\hat{F}(x,I)-F(x))/\varepsilon}.$$

Let us now describe the main results of this paper. We have to introduce some more notation. For a regular domain $\Sigma$ let $L_\varepsilon^\Sigma$ denote the self-adjoint operator with Dirichlet boundary conditions corresponding to the quadratic form

$$(1.12) \qquad q_\varepsilon^\Sigma(h) \equiv \varepsilon \int_\Sigma e^{-F/\varepsilon} |\nabla h|^2 \, dx$$

of the operator $L_\varepsilon$ on $L^2(\Sigma, e^{-F/\varepsilon} dx)$ with domain $W_0^{1,2}(\Sigma, e^{-F/\varepsilon} dx)$. Denote the principal eigenvalue of the Dirichlet operator $L_\varepsilon^\Sigma$ by

$$(1.13) \qquad \lambda(\Sigma) \equiv \inf \sigma(L_\varepsilon^\Sigma),$$

where $\sigma(L_\varepsilon^\Sigma)$ is the spectrum of $L_\varepsilon^\Sigma$. In the sequel we impose the following conditions on $F$.

ASSUMPTION 1.2. $F \in W_{\text{loc}}^{2,\infty}(\mathbb{R}^d) \cap \mathcal{C}^1(\mathbb{R}^d)$ and $\nabla F$ is locally Hölder continuous. There are constants $c > 0$ and $C_1$ satisfying $\inf_{\{F>C_1\}} |\nabla F| \geq c$. Moreover, $\mu_\varepsilon \equiv \lambda(\{F > C_1\}) \geq \delta$ for some $\delta > 0$ independent of small $\varepsilon > 0$.

Let us remark that the conditions under which (1.11) holds are not at all borderline to our approach. In fact, the only condition we need is that either $\delta T_{x,I\setminus x} > T_{y,I\setminus y}$ or $T_{x,I\setminus x} < \delta T_{y,I\setminus y}$, where $x, y \in \mathcal{M}$, $x \neq y$, $I \subset \mathcal{M}$, for some $\delta > 0$. In particular, as long as this condition is valid we can deal with all kinds of degenerate situations of $F$ in the relevant regions. This just leads to different asymptotic behaviors of $T_{x,I}$. We refer the reader to [17, 18] where in the context of Markov chains precise, minimal conditions on the time scales are given. We would also like to mention that Assumption 1.2 can be weakened in several directions. The condition $\inf_{\{F>C_1\}} |\nabla F| \geq c$ can be relaxed to the existence of a subset $\tilde{\mathcal{M}}$ of the set of local minima $\mathcal{M}$ of $F$ with the property that $\delta \min_{x \in \tilde{\mathcal{M}}} T_{x,\mathcal{M}\setminus\tilde{\mathcal{M}}} > \sup_{y \in \mathcal{M}\setminus\tilde{\mathcal{M}}} T_{x,\tilde{\mathcal{M}}}$ for some $\delta > 0$. In particular, $F$ may have infinitely many local minima where the minima in $\mathcal{M}\setminus\tilde{\mathcal{M}}$ are not as "deep" as those in $\tilde{\mathcal{M}}$. Moreover, the analysis works for a large class of functions $F \equiv F_\varepsilon$ depending on $\varepsilon$ also [for further comment concerning this point see the remark after (4.47)]. One could further considerably relax the regularity assumptions on $F$. It is also possible to study the irreversible situation where $\nabla F$ is replaced by a general vector field $b$. Finally, a generalization to Riemannian manifolds is straightforward.

The condition on the principal eigenvalue is quite natural and flexible. If, for example, $F$ is in addition $\mathcal{C}^2$ and $\limsup_{|x|\to\infty} |\Delta F(x)|/|\nabla F(x)|^2 < \infty$, it is easy to see that $\mu_\varepsilon > \delta/\varepsilon$ for some $\delta > 0$. For $F \equiv F_\varepsilon$ depending on $\varepsilon$ the bound on $\mu_\varepsilon$ can be replaced by, for example, $\varepsilon^M$ for some constant $M$ or (even exponentially small in $1/\varepsilon$ with small rate depending on the geometry of $F$ in $\{F < C_1\}$). If $F$ is uniformly strictly convex outside some convex set, one could use Brascamp–Lieb's inequality to show that $\mu_\varepsilon \geq \inf_{\{F>C_1\}} \min(\sigma(\nabla\nabla^t F))$. As we only focus on the new technical complications in the continuous state space setting we do not aim at the most general conditions under which the analysis works.

Assumption 1.2 implies that $F$ has local uniform, exponentially tight level sets, that is, $\int_{\{F>\alpha\}} e^{-F/\varepsilon} dx \leq C e^{-\alpha/\varepsilon}$ for some constant $C \equiv C(d, |\{F \leq \alpha\}|)$. Indeed, for a point $x \in \{F > \alpha\}$ the solution $\gamma$ to $\dot\gamma(t) = \nabla F(\gamma(t))$ with $F(\gamma(0)) = \alpha$ and $\gamma(T) = x$ we may estimate $F(x) - \alpha = \int_0^T |\nabla F(\gamma(t))| |\dot\gamma(t)| dt \geq c \operatorname{dist}(x, \{F > \alpha\})$ for $\alpha \geq C_1$. Therefore, in this work we may use in compact ($\varepsilon$-independent) sets (obvious generalization from $F$ being $\mathcal{C}^2$ to $F$ being $\mathcal{C}^1$) the results given in [4].

The first result, stated in Theorem 4.2 and referred to as the sharp uncertainty principle, is strikingly reminiscent of the uncertainty principle in quantum mechanics. We recall that the tunneling time of a quantum-mechanical particle moving in a double-well potential approximately is given by the inverse of the spectral gap. Let $\lambda(\Omega)$ be the principal eigenvalue of the Dirichlet operator $L_\varepsilon^\Omega$ with zero boundary conditions on $\mathbb{R}^d \setminus \Omega$, where $\Omega$ is an open, regular set. Furthermore, introduce for a Borel set $B \subset \mathbb{R}^d$ the transition time

$$(1.14) \quad \tau_B^x \equiv \inf\{t \geq 0 | X_t^x \in B\} \quad \text{and write shorthand} \quad \tau_I^x \equiv \tau_{B_I}^x$$

of the diffusion given by (1.2) from $x$ to the union $B_I$ of small balls, defined in (1.7), which are centered at the points in $I$.

THEOREM 1.3. *Assume that $F$ satisfies Assumption* 1.2. *Then there exists $N \equiv N(d) \geq 0$ such that for all $\rho > N\varepsilon \log(1/\varepsilon)$, $x \in \mathcal{M}$, $I \subset \mathcal{M} \setminus x$ satisfying $T_{x,I} = T_I \equiv \max_{y \in \mathcal{M} \setminus I} T_{y,I} \geq e^{-\rho/\varepsilon} T_{I \cup x}$,*

$$(1.15) \quad \mathbb{E}[\tau_I^x] = (1 + \mathcal{O}(1)\varepsilon^{-N} e^{-\rho/\varepsilon}) \frac{1}{\lambda(\mathbb{R}^d \setminus B_I)} = (1 + \mathcal{O}(1)\varepsilon^{-N} e^{-\rho/\varepsilon}) T_I.$$

*Here the modulus of the Landau symbol is dominated by a constant $C \equiv C(d, N, F)$.*

We also are able to compute the limit law of the distribution of the rescaling $\tau(x)/\mathbb{E}[\tau(x)]$ of a *metastable transition time* $\tau(x)$, $x \in \mathcal{M}$, defined by

$$(1.16) \qquad \tau(x) \equiv \tau^x_{\mathcal{M}(x) \cup \Omega^c}, \qquad \mathcal{M}(x) \equiv \{y \in \mathcal{M} | F(y) < F(x)\},$$

where $\Omega \subset \mathbb{R}^d$ is a regular domain containing the set $\{F < C_1\}$. Let us define $\rho \equiv \rho(F, \varepsilon)$ by

$$(1.17) \quad e^{\rho/\varepsilon} \equiv \min\left\{\frac{T_{x,I\setminus x}}{T_{y,I\setminus y}} \Big| x, y \in \mathcal{M}, x \neq y, I \subset \mathcal{M}, T_{x,I\setminus x} \geq T_{y,I\setminus y}\right\}.$$

Then we have:

THEOREM 1.4. *Suppose that $F$ satisfies Assumption 1.2 with $\mu_\varepsilon \geq \delta\varepsilon$ for some $\delta > 0$. Assume furthermore that either $\Omega$ is bounded or $\int_{\{F>C_1\}} |\nabla F|^d \times e^{-(F-C_1)/\gamma} \, dy < \infty$ for some $\gamma > 0$. There exist $N \equiv N(d) \geq 0$ and $C \equiv C(d, F)$ such that for all $\rho > N\varepsilon \log(1/\varepsilon)$*

$$(1.18) \quad \begin{aligned} &\mathbb{P}[\tau(x) > T\mathbb{E}[\tau(x)]] \\ &= (1 + \mathcal{O}(1)\varepsilon^{-N} e^{-\rho/\varepsilon}) \exp(-(1 + \mathcal{O}(1)\varepsilon^{-N} e^{-\rho/\varepsilon})T), \end{aligned}$$

*where the modulus of the Landau symbol is bounded by $C$ uniformly in $\varepsilon$ and $T$.*

A more detailed version of this theorem is Theorem 5.2 [see also the remark following (5.2)].

The main ingredient to prove (1.18) is that besides principal eigenvalues we are able to analyze all other exponentially small eigenvalues and relate them to the metastable structure given by $F$. We then have [see (4.6) for a more detailed version]

THEOREM 1.5. *Assume that $F$ satisfies Assumption 1.2. There exist $N \equiv N(d) \geq 0$ and a constant $C \equiv C(d, F)$ such that for all $\rho > N\varepsilon \log(1/\varepsilon)$ the following holds:*

(i) *For every $x \in \mathcal{M}$ there exists a simple eigenvalue $\lambda_x$ of $L_\varepsilon$ such that*

$$(1.19) \qquad \lambda_x = (1 + \mathcal{O}(1)\varepsilon^{-N} e^{-\rho/\varepsilon}) \lambda(\mathbb{R}^d \setminus B_{\mathcal{M}(x)}),$$

*where $\mathcal{M}(x)$ is defined in (1.16).*

(ii) *Let $\mathcal{M}_x \equiv \{y \in \mathcal{M} | \lambda_y < \lambda_x\}$. There is an eigenfunction $\phi_x$ corresponding to $\lambda_x$, normalized by $\phi_x(x) \equiv 1$ and a set $\tilde{\mathcal{M}}_x$ of cardinality $|\mathcal{M}_x|$ such that $B(y, \sqrt{\varepsilon}) \cap \tilde{\mathcal{M}}_x$ is a singleton for all $y \in \mathcal{M}_x$ and for all $z \in \{F < C_1\}$*

$$(1.20) \qquad \begin{aligned} \phi_x(z) &= (1 + \mathcal{O}(1)\varepsilon^{-N} e^{-\rho/\varepsilon}) \mathbb{P}[\tau^z_x < \tau^z_{\tilde{\mathcal{M}}_x}] \\ &\quad + \mathcal{O}(1)\varepsilon^{-N} e^{-\rho/\varepsilon} \mathbb{P}[\tau^z_{\tilde{\mathcal{M}}_x} < \tau^z_x]. \end{aligned}$$

*Here the Landau symbols are bounded by $C$ in absolute value.*

(iii)

(1.21) $$\sigma(L_\varepsilon) \cap [0, \varepsilon^N) = \{\lambda_x | x \in \mathcal{M}\}.$$

Equation (1.19) in combination with (1.15) and (1.18) relates exponentially small eigenvalues of $L_\varepsilon$ to the metastable structure of the diffusion $X$. Furthermore, under the conditions required for (1.11) we have determined the leading asymptotic in (1.15).

Let us finally describe the organization of the paper. Using sharp Harnack- and Hölder-type estimates, in Section 2 we derive analogous estimates for a priori nonpositive harmonic functions. As a result we gain in Lemma 2.3 pointwise control on the oscillation of eigenfunctions corresponding to small eigenvalues in terms of suprema over suitable small balls close to the local minima of $F$. In Section 3 we prove bounds of those suprema by exploiting the strong drift of the diffusion toward local minima of $F$. The a priori input then gives precise control of eigenfunctions in compact sets. As soon as we have established this structural information we are in a position to generalize the analysis developed in the discrete state space setting to the diffusion setting. In particular, in Section 4 we relate the low-lying spectrum to the capacity matrix introduced in [28] and derive the asymptotic information in terms of the time scales introduced in (1.9). As a consequence we obtain the limit law of metastable transition times defined in (1.16).

**2. Pointwise asymptotics in bounded sets.** Fix an open, connected, regular set $\Omega$ and recall the definition $W_0^{1,2}(\Omega, e^{-F/\varepsilon} dx)$. This section is devoted to the following simple idea. A weak solution $\phi \in W_0^{1,2}(\Omega, e^{-F/\varepsilon} dx)$ of the eigenvalue problem

(2.1) $$(L_\varepsilon - \lambda)\phi(x) = 0, \qquad x \in \Omega,$$

with small energy $\lambda$ cannot create large oscillations everywhere in a region where $F$ is small. We start with the following.

2.1. *A priori bounds on principal eigenvalues.* Recall the definition of the $(L_\varepsilon - \lambda)$-equilibrium potential $h_{A,B}^\lambda$, $A$, $B$ closed and regular with connected complement $\mathbb{R}^d \backslash A \backslash B$, $\lambda \geq 0$, introduced in (1.4). Furthermore, let $w_{A,B}^\lambda$ be the solution of the Poisson problem (1.4) with $f \equiv 0$ and $g \equiv h_{A,B}^\lambda$. We also shall use the convention $h_A^\lambda \equiv h_{A,A}^\lambda$ and $w_A^\lambda \equiv w_{A,A}^\lambda$. Since $h_{A,B}^\lambda$ and $w_{A,B}^\lambda$ are weak solutions of the corresponding problem for the operator $\Lambda_\varepsilon - e^{-F/\varepsilon}\lambda$ defined in (1.5), Theorem 8.8 in [24] in combination with Theorem 9.19 in [24] show that the unique solutions if they exist are locally $\mathcal{C}^{2,\alpha}$ up to the boundary. Define for $K \subset (A \cup B)^c$

(2.2) $$s_K^\lambda(A, B) \equiv \sup_K \frac{w_{A,B}^\lambda}{h_{A,B}^\lambda}.$$

We abbreviate

$$
\begin{aligned}
s^\lambda(A,B) &\equiv s^\lambda_{(A\cup B)^c}(A,B), \\
s^\lambda_K(A) &\equiv s^\lambda_K(A, B\equiv A), \\
s^\lambda_{A^c} &\equiv s^\lambda_{A^c}(A, B\equiv A).
\end{aligned}
\qquad(2.3)
$$

Recall the definition of the self-adjoint operator $L_\varepsilon^\Sigma$ with Dirichlet boundary conditions at $\partial\Sigma$ corresponding to the quadratic form defined in (1.12) and its principal eigenvalue $\lambda(\Sigma)\equiv\inf\sigma(L_\varepsilon^\Sigma)$. For $\lambda\notin\sigma(L_\varepsilon^\Sigma)$ we denote by

$$G^\lambda_\Sigma \equiv (L_\varepsilon^\Sigma - \lambda)^{-1} \qquad(2.4)$$

the resolvent operator. A priori we have that positive kernel of the resolvent $G^\lambda_\Sigma$, defined by the semigroup of the solution $X^x$ of (1.2) for $\lambda < \lambda(\Sigma)$, is in $L^2(\Sigma^2, e^{-(F(x)+F(y))/\varepsilon}\, dx\, dy)$.

We refer to the lower bound in (2.5) on the principal eigenvalue as the uncertainty principle.

LEMMA 2.1. *Let $\Sigma$ be a bounded, regular, open, connected set. Then for all regular, closed sets $A$, $B$, such that $A\cup B = \Sigma^c$ it follows that*

$$\lambda(\Sigma) \geq \frac{1}{s^0(A,B)}. \qquad(2.5)$$

PROOF. We claim that the following variational formula of Donsker and Varadhan (see [15] or [36]) for the principal eigenvalue holds:

$$\lambda(\Sigma) = \inf_{\substack{f\in\mathcal{C}^1(\overline{\Sigma})\\ f|\partial\Sigma=0, \int_\Sigma f^2=1}} \sup_{\substack{u\in\mathcal{C}^2(\overline{\Sigma}), u|\partial\Sigma=0\\ u(x)>0, x\in\Sigma}} \int_\Sigma \frac{L_\varepsilon u(x)}{u(x)} f(x)^2\, dx. \qquad(2.6)$$

Since $\Lambda_\varepsilon$ satisfies the conditions of Theorem 8.6 in [24], we have $\lambda(\Sigma) > 0$. By the weak maximum principle Theorem 8.1 in [24] it follows that $G^0_\Sigma$ is a positive operator, that is, the kernel is nonnegative and thus strictly positive since $G^0_\Sigma$ is injective. Theorem XIII.44 in [37] tells us that $\lambda(\Sigma)$ is a simple eigenvalue and that an eigenfunction $\phi\in W^{1,2}_0(\Sigma, e^{-F/\varepsilon}\, dx)$ almost surely does not change sign. By the same arguments given before (2.2) this function is in $\mathcal{C}^{2,\alpha}(\overline{\Sigma})$. Inserting $u\equiv\phi$ on the right-hand side of (2.6) yields one inequality. On the other hand, for every $u$ in the class of functions the supremum is taken over, we may choose $f\equiv C(ue^{-F/\varepsilon}\phi)^{1/2}$ with normalizing $C$ such that $f^2$ is a density. We obtain the remaining assertion by inserting $f$ on the right-hand side of (2.6) since the integral equals $C^2\int_\Sigma L_\varepsilon u\phi e^{-F/\varepsilon}\, dx = C^2\lambda(\Sigma)\int_\Sigma u\phi e^{-F/\varepsilon}\, dx = \lambda(\Sigma)$. Here we have used that $L_\varepsilon$ is symmetric on $\mathcal{C}^{2,\alpha}(\overline{\Sigma})$ and that $L_\varepsilon\phi(x) = L_\varepsilon^\Sigma\phi(x)$.

To obtain (2.5), we simply insert $u\equiv w^0_{A,B}\in\mathcal{C}^{2,\alpha}(\overline{\Sigma})$ and use $L_\varepsilon u = h^0_{A,B}\in\mathcal{C}^{2,\alpha}(\overline{\Sigma})$ on $\Sigma$, using that both functions exist by Theorem 8.3 in [24]. □

From the variational principle, Theorems 4.5.2 and 4.5.1 in [12], we also obtain the following sharp upper bound as we shall see in Theorem 4.2.

LEMMA 2.2. *Let $\Sigma$ be a regular, open set such that* $\operatorname{dist}(\mathcal{M} \cap \Sigma, \partial \Sigma) > \rho$ *for some $\rho > 0$. Then for some $C \equiv C(d, F|\Sigma, \rho)$ and all $x \in \mathcal{M} \cap \Sigma$,*

$$\lambda(\Sigma) \leq (1 + Ce^{-\beta/\varepsilon}/\varepsilon) \frac{\operatorname{cap}^0_{B(x,\varepsilon)}(\Sigma^c)}{\int_{A^\beta_{x,\Sigma^c}} e^{-F/\varepsilon} \, dx}, \tag{2.7}$$

*where we have defined $A^\beta_{x,\Sigma^c} \equiv \{y \in \Sigma \mid \hat{F}(y, x) < \hat{F}(y, \Sigma^c) - \beta\}$ for all $\beta > 0$ such that $B(x, \rho/2) \subset A$. Here $\hat{F}$ denotes the communication height introduced in (1.6).*

PROOF. Insert $u \equiv h^0_{B(x,\varepsilon),\Sigma^c}$ with $x \in \mathcal{M} \cap \Sigma$ and by convention $h^0_{B(x,\varepsilon),\Sigma^c} \equiv 1$ on $B(x,\varepsilon)$ into the variational principle Theorems 4.5.2 and 4.5.1 in [12]:

$$\lambda(\Sigma) = \inf_{u \in W^{1,2}_0(\Sigma, e^{-F/\varepsilon} \, dx) \setminus 0} \frac{\int_\Sigma e^{-F/\varepsilon} |\nabla u|^2 \, dx}{\int_\Sigma e^{-F/\varepsilon} |u|^2 \, dx} \tag{2.8}$$

to obtain by Green's first formula

$$\lambda(\Sigma) \leq \frac{\operatorname{cap}^0_{B(x,\varepsilon)}(\Sigma^c)}{\int_{A^\beta_{x,\Sigma^c}} e^{-F/\varepsilon} (h^0_{B(x,\varepsilon),\Sigma^c})^2 \, dx}. \tag{2.9}$$

Invoking Corollary 4.8 in [4], we derive $h^0_{\Sigma^c, B(x,\varepsilon)}(y) \leq Ce^{-\beta/\varepsilon}/\varepsilon$ for some $C \equiv C(d, F | A^\beta_{x,\Sigma^c})$ and all $y \in A^\beta_{x,\Sigma^c} \setminus B(x, 2\varepsilon)$. This estimate in combination with the maximum principle shows $h^0_{B(x,\varepsilon),\Sigma^c}(y) \geq 1 - Ce^{-\beta/\varepsilon}/\varepsilon$ for all $y \in A^\beta_{x,\Sigma^c}$. This establishes (2.7). □

2.2. *Uniform regularity estimates for $(L_\varepsilon - \lambda)$-harmonic functions changing sign.* For a regular domain $\Sigma \subset \mathbb{R}^d$ and a function $f : \mathbb{R}^d \to \mathbb{R}$ we define the oscillation of $f$ in $\Sigma$ as

$$\operatorname{osc}_\Sigma f \equiv \sup_\Sigma f - \inf_\Sigma f. \tag{2.10}$$

We are now in a position to turn the idea mentioned in the beginning of this section into

LEMMA 2.3. *Let $\beta(F) > 0$ be the Hölder exponent of $F$ locally around $\mathcal{M}$. There exists a constant $C \equiv C(d, F)$ with the following property. Let $h \in W^{2,d}(\mathbb{R}^d)$ be a strong solution of the equation $(L_\varepsilon - \lambda)h = 0$ in $B(x, \varepsilon^{1/(1+\beta)})$, $\beta \in (\beta(F)/2, \beta(F))$, where $x \in \mathcal{M}$ and $0 \leq \lambda \leq \varepsilon$. Then there exists $\tilde{x} \in B(x, \varepsilon^{1/(1+\beta)})$ such that $h$ does not change sign in $B(\tilde{x}, \varepsilon)$.*

Let $\Gamma$, $\Sigma \subset \Gamma$, be regular domains. Let $g \in L^\infty_{\text{loc}}(\mathbb{R}^d)$ and let $h$ be a nonnegative, strong (i.e., twice weakly differentiable) solution of the equation $(L_\varepsilon - \lambda) \times h = g$ in $\Gamma$. Assume that there are $0 < r < 1/2$ and $B(x, 2\sqrt{\varepsilon}) \subset \Sigma$ such that $(1-r)\sup_\Sigma h \leq \sup_{B(x,\varepsilon)} h$. Then for all $0 \leq \lambda < \lambda(\Gamma \setminus B(x,\varepsilon))$ there is $C \equiv C(d, F|B(x, 2\sqrt{\varepsilon}))$ such that

$$\text{(2.11)} \quad \text{osc}_\Sigma h \leq \left(4r + C\varepsilon^{d/2}\lambda + 4 \sup_{\Sigma \setminus B(x,\varepsilon)} h^\lambda_{\mathbb{R}^d \setminus \Gamma, B(x,\varepsilon)}\right) \sup_{B(x,\varepsilon)} h + C\varepsilon^{d/2} \sup_{B(x, 2\sqrt{\varepsilon})} |g|.$$

Having established positivity of eigenfunctions in vicinities of the local minima of $F$, we may use strong pointwise regularity such as the local (boundary) maximum principle Theorem 9.20 in [24] (Theorem 9.26 in [24]), the Harnack inequality Theorem 8.20 or 9.22 in [24] and the (boundary) Hölder estimates Corollary 9.24 in [24] (Corollary 9.28 in [24]).

For later purpose also let us define for $x \in \mathbb{R}^d$

$$\text{(2.12)} \quad \delta(x) \equiv \delta_{F,\varepsilon}(x) \equiv \sup\left\{\delta > 0 \,\Big|\, 8\delta \sup_{B(x, 8\varepsilon\delta)} |\nabla F| < 1\right\}.$$

Clearly, $\varepsilon\delta$ only depends on $|\nabla F|/\varepsilon$ and $8\varepsilon\delta(x) \sup_{B(x, 4\varepsilon\delta(x))} |\nabla F|/\varepsilon = 1$. Combination of Harnack's and Hölder's principles gives:

THEOREM 2.4. *Assume that $\nabla F$ is locally Hölder continuous. Fix $0 < \rho \leq \varepsilon$. Let $0 \leq h \in W^{2,d}(\mathbb{R}^d)$ be a strong, nonnegative solution of the equation $(L_\varepsilon - \lambda)h = 0$ in $B(x, 2\sqrt{\rho})$, where $x \in \mathcal{M}$ and $0 \leq \lambda \leq 1$. Then there exists $C = C(F|B(x, 2\sqrt{\rho}))$ and $\alpha = \alpha(F|B(x, 2\sqrt{\rho})) > 0$ such that for all $0 < r < \sqrt{\rho}$*

$$\text{(2.13)} \quad \text{osc}_{B(x,r)} h \leq C(r/\sqrt{\rho})^\alpha \inf_{B(x,r)} h.$$

*Assume that $0 \leq h \in W^{2,d}(\mathbb{R}^d)$ is a strong, nonnegative solution of the equation $(L_\varepsilon - \lambda)h = f$ in $\Sigma$, where $\Sigma$ is an open, regular set, $0 \leq \lambda \leq 1$ and $f$ is in $L^d(\Sigma)$. There are constants $C = C(d)$ and $\alpha = \alpha(d) > 0$ such that for all $x \in \Sigma$ and all $0 < \rho < \varepsilon\delta(x)$ satisfying $B(x, 4\rho) \subset \Sigma$ and all $0 < r < \rho$,*

$$\text{(2.14)} \quad \text{osc}_{B(x,r)} h \leq C(r/\rho)^\alpha \left(\inf_{B(x,r)} h + \|f\|_{L^d(B(x,r) \cap \Sigma)}\right).$$

*For $x \in \partial\Sigma$ let $V_x$ be the exterior cone at $x$. We still have for some constant $C \equiv C(d, V_x)$ and $\alpha \equiv \alpha(d, V_x) > 0$ and all $0 < r < \rho < \varepsilon\delta(x)$*

$$\text{(2.15)} \quad \begin{aligned}\text{osc}_{B(x,r) \cap \Sigma} h &\leq C(r/\rho)^\alpha (\text{osc}_{B(x,\rho) \cap \Sigma} h + \|f - \lambda h\|_{L^d(B(x,r) \cap \Sigma)}) \\ &\quad + C\, \text{osc}_{B(x, \sqrt{r\rho}) \cap \partial\Sigma} h,\end{aligned}$$

*where $\text{osc}_{B(x,r) \cap \partial\Sigma} h \equiv \limsup_{y \to B(x,r) \cap \partial\Sigma} h - \liminf_{y \to B(x,r) \cap \partial\Sigma} h$.*

We also need the boundary Harnack inequality, which is a consequence of Theorem 8.0.1 in [36].

THEOREM 2.5. *Assume that $\nabla F$ is locally Hölder continuous. Let $\Sigma$ be an open set with uniformly Lipschitz continuous boundary. There exist $C \equiv C(d)$, $\rho \equiv \rho(d) > 0$ and a function $R : \partial \Sigma \to (0, \infty)$, $R \leq \delta$, with the following properties. Fix $z \in \partial \Sigma$ and write $\partial \Sigma \cap B_z = \operatorname{graph} \chi$ for some ball $B_z$ around $z$ and some function $\chi$. Fix $0 \leq r \leq \rho$ and let $0 < u, v \in W^{2,d}(\Sigma)$ be positive solutions of $L_\varepsilon h = 0$ in $\Sigma \cap B_z \cap B(z, 8r\varepsilon R(z))$ and $h = 0$ on $\partial \Sigma \cap B_z \cap B(z, 8r\varepsilon R(z))$. Then*

$$(2.16) \qquad \frac{u(x)}{v(x)} \leq C \frac{u(y)}{v(y)}, \qquad x, y \in \Sigma \cap B_z \cap B(z, r\varepsilon R(z)).$$

PROOF. Denote by $1/\gamma(z)$ the best Lipschitz constant of $\chi$ at $z$ in $B(z, 8\varepsilon)$ and let $1/\beta(z)$ be the best Hölder constant of $\nabla F$ in $B(z, 8\varepsilon)$. Define $R(z) \equiv \min(\beta(z), \gamma(z), \delta(z))$, where $\delta(z)$ is given in (2.12). Let us introduce the function $\tilde{u}(\tilde{x}) \equiv u(x)$, $\tilde{x} \equiv (x - z)/(\varepsilon R(z))$, and likewise $\tilde{v}$. Furthermore, let $\tilde{L} \equiv -\Delta + \tilde{b} \cdot \nabla$, where $\tilde{b}(\tilde{x}) \equiv \tilde{b}_{\varepsilon R(z)}(\tilde{x}) \equiv R(z) \nabla F(x)$. Fix $r > 0$ and let $u$ and $v$ be $L_\varepsilon$-harmonic in $B(z, 8r\varepsilon R(z)) \cap B_z \cap \Sigma$, vanishing identically on $B(z, 8r\varepsilon R(z)) \cap B_z \cap \partial \Sigma$. We then compute $\tilde{L}\tilde{u} = \tilde{L}\tilde{v} = 0$ in $B(0, 8r) \cap \tilde{B}_z \cap \tilde{\Sigma}$, where $\tilde{\Sigma} \equiv \tilde{\Sigma}_{\varepsilon R(z)} \equiv \{\tilde{x} | x \in \Sigma\}$ and likewise $\tilde{B}_z$, and clearly $\tilde{u} = \tilde{v} = 0$ on $B(0, 8r) \cap \tilde{B}_z \cap \partial \tilde{\Sigma}$. Note that by definition of $R(z)$ under this transformation, the best Lipschitz constant of $\tilde{\chi}(\tilde{x}) \equiv \chi(x)$ at $z$ and the best Hölder constant of $\tilde{b}$ in $B(0, 8)$ are bounded by 1. Moreover, the supremum norm of $\tilde{b}$ in $B(0, 8)$ is that of $\nabla F$ in $B(z, 8\varepsilon R(z))$ and hence is bounded by 1 by definition of $\delta(z)$. The boundary Harnack principle Theorem 8.0.1 in [36] applied to $D \equiv B(0, 8) \cap \tilde{B}_z \cap \tilde{\Sigma}$ gives the existence of $C \equiv C(d)$ and $\rho \equiv \rho(d) > 0$ such that $\tilde{u}(x)/\tilde{v}(x) \leq C\tilde{u}(y)/\tilde{v}(y)$ for all $0 < r \leq \rho$ and $\tilde{x}, \tilde{y} \in B(0, r) \cap \tilde{B}_z \cap \tilde{\Sigma}$. □

On several occasions we shall meet the following obvious representation formula. The solution $h$ of the Poisson–Dirichlet problem (1.4) for an open, connected, regular set $\Gamma$ in a relatively compact, open, connected, regular set $\Sigma \subset\subset \Gamma$ is given by

$$(2.17) \qquad h(x) = G_\Sigma^\lambda g(x) + H_\Sigma^\lambda h(x), \qquad x \in \Sigma,$$

where $H_\Sigma^\lambda f$ is the $(L_\varepsilon - \lambda)$-harmonic extension of $f$ to $\Sigma$ and where the resolvent $G_\Sigma^\lambda$ is defined in (2.4). Several times in the sequel we shall use the following obvious consequence of (2.17) and the weak maximum principle:

$$(2.18) \quad \sup_K |h| \leq s_K^0(\partial \Sigma) \sup_\Sigma (\lambda |h| + |g|) + \sup_K H_\Sigma^0 |h|, \qquad K \subset \Sigma \subset \Gamma.$$

Let $G_\Sigma^\lambda(x,y)e^{F(y)/\varepsilon}$ be the (symmetric) kernel of $G_\Sigma^\lambda$ in $L^2(\Sigma^2, e^{-(F(x)+F(y))/\varepsilon}\,dx\,dy)$. It is easy to see that $(\Lambda_\varepsilon - e^{-F/\varepsilon}\lambda)G_\Sigma^\lambda e^{F/\varepsilon}f = f$ weakly for all $f \in L^2(\Sigma)$. Since $G_\Sigma^\lambda e^{F/\varepsilon}f(x) = \int_\Sigma G_\Sigma^\lambda(x,y)f(y)\,dy$ by definition, and since $L^2(B(y,r))$, $y \in \Sigma$, $r > 0$, is separable, $G_\Sigma^\lambda(\cdot,y)$ is $(\Lambda_\varepsilon - e^{-F/\varepsilon}\lambda)$-harmonic in $\Sigma\setminus\overline{B(y,r)}$ and almost all $z \in B(y,r) \cap \Sigma$. Theorem 8.8 in [24] and Theorem 9.19 in [24] imply that $G_\Sigma^\lambda(\cdot,z)$ is $\mathcal{C}^{2,\alpha}(\Sigma\setminus\overline{B(y,r)})$ for those $z$. Symmetry of $G_\Sigma^\lambda(x,z)e^{F(z)/\varepsilon}$ implies the same assertion for all $z \in B(y,r)$. Therefore, $G_\Sigma^\lambda(x,y)e^{F(y)/\varepsilon}$ is in $\mathcal{C}^{2,\alpha}(\Sigma^2\setminus D)$, where $D \equiv \{(x,x)|x \in \mathbb{R}^d\}$. We recall from, for example, [4] that the $(L_\varepsilon - \lambda)$-harmonic extension $H_\Sigma^\lambda f$, $\Sigma$ regular, open and connected, of a function $f \in L^\infty(\partial\Sigma)$ is given by

$$
\begin{aligned}
(2.19) \quad h(x) = H_\Sigma^\lambda f(x) &= -\varepsilon \int_{\partial\Sigma} f(y)\,\partial_{n(y)} G_\Sigma^\lambda(y,x) e^{(F(x)-F(y))/\varepsilon}\,d\sigma(y) \\
&= -\varepsilon \int_{\partial\Sigma} f(y)\,\partial_{n(y)} G_\Sigma^\lambda(x,y)\,d\sigma(y),
\end{aligned}
$$

where $n(y)$ denotes the outer unit normal at $y \in \partial\Sigma$ and the normal derivative is taken from the inside of $\Sigma$. Here we have used that $e^{-F(x)/\varepsilon}G_\Sigma^\lambda(x,y)$ is symmetric in $x$ and $y$, that the normal derivative at $\partial\Sigma$ exists and that $G_\Sigma^\lambda(x,y)$ vanishes on the boundary.

As already pointed out, the problem is that a priori we cannot apply Theorem 2.4 to an eigenfunction $\phi$. However, by combination of Theorem 2.4 with the Poisson representation formula (2.19) we still can control the regularity of $\phi$.

PROOF OF LEMMA 2.3. By standard comparison arguments with the ordinary Laplace operator in $B \equiv B(x,R)$ as can be found, for example, in the proof of Theorem 2.1(i) in [36] one finds $\delta = \delta(F|B) > 0$ and $C(F|B)$ such that

$$
(2.20) \quad s_B^0 \leq C\varepsilon^{(1-\beta)/(1+\beta)} \quad \text{for } R \equiv \delta\varepsilon^{1/(1+\beta)},
$$

where $s_B^0$ is defined in (2.3) and where $\beta \equiv \beta(F) > 0$ is smaller than or equal to the optimal Hölder exponent of $F$ locally around $x$. For the convenience of the reader we shall formulate the details of the proof in our situation. Define $v_R(y) \equiv (R^2 - |y-x|^2)/(2d\varepsilon)$ for $|y-x| \leq R$. We compute $-\varepsilon\Delta v_R = 1$ for $|y-x| \leq R$. Since $|\nabla F(y)\cdot\nabla v_R(y)| \leq \sup_B |\nabla F||y-x|/(d\varepsilon) \leq \delta\sup_B |\nabla F|/(d\varepsilon^{1-1/(1+\beta)})$ and since $x \in \mathcal{M}$, it follows that $L_\varepsilon v_R \geq 1/2$ for $\delta \equiv \sup\{r \in (0,1)\,|\,r\sup_{B(x,r\varepsilon^{1/(1+\beta)})}|\nabla F| \leq d\varepsilon^{1-1/(1+\beta)}/2\} > 0$ and $|y-x| \leq R$ so that $v_R(y) \geq (1/2)w_B(y)$ in $B$. Recall the notion of the principal eigenvalue $\lambda(\Sigma)$ of the Dirichlet operator $L_\varepsilon^\Sigma$ introduced in (1.13). For the purpose of (5.24) and (5.32) we note that on the other hand the same arguments show $v_R(y) \leq (3/2)w_B(y)$ in $B$ and therefore for some $C \equiv C(d,F)$ and all

$\delta \in (1/8, 1)$

(2.21)
$$w_B(y) = e^{\mathcal{O}(1)}\varepsilon$$
for $R \equiv \delta\varepsilon$ and $y \in B(x, R(1 - 1/100))$ and $\lambda(B) \geq 1/(C\varepsilon)$,

where the last inequality is a consequence of (2.20) and (2.5).

Since the uncertainty principle (2.5) tells us $\lambda(B) \geq 1/s_B^0$, the condition on $\lambda$ ensures that $G_B^\lambda$ exists and that $h$ satisfies (2.19). Choose a ball $\tilde{B} \subset B$ of radius $0 < \rho < \varepsilon$ such that $\sup_B |h| = \sup_{\tilde{B}} |h|$. Since $-\partial_{n(y)} G_B^\lambda(x, y)$ is a positive strong solution for every $y \in \partial\tilde{B}$, we may apply (2.15) and obtain for some $C \equiv C(d, F|B)$, $\beta = \beta(F|B)$ and all $y_0, y_1, y_2 \in \tilde{B}$,

(2.22)
$$\begin{aligned}|h(y_1) - h(y_2)| &\leq \varepsilon \int_{\partial B} |h(z)||\partial_{n(z)} G_B^\lambda(y_1, z) - \partial_{n(z)} G_B^\lambda(y_2, z)| \, d\sigma(z) \\ &\leq C(\rho/\varepsilon^{1/(1+\beta)})^\alpha \sup_{\partial B} |h| \varepsilon \int_{\partial B} -\partial_{n(z)} G_B^\lambda(y_0, z) \, d\sigma(z) \\ &\leq C(\rho/\varepsilon^{1/(1+\beta)})^\alpha \sup_{\tilde{B}} |h| \sup_B h_B^\lambda,\end{aligned}$$

where $h_B^\lambda \equiv H_B^\lambda \mathbb{1}_{\partial B}$. Applying (2.17) to $h_B^\lambda$ and $\Sigma \equiv B$, we obtain

(2.23)
$$\sup_B h_B^\lambda \leq \lambda s_B^0 \sup_B h_B^\lambda + 1.$$

Combination of (2.23) with (2.22) implies

(2.24)
$$|h(y_1) - h(y_2)| \leq \frac{C(\rho/\varepsilon^{1/(1+\beta)})^\alpha \sup_{\tilde{B}} |h|}{1 - \lambda s_B^0}.$$

The bounds on $\lambda$ and $s_B^0$ show that the denominator can be absorbed in the constant. We thus have proven for some $C$ and $\alpha > 0$ only depending on $F|B$ that for all $y_1, y_2 \in \tilde{B}$,

(2.25)
$$\sup_{y_1, y_2 \in \tilde{B}} |h(y_1) - h(y_2)| \leq C(\rho/\varepsilon^{1/(1+\beta)})^\alpha \sup_{\tilde{B}} |h|.$$

Now let us assume that there is $y \in \tilde{B}$ such that $h(y) = 0$. We apply (2.18) for $\Sigma \equiv \tilde{B}$ and deduce from (2.25), using the condition on $\tilde{B}$ and choosing $\rho \equiv \varepsilon$,

(2.26)
$$\sup_{\tilde{B}} |h| \leq \lambda s_B^0 \sup_{\tilde{B}} |h| + C\varepsilon^{\alpha\beta/(1+\alpha)} \sup_{\tilde{B}} |h|.$$

It follows that $h$ vanishes identically in $\tilde{B}$ for small $\varepsilon > 0$ and hence by analytic continuation everywhere in $\mathbb{R}^d$.

For the proof of (2.11) let $c_n \equiv (7C\varepsilon^\alpha)^n (C\varepsilon^\alpha + C|B|(2\lambda + \sup_B |g|/\sup_{\tilde{B}} h))$, where $B \equiv B(x, 2\varepsilon^{1/(1+\beta)})$. We claim the existence of $C \equiv C(d, F|B)$, $\alpha \equiv$

$\alpha(d, F|B) > 0$, such that for all $n$ the inequality $c_{n-1} > M \equiv \max(r, C|B|(2\lambda + \sup_B |g|/\sup_{\tilde{B}} h), \sup_{\Sigma \setminus \tilde{B}} h^\lambda_{\mathbb{R}^d \setminus \Gamma, \tilde{B}})$, $\tilde{B} \equiv B(x, \varepsilon)$, implies

$$\operatorname{osc}_{\tilde{B}} h \leq c_n \sup_{\tilde{B}} h. \tag{2.27}$$

For $n \equiv 1$ this is nothing more than (2.13). Assume (2.27) for some $n \geq 1$. It follows from $h \geq 0$ and (2.17) applied to $h$ for $K \equiv \Sigma \setminus \tilde{B}$ in $\Sigma \equiv \Gamma \setminus \tilde{B}$—in slight abuse of notation—that

$$\begin{aligned}\inf_\Sigma h &\geq 0 + (1 - c_n)\left(1 - \sup_{\Sigma \setminus \tilde{B}} h^\lambda_{\mathbb{R}^d \setminus \Gamma, \tilde{B}}\right) \sup_{\tilde{B}} h \\ &\geq (1 - c_n)(1 - r)\left(1 - \sup_{\Sigma \setminus \tilde{B}} h^\lambda_{\mathbb{R}^d \setminus \Gamma, \tilde{B}}\right) \sup_\Sigma h,\end{aligned} \tag{2.28}$$

where we use the convention that $h^\lambda_{\mathbb{R}^d \setminus \Gamma, \tilde{B}} \equiv 0$ in $\tilde{B}$. Thus

$$\begin{aligned}\operatorname{osc}_\Sigma h &\leq \left(c_n + r + \sup_{\Sigma \setminus \tilde{B}} h^\lambda_{\mathbb{R}^d \setminus \Gamma, \tilde{B}}\right) \sup_\Sigma h \\ &\leq \left(\left(c_n + r + \sup_{\Sigma \setminus \tilde{B}} h^\lambda_{\mathbb{R}^d \setminus \Gamma, \tilde{B}}\right)/(1 - r)\right) \sup_{\tilde{B}} h \\ &\leq 2\left(c_n + r + \sup_{\Sigma \setminus \tilde{B}} h^\lambda_{\mathbb{R}^d \setminus \Gamma, \tilde{B}}\right) \sup_{\tilde{B}} h.\end{aligned} \tag{2.29}$$

From (2.13) again we hence obtain

$$\begin{aligned}\operatorname{osc}_{\tilde{B}} h &\leq C\varepsilon^\alpha 2\left(c_n + r + \sup_{\Sigma \setminus \tilde{B}} h^\lambda_{\mathbb{R}^d \setminus \Gamma, \tilde{B}}\right) \\ &\quad + C|B|\left(\lambda/(1-r) + \sup_B |g|/\sup_{\tilde{B}} h\right) \sup_{\tilde{B}} h \\ &\leq c_{n+1} \sup_{\tilde{B}} h\end{aligned} \tag{2.30}$$

since $c_n > M$. Choosing $n$ maximal in (2.27), from (2.29) we obtain the estima- te since $c_n \leq M$. □

2.3. *A priori bounds on conditioned, expected exit times from bounded sets.* In this section we prove an estimate on the suprema $s^0_K(A, B)$ for regular, closed sets $A$ and $B$ with bounded complement of their union.

For the sake of convenience we set

$$T_J \equiv \max_{y \in \mathcal{M} \setminus J} T_{x, J}, \qquad J \subset \mathcal{M}, J \neq \mathcal{M}, \tag{2.31}$$

where the time scale $T_{x,J}$ is defined in (1.9). In the case $J \equiv \mathcal{M}$ we use the convention that $T_\mathcal{M} \equiv 1/\varepsilon^{d-1}$. For every finite set of points $I \subset \mathbb{R}^d$ such that

$\min_{x,y \in I, x \neq y} \text{dist}(x,y) > 2\delta$ we set $B_I \equiv B_I(\varepsilon/4)$, where

$$B_I(\delta) \equiv \bigcup_{x \in I} B(x, \delta). \tag{2.32}$$

We then have:

LEMMA 2.6. *Fix disjoint, regular, nonempty, closed sets $A, B \subset \mathbb{R}^d$ such that $B_I \subset A$ and $B_J \subset B$, where $I \equiv \mathcal{M} \cap A$ and $J \equiv \mathcal{M} \cap B$. There are $N \equiv N(d)$ and $C = C(d)$ such that*

$$s_K^0(A, B) \leq C\varepsilon^{-N}(T_{I \cup J} + |\mathbb{R}^d \backslash A \backslash B|). \tag{2.33}$$

We start with the following bound on the Green function $G_\Sigma^0(x,y)$ defined in (2.4).

LEMMA 2.7. *For all regular, open, bounded sets $\Gamma$ there exists $C = C(d)$ such that for all $x, y \in \Gamma$ and all $0 < \rho < \delta(y)$ satisfying $|x - y| > \rho\varepsilon$ and $\text{dist}(x \cup y, \partial \Gamma) \geq 4\rho\varepsilon$,*

$$G_\Gamma^0(x, y) \leq \frac{C h_{B(y,\rho\varepsilon), \Gamma^c}(x) e^{-F(y)/\varepsilon}}{\text{cap}_{B(y,\rho\varepsilon)}(\Gamma^c)}. \tag{2.34}$$

*For $|x - y| < \rho\varepsilon$, $0 < \rho \leq \delta(x)\varepsilon$ and $\text{dist}(x, \partial\Gamma) > 4\rho\varepsilon$ we have*

$$G_\Gamma^0(x, y) \leq \frac{C}{\varepsilon} G^\Delta(|x-y|) + \frac{C e^{-F(x)/\varepsilon}}{\text{cap}_{B(x,\rho\varepsilon)}(\Gamma^c)}, \tag{2.35}$$

*where $G^\Delta$ is the Green function of the Laplace operator in $\mathbb{R}^d$.*

PROOF. Let $h_y \equiv h^0_{B(y,\rho\varepsilon),\Gamma^c}$. The second Green formula as, for example, in (2.8) in [4] for $\Gamma \equiv \Gamma \backslash \overline{B}(y,\rho\varepsilon) \backslash \overline{B}(x,r)$ and for $\Gamma \equiv B(y,\rho\varepsilon)$ shows for all $0 < r < |x - y| - \rho\varepsilon$,

$$\begin{aligned}
\varepsilon &\int_{\partial B(y,\rho\varepsilon)} e^{-F/\varepsilon} G_\Gamma(\cdot, x) \, \partial_n h_y \, d\sigma \\
&= -\varepsilon \int_{\partial B(x,r)} h_y \, \partial_n G_\Gamma(\cdot, x) + \varepsilon \int_{\partial B(y,\rho\varepsilon)} e^{-F/\varepsilon} \, \partial_n G_\Gamma(\cdot, x) \, d\sigma \\
&= -e^{\mathcal{O}(1)(r/(|x-y|-\rho\varepsilon))^\alpha} h_y(x) e^{-F(x)/\varepsilon} \varepsilon \int_{\partial \Gamma} \partial_n G_\Gamma(x, \cdot) \, d\sigma \\
&= e^{-F(x)/\varepsilon} h_y(x),
\end{aligned} \tag{2.36}$$

where $\partial_n$ is the normal derivative taken from the interior with respect to the outer unit normal at the boundary. The last equation uses (2.19) and the fact that $r > 0$ can be chosen arbitrarily small. Invoking the Harnack inequality Corollary 9.25 in [24] on the left-hand side of (2.36), we thus have for some $C = C(d)$

$$G_\Gamma(y,x)\,\mathrm{cap}_{B(y,\rho\varepsilon)}(\Gamma^c) \leq C e^{-F(x)/\varepsilon} h_y(x). \tag{2.37}$$

Equation (2.34) now follows from the symmetry of $G_\Gamma(y,x)e^{F(x)/\varepsilon}$ in $x$ and $y$.

For the proof of (2.35) we first observe that considering $y \equiv x$ in (2.36) a similar calculation gives for all $0 < \rho \leq 1$

$$\varepsilon \int_{\partial B(x,\rho\varepsilon)} e^{-F/\varepsilon} G_\Gamma(\cdot,x)\,\partial_n h_x = e^{-F(x)/\varepsilon}. \tag{2.38}$$

Analogously to (2.37), we find $C$ independent of $\varepsilon$ such that for all $y \in \partial B(x,\rho\varepsilon)$,

$$G_\Gamma(y,x) \leq C \frac{e^{-F(x)/\varepsilon}}{\mathrm{cap}_{B(x,\rho\varepsilon)}(\Gamma)}. \tag{2.39}$$

For, choose a sequence of points $y_0 = y, \ldots, y_k = z \in \partial B(x,2\rho\varepsilon)$ such that $\rho\varepsilon/100 < |y_i - y_{i+1}| < \rho\varepsilon/3$ and balls $B_i$ of radii $\rho\varepsilon/3$ such that $y_{i-1}, y_i \in B_i$. Applying the Harnack inequality to each ball, we derive $G_\Gamma(y_i,x)/G_\Gamma(y_{i+1},x) \leq C$ for some $C \equiv C(d)$ and $0 < \rho \leq \delta(x)$. Since the arclength of a ball depends linearly on the distance, we get that $k$ is bounded independent of $\varepsilon$ and $y,z \in \partial B(x,2\rho\varepsilon)$ and thus

$$G_\Gamma(y,x)/G_\Gamma(z,x) \leq C^k, \tag{2.40}$$

from which (2.39) follows. Assume first that $\Gamma = B(x,\rho\varepsilon)$. Invoking the Dirichlet principle for the capacity, we derive for $0 \leq r \leq \rho$

$$\mathrm{cap}_{B(x,r\varepsilon)}(B(x,\rho\varepsilon)^c) \leq \varepsilon C e^{-F(x)/\varepsilon}\,\mathrm{cap}^\Delta_{B(x,r\varepsilon)}(B(x,\rho\varepsilon)^c), \tag{2.41}$$

where $\mathrm{cap}^\Delta$ denotes the capacity with respect to the Laplace operator. Since $G^\Delta$ is rotationally invariant, it follows that

$$\mathrm{cap}^\Delta_{B(x,r\varepsilon)}(B(x,\rho\varepsilon)^c) = 1/(G^\Delta(\rho\varepsilon) - G^\Delta(r\varepsilon)) \geq 1/G^\Delta(\rho\varepsilon). \tag{2.42}$$

Combination of (2.41) and (2.39) shows

$$G_{B(x,\rho\varepsilon)}(y,x) \leq \frac{C}{\varepsilon} G^\Delta(|x-y|), \qquad |x-y| < \rho\varepsilon. \tag{2.43}$$

To obtain the full estimate we note that the function $G_\Omega(\cdot,x) - G_{B(x,\rho\varepsilon)}(\cdot,x) - h$ is a weakly $L_\varepsilon$-harmonic function in $B(x,\rho\varepsilon)$ and equals zero on $\partial B(x,\rho\varepsilon)$, where $h$ is the solution of the Dirichlet problem in $B(x,\rho\varepsilon)$ with boundary values $G_{B(x,\rho\varepsilon)}(\cdot,x)$. By (2.39) and (2.43) we thus have proven for $|x-y| < \rho\varepsilon$

$$G_\Gamma(y,x) \leq \frac{C}{\varepsilon} G^\Delta(|x-y|) + \frac{C e^{-F(x)/\varepsilon}}{\mathrm{cap}_{B(x,\rho\varepsilon)}(\Gamma)} \tag{2.44}$$

which gives (2.35). □

For later purpose we notice that the definition of $\delta$ in (2.12) implies

(2.45) $$\inf_{\Sigma} F - \varepsilon \leq \inf_{\Sigma^{\varepsilon\delta}} F \leq \sup_{\Sigma^{\varepsilon\delta}} F \leq \sup_{\Sigma} F + \varepsilon.$$

Indeed, fix arbitrary $x \in \Sigma$ and let $y \in \partial B(x, \varepsilon\delta(x))$. We then obtain, using (2.12),

(2.46) $$F(y) - F(x) = \int_0^1 \nabla F((1-t)x + ty) \cdot (y - x)\, dt \leq |x - y|/\delta(x) = \varepsilon$$

and $F(y) - F(x) \geq -\varepsilon$ by replacing the roles of $x$ and $y$.

PROOF OF LEMMA 2.6. Applying (2.16), respectively (2.15), to $w_{A,B}^0$, respectively $h_{A,B}^0$, with the obvious choice $z \in \partial B$, respectively $z \in \partial A$, we may assume that $x \in \mathbb{R}^d \setminus \tilde{A} \setminus \tilde{B}$, where $\tilde{A} \equiv A \cup (\partial A)^{R\varepsilon}$, $\tilde{B} \equiv B \cup (\partial B)^{R\varepsilon}$ and where $R: \partial A \cup \partial B \to (0, \infty)$ is as in Theorem 2.5. Let $\varepsilon r: \mathbb{R}^d \setminus A \setminus \tilde{B} \to (0, \infty)$ be the maximum of $\varepsilon\delta$ and the distance from $B$. We may assume that $R$ is bounded by 1. We now may write, using (2.35) for $\Gamma \equiv \mathbb{R}^d \setminus A \setminus B$ and all $x \in \mathbb{R}^d \setminus \tilde{A} \setminus \tilde{B}$,

(2.47) $$\begin{aligned}
\frac{w_{A,B}^0(x)}{h_{A,B}^0(x)} &= \int_{\mathbb{R}^d \setminus A \setminus B} \frac{G_{\mathbb{R}^d \setminus A \setminus B}^0(x,y) h_{A,B}^0(y)}{h_{A,B}^0(x)}\, dy \\
&\leq C \int_{\tilde{B} \setminus B} \frac{G_{\mathbb{R}^d \setminus A \setminus B}^0(x,y) h_{A,B}^0(y)}{h_{A,B}^0(x)}\, dy \\
&\quad + \frac{C}{\varepsilon} + \frac{Ce^{-F(x)/\varepsilon}}{\mathrm{cap}_{B(x,\varepsilon r(x))}(A \cup B)} |B(x,\varepsilon)| \\
&\quad + C \int_{\substack{\mathbb{R}^d \setminus A \setminus \tilde{B} \\ y:\, |y-x|>\varepsilon}} \frac{h_{B(y,r(y)\varepsilon),A \cup B}^0(x) h_{A,B}^0(y)}{h_{A,B}^0(x) \mathrm{cap}_{B(y,r(y)\varepsilon)}(A \cup B)} e^{-F(y)/\varepsilon}\, dy.
\end{aligned}$$

Using the Harnack inequality for $h_{A,B}^0$ on $B(y, \varepsilon r(y))$, we observe for $x \in (A \cup B \cup B(y, \varepsilon r(y)))^c$ that

(2.48) $$\begin{aligned}
h_{B(y,\varepsilon r(y)),A \cup B}^0(x) h_{A,B}^0(y) &\leq C H_{(A \cup B \cup B(y,\varepsilon r(y)))^c}^0 \mathbb{1}_{\partial B(y,\varepsilon r(y))} h_{A,B}(x) \\
&\leq C h_{A,B}^0(x).
\end{aligned}$$

By (2.16) for $u \equiv G_{\mathbb{R}^d \setminus A \setminus B}^0(x, \cdot)$ and $v \equiv h_{A,B}^0$ and by (2.34) once more we have for all $y \equiv z + tR(z)\varepsilon n(z)$, $z \in \partial B$, $0 < t \leq 1$, the bound

(2.49) $$\begin{aligned}
G_{\mathbb{R}^d \setminus A \setminus B}^0&(x,y) h_{A,B}^0(y) \\
&\leq C G_{\mathbb{R}^d \setminus A \setminus B}^0(x, z + R(z)\varepsilon n(z)) h_{A,B}^0(z + R(z)\varepsilon n(z)) \\
&\leq C \frac{h_{B(z+R(z)\varepsilon n(z), R(z)\varepsilon/4), A \cup B}^0(x) h_{A,B}^0(z + R(z)\varepsilon n(z))}{\mathrm{cap}_{B(z+R(z)\varepsilon n(z), R(z)\varepsilon/4)}(A \cup B)^c},
\end{aligned}$$

where $n$ is the outer unit normal vector field at the boundary of $B$. We hence may apply the Harnack inequality again to the right-hand side of (2.49), proving

$$(2.50) \quad G^0_{\mathbb{R}^d \setminus A \setminus B}(x,y) h^0_{A,B}(y) \leq C \frac{h^0_{A,B}(x)}{\text{cap}_{B(z+R(z)\varepsilon n(z), R(z)\varepsilon/4)}(A \cup B)^c}.$$

Invoking Proposition 4.7 in [4], we have for all $0 < \rho < 1$

$$(2.51) \quad \text{cap}_{B(y,\rho\varepsilon)}(A \cup B) \geq e^{-\hat{F}(y, A \cup B)/\varepsilon} (\rho\varepsilon)^d / (C\rho\varepsilon).$$

Inserting (2.50) and (2.51) into the integrals on the right-hand side of (2.47), we thus may bound the right-hand side by $C$ times

$$(2.52) \quad \begin{aligned} &\frac{1}{\varepsilon} + \frac{1}{\varepsilon} e^{(\hat{F}(x, I \cup J) - F(x))/\varepsilon} \\ &+ |\{F = \hat{F}(\cdot, A \cup B)\} \setminus A \setminus B| \\ &+ |\{\text{dist}(\cdot, B) \leq R\varepsilon\}| \sup_{z \in \partial B, 0 < t < 1} e^{|F(z + tR(z)\varepsilon n(z)) - F(z + R(z)\varepsilon n(z))|/\varepsilon} \\ &+ \frac{1}{\varepsilon^{d-1}} \int_{\{F < \hat{F}(\cdot, A \cup B)\} \setminus A \setminus \tilde{B}} e^{(\hat{F}(y, I \cup J) - F(y))/\varepsilon} \, dy. \end{aligned}$$

Since $R(z) \leq \delta(z)$, (2.45) tells us that the supremum appearing in the fourth term is bounded by $\varepsilon$. We readily verify (2.33) by computation of a Laplace-type integral. $\square$

**3. Growth estimates at infinity.** Because of the strong drift of $-\nabla F$ toward the local minima of $F$, the influence of the values of a solution $\phi$ of (2.1) at infinity on its values in compact sets can be neglected. Technically, this will be achieved by weighted $L^2$-estimates near infinity in the spirit of Agmon and Helffer and Sjöstrand (see [1] and [26]) in combination with pointwise estimates based on the maximum principle in compact sets.

3.1. *Laplace transforms in compact sets.* The following lemma provides us good control on Laplace transforms $h^\lambda_{A,B}$ in compact sets away from its first pole $\lambda(\mathbb{R}^d \setminus A \setminus B)$ in terms of the maximal conditioned expected exit time from $\mathbb{R}^d \setminus A \setminus B$.

LEMMA 3.1. *Fix regular, closed, disjoint, nonempty sets $A$ and $B$ with bounded complement $\mathbb{R}^d \setminus A \setminus B$. Assume that $B_I \subset A$ and $B_J \subset B$, where $I \equiv \mathcal{M} \cap A$ and $J \equiv \mathcal{M} \cap B$. Assume that $0 \leq \lambda s^0(A,B) \leq 1/2$. Then for some $N \equiv N(d)$ and $C \equiv C(d)$ and all $x \notin A \cup B$,*

$$(3.1) \quad \frac{h^\lambda_{A,B}(x)}{h^0_{A,B}(x)} \leq 1 + \lambda C \varepsilon^{-N}(T_{I \cup J} + |\mathbb{R}^d \setminus A \setminus B|).$$

*Moreover,*

$$\text{(3.2)} \qquad \frac{w_{A,B}^\lambda(x)}{w_{A,B}^0(x)} \leq 1 + \lambda C \varepsilon^{-N}(T_{I \cup J} + |\mathbb{R}^d \backslash A \backslash B|).$$

PROOF. Equation (2.5) and the condition on $\lambda$ show that $G_\Sigma^\lambda$ and $H_\Sigma^\lambda$, where $\Sigma \equiv \{\alpha < F < \beta\}$, exist. The Harnack inequality, Theorem 8.20 in [24], the weak maximum principle and (2.17) applied to $h \equiv h_{A,B}^\lambda - h_{A,B}^0$ yield for all $x \in \Sigma \equiv \mathbb{R}^d \backslash A \backslash B$

$$\text{(3.3)} \qquad \begin{aligned} \frac{h_{A,B}^\lambda(x)}{h_{A,B}^0(x)} - 1 &= \lambda \frac{1}{h_{A,B}^0(x)} G_\Sigma^0 \left( h_{A,B}^0 \left( \frac{h_{A,B}^\lambda}{h_{A,B}^0} - 1 \right) \right)(x) + \lambda \frac{G_\Sigma^0 h_{A,B}^0(x)}{h_{A,B}^0(x)} \\ &\leq \lambda s^0(A, B) \sup_\Sigma \left( \frac{h_{A,B}^\lambda}{h_{A,B}^0} - 1 \right) + \lambda s^0(A, B). \end{aligned}$$

Taking the supremum on the left-hand side and assuming it is finite, we have proven

$$\text{(3.4)} \qquad \sup_\Sigma \frac{h_{A,B}^\lambda}{h_{A,B}^0} \leq 1 + 2\lambda s^0(A, B).$$

Simply by continuity at the boundary, the supremum stays finite near boundary points $x_0 \in \partial A$. Since $h_{A,B}^0$ takes its minimal value at zero by the Hopf maximum principle Theorem 3.2.5 in [36], it follows that

$$\text{(3.5)} \qquad \lim_{\Sigma \ni x \to x_0} \frac{h_{A,B}^\lambda(x) - 0}{h_{A,B}^0(x) - 0} = \frac{\partial_{n(x_0)} h_{A,B}^\lambda(x_0)}{\partial_{n(x_0)} h_{A,B}^0(x_0)} < \infty.$$

Equation (3.1) now follows from (2.33).

For the proof of (3.2) we apply (2.17) to $h \equiv w_{A,B}^\lambda - w_{A,B}^0$ and obtain for all $x \notin A \cup B$

$$\text{(3.6)} \qquad \begin{aligned} \frac{w_{A,B}^\lambda(x)}{w_{A,B}^0(x)} - 1 &= \frac{\lambda}{w_{A,B}^0(x)} G_{\mathbb{R}^d \backslash A \backslash B}^0 \left( h_{A,B}^0 \frac{w_{A,B}^0}{h_{A,B}^0} \frac{w_{A,B}^\lambda}{w_{A,B}^0} \right)(x) \\ &\quad + \frac{1}{w_{A,B}^0(x)} G_{\mathbb{R}^d \backslash A \backslash B}^0 \left( h_{A,B}^0 \left( \frac{h_{A,B}^\lambda}{h_{A,B}^0} - 1 \right) \right)(x) \\ &\leq \lambda s_{A,B}^0 \sup_{\mathbb{R}^d \backslash A \backslash B} \frac{w_{A,B}^\lambda}{w_{A,B}^0} + \sup_{\mathbb{R}^d \backslash A \backslash B} \frac{h_{A,B}^\lambda}{h_{A,B}^0} - 1. \end{aligned}$$

By the Hopf maximum principle we again may take the supremum over $\mathbb{R}^d \backslash A \backslash B$ in this inequality. The assertion thus follows from (3.2) and (2.33). □

3.2. *Weighted estimates.* Let $\tilde{F}$ be a $\mathcal{C}^\infty$-function on a regular domain $\Sigma$. Denoting by $L_\varepsilon^{\tilde{F}}$ the operator defined in (1.1), we have that the $e^{\tilde{F}/(2\varepsilon)}$-transform $H_\varepsilon^{\tilde{F}} \equiv e^{-\tilde{F}/(2\varepsilon)} L_\varepsilon^{\tilde{F}} e^{\tilde{F}/(2\varepsilon)}$ equals the Schrödinger operator

$$(3.7) \qquad H_\varepsilon^{\tilde{F}} = -\varepsilon \Delta + V_\varepsilon^{\tilde{F}}, \qquad V_\varepsilon^{\tilde{F}} \equiv |\nabla \tilde{F}|^2/(4\varepsilon) - \Delta \tilde{F}/2.$$

Fix $u \in \mathcal{C}^2(\Sigma) \cap \mathcal{C}^1(\overline{\Sigma})$. The well-known basic identity (see, e.g., Theorem 3.1.1 in [25])

$$(3.8) \quad \begin{aligned} &\varepsilon \int_\Sigma |\nabla e^{\varphi/\varepsilon} u|^2 \, dx + \int_\Sigma (V_\varepsilon^{\tilde{F}} - |\nabla \varphi|^2/\varepsilon)|e^{\varphi/\varepsilon} u|^2 \, dx \\ &\qquad = \frac{\varepsilon}{2} \int_{\partial \Sigma} \partial_n |e^{\varphi/\varepsilon} u|^2 \, d\sigma + \int_\Sigma e^{2\varphi/\varepsilon} u H_\varepsilon^{\tilde{F}} u \, dx \end{aligned}$$

for $L^2$-decay estimates is a consequence of Green's first formula and Gauss's divergence theorem. Equation (3.7) holds for all Lipschitz continuous functions $\varphi$ on $\overline{\Sigma}$. Fix $C_3 > C_2 > C_1$ and let us now assume that $\tilde{F}$ is close to $F$ in $\mathcal{C}^1(\Sigma)$ such that $\sup_{\{F<C_3\} \cap \Sigma} |\Delta \tilde{F}| \leq \text{ess} - \sup_{\{F<C_3\} \cap \Sigma} |\Delta F|$ and such that the conditions in Assumption 1.2 are also satisfied by $\tilde{F}$ for slightly modified constants $\tilde{C}_1 < C_2$, $\tilde{c}$ and $\tilde{\mu}_\varepsilon$ defined with respect to $\tilde{F}$ instead [e.g., let $\tilde{F}(x) \equiv \int \varphi_\delta(x-y) F(y) \, dy$ for $\delta > 0$ sufficiently small, where $\varphi_\delta$ is the density of the centered normal distribution with covariance matrix $(\delta \, \delta_{ij})_{i,j \leq d}$]. Being only interested in bounds on eigenfunctions in compact sets, we can bypass conditions like a uniform lower bound on $V_\varepsilon^{\tilde{F}}$. For, set $R_2 \equiv \sup\{|x| \, | \, \tilde{F}(x) < C_2\}$, let $L \geq R_2$ and assume that $B(0, R_2 + L) \subset \{\tilde{F} < C_3\}$. Let $\chi : \mathbb{R} \to [0,1]$ be a smooth, decreasing cut-off function with $\chi = 1$ on $(-\infty, 1]$, $\chi = 0$ on $[9, \infty)$ and $\chi(4) = 1/2$ and introduce nonnegative functions $J_1(x) \equiv (1 - J_2(x)^2)^{1/2} \equiv J(x) \equiv \chi((|x|^2 - (R_2)^2)/L^2)$. The IMS localization formula (see Theorem 3.2 in [8]) reads

$$(3.9) \quad \begin{aligned} H_\varepsilon^{\tilde{F}} &= J_1 H_\varepsilon^{\tilde{F}} J_1 + J_2 H_\varepsilon^{\tilde{F}} J_2 - |\nabla J_1|^2 - |\nabla J_2|^2 \\ &= J_1 H_\varepsilon^{\tilde{F}} J_1 + J_2 H_\varepsilon^{\tilde{F}} J_2 - \frac{|\nabla J|^2}{1 - J^2}. \end{aligned}$$

Since on the left-hand side of (3.8) there appears the quadratic form of $H_\varepsilon^{\tilde{F}}$ applied to $e^{\varphi/\varepsilon} u$ and since by Assumption 1.2 and monotonicity in volume of the principal eigenvalue for $L > R_2$, $(1280/15)^{1/2} \sup |\chi'|/(\tilde{\mu}_\varepsilon)^{1/2}$

$$(3.10) \qquad J_2 H_\varepsilon^{\tilde{F}} J_2 - \mathbb{1}_{\{J<1/2\}} |\nabla J|^2/(1-J^2) \geq \tilde{\mu}_\varepsilon \mathbb{1}_{\{J<1/2\}}/8$$

on $W_0^{1,2}(\{J \neq 1\})$, we have

$$(3.11) \quad \begin{aligned} &\varepsilon \int_\Sigma |\nabla J e^{\varphi/\varepsilon} u|^2 \, dx \\ &\quad + \int_\Sigma \left( J^2 V_\varepsilon^{\tilde{F}} + \frac{\tilde{\mu}_\varepsilon}{8} \mathbb{1}_{\{J<1/2\}} - \frac{1}{\varepsilon} |\nabla \varphi|^2 - \mathbb{1}_{\{J>1/2\}} \frac{|\nabla J|^2}{1-J^2} \right) |e^{\varphi/\varepsilon} u|^2 \, dx \\ &\quad \leq \frac{\varepsilon}{2} \int_{\partial \Sigma} \partial_n |e^{\varphi/\varepsilon} u|^2 \, d\sigma + \int_\Sigma e^{2\varphi/\varepsilon} u H_\varepsilon^{\tilde{F}} u \, dx. \end{aligned}$$

Choose $\Sigma \equiv \{\tilde{F} > \tilde{C}_1, J > 0\}$ and let $\tilde{\varphi}_J$ be the solution to the eiconal equation

$$(3.12) \quad |\nabla \tilde{\varphi}|^2 = J^2 |\nabla \tilde{F}|^2, \qquad \tilde{\varphi} = \tilde{C}_1 \text{ on } \{\tilde{F} = \tilde{C}_1\}, \qquad \nabla \tilde{\varphi}(x_0) = \nabla \tilde{F}(x_0),$$

for some $x_0 \in \{\tilde{F} = \tilde{C}_1\}$. By Theorem 5.5 in [23] and by local flattening of the set $\{\tilde{F} = \tilde{C}_1\}$, we can construct a unique, smooth solution defined on a neighborhood of this level set. In fact, we may assume that $\{J > 0, \tilde{F} > \tilde{C}_1\}$ is contained in the domain of $\tilde{\varphi}_J$. Moreover, as in Lemma 3.2.1 in [25] the solution can be identified with the Agmon distance corresponding to the potential $J^2 |\nabla \tilde{F}|^2$. More precisely, for a point $x \in \{J > 0, \tilde{F} > \tilde{C}_1\}$

$$(3.13) \quad \tilde{\varphi}_J(x) - \tilde{C}_1 = \rho(x) \equiv \inf_{\substack{c \colon [0,1] \to \operatorname{supp} J \\ c(0) = x, \tilde{F}(c(1)) = \tilde{C}_1}} \int_0^1 J(c(t)) |\nabla \tilde{F}(c(t))| |\dot{c}(t)| \, dt,$$

where the infimum is taken over all continuously differentiable curves. The proof of the upper bound $\tilde{\varphi}_J(x) - \tilde{C}_1 \leq \rho(x)$ is the same as that in Lemma 3.2.1 in [25] while the proof of the lower bound is a slight modification of the corresponding assertion. For convenience of the reader we shall give the details of this modification. Let $X_p \equiv (\nabla_\xi p, -\nabla_x p)$ be the Hamiltonian vector field corresponding to the Hamiltonian $p(x, \xi) \equiv |\xi|^2 - J(x)^2 |\nabla \tilde{F}(x)|^2$, $(x, \xi) \in \mathbb{R}^{2d}$. As in Proposition 5.4 in [23] we define $\Lambda$ to be the set of points $(x, \xi)$ such that there is an integral curve $\gamma(t) \equiv (x(t), \xi(t))$ of $X_p$ satisfying $\tilde{F}(x(0)) = \tilde{C}_1$, $\xi(0) = \nabla \tilde{F}(x(0))$ and $(x(T), \xi(T)) = (x, \xi)$. Moreover, the proof of Theorem 5.5 in [23] shows $\{(x, \nabla \tilde{\varphi}_J(x)) | \tilde{F}(x) > \tilde{C}_1\} \subset \Lambda$. Replacing $\Lambda_+$ in the proof of Lemma 3.2.1 in [25] by $\Lambda$, we compute analogously $(d/dt)\tilde{\varphi}_J(x(t)) = \nabla \tilde{\varphi}_J(x(t)) \cdot \dot{x}(t) = 2|\xi(t)|^2 = J |\nabla \tilde{F}|(x(t)) |\dot{x}(t)|$, where we use that $\gamma(t)$ is an integral curve of $X_p$, $\xi(t) = \nabla \tilde{\varphi}_J(x(t))$ and that $\tilde{\varphi}_J$ satisfies the eiconal equation. The latter equation now gives $\int_0^T J |\nabla \tilde{F}|(x(t)) |\dot{x}(t)| \, dt = \tilde{\varphi}_J(x) - \tilde{\varphi}_J(x(0))$. As $x(t) \in \operatorname{supp} J$ for all $t \leq T$, this clearly implies $\rho(x) \leq \tilde{\varphi}_J(x) - \tilde{C}_1$. Since $J = 1$ on $\{\tilde{F} < C_2\}$, it is easy to see that $\rho(x)$ and therefore $\tilde{\varphi}_J(x) - \tilde{C}_1$ equals $\tilde{F}(x) - \tilde{C}_1$ for all $x \in \{\tilde{F} < C_2\}$. Thus for $\tilde{c}^2 R/24$ larger than $\sup_{\{\tilde{C}_1 < \tilde{F} < C_3\}} |\Delta \tilde{F}|/|\nabla \tilde{F}|^2$ and $\tilde{c}^2 R/48$ larger than $20 \sup(\chi')^2/((\tilde{c}L)^2(1 - \chi^2))$, the choice $\varphi \equiv (1 - R\varepsilon)\tilde{\varphi}_J/2$ shows that the second term on the left-hand side of (3.11) is bounded below by

$$(3.14) \quad \begin{aligned} \int_{\{J > 1/2, \tilde{F} > \tilde{C}_1\}} &\left( \frac{R}{16} |\nabla \tilde{F}|^2 - \frac{1}{2} \Delta \tilde{F} - \frac{|\nabla J|^2}{1 - J^2} \right) |e^{\varphi/\varepsilon} u|^2 \, dx \\ &\geq \frac{R}{48} \int_{\{\tilde{C}_1 < \tilde{F} < C_2\}} |\nabla \tilde{F}|^2 |e^{\varphi/\varepsilon} u|^2 \, dx. \end{aligned}$$

We therefore obtain for $C_3$ sufficiently large depending on $\tilde{\mu}_\varepsilon$ and $R_2 + L$ the existence of a constant $C$ depending on $\tilde{c}$ and $\sup_{\{\tilde{C}_1 < \tilde{F} < C_3\}} |\Delta \tilde{F}|/|\nabla \tilde{F}|^2$

satisfying

$$\begin{aligned}(3.15)\quad &\varepsilon \int_{\{\tilde{C}_1<\tilde{F}<C_2\}}|\nabla e^{(1-R\varepsilon)\tilde{F}/(2\varepsilon)}u|^2\,dx \\ &+ (1/C)\int_{\{\tilde{C}_1<\tilde{F}<C_2\}}e^{(1-R\varepsilon)\tilde{F}/\varepsilon}|u|^2\,dx \\ &\leq \frac{\varepsilon}{2}\int_{\partial\Sigma}\partial_n|e^{(1-R\varepsilon)\tilde{F}/(2\varepsilon)}u|^2\,d\sigma + \int_\Sigma e^{(1-R\varepsilon)\tilde{F}/\varepsilon}u\,H_\varepsilon^{\tilde{F}}u\,dx.\end{aligned}$$

This estimate readily implies:

LEMMA 3.2. *There are constants $C \equiv C(F|\{F>C_1\})$, $C_1$ introduced in Assumption 1.2, $R \equiv R(F|\{F>C_1\})$ such that for every $C_2 > C_1$ and $C_3 > C_2, R$ and for every function $h \in \mathcal{C}^2(\Sigma)\cap\mathcal{C}^1(\overline{\Sigma})$ we have*

$$\begin{aligned}(3.16)\quad &\int_{\{C_1<F<C_2\}}\varepsilon|\nabla e^{-CF/2}h|^2 + ((1/C)-\lambda)e^{-CF}|h|^2\,dx \\ &\leq (\varepsilon/2)\int_{\partial\Sigma}\partial_n|e^{-CF/2}h|^2\,d\sigma + \int_\Sigma e^{-CF}h(L_\varepsilon - \lambda)h\,dx,\end{aligned}$$

*where $\Sigma \equiv \{C_1 < F < C_3\}$, provided Assumption 1.2 holds.*

PROOF. Inserting $h \equiv e^{\tilde{F}/(2\varepsilon)}u$ and the definition of $H_\varepsilon^{\tilde{F}}$, we obtain (3.16) with $\tilde{F}$ in place of $F$. Approximating $F$ in $\mathcal{C}^1(\Sigma)$ by a sequence $\tilde{F}_n$ of functions in $\mathcal{C}^\infty(\Sigma)$ and observing that the analogous quantities $\tilde{\mu}_n, \tilde{c}_n, \tilde{C}_{1,n}$ corresponding to $\tilde{F}_n$ tend to $\mu, c, C_1$, respectively, we derive the assertion. □

For a subset $\Sigma \subset \mathbb{R}^d$ and a function $r: \mathbb{R}^d \to (0, \infty)$ we introduce its $r$-neighborhood by

$$(3.17)\qquad \Sigma^r \equiv \{x \in \mathbb{R}^d | \operatorname{dist}(x, \Sigma) < r(x)\}.$$

Recall definition (2.32) of $B_I(\gamma)$ and (2.12) of $\delta$. Combination of the Harnack inequality, Theorem 8.20 in [24] with (3.16) gives:

PROPOSITION 3.3. *Let $F$ satisfy Assumption 1.2. For every regular, open domain $\Omega \subset \mathbb{R}^d$ containing $\{F < C_1\}$ and every bounded, regular subset $\Sigma \subset \Omega \cap \{F > C_1\}$ there exists a constant $C \equiv C(d, F|\Omega \cap \{F > C_1\}, \Sigma)$ such that for every nonnegative function $\phi \in \mathcal{C}_0^2(\Omega)$ the solution $h \in W^{1,2}(\Omega \cap \{F > C_1\}, e^{-F/\varepsilon}\,dx)$ to the boundary value problem*

$$(3.18)\qquad \begin{aligned}(L_\varepsilon - \lambda)h &= 0, \quad 0 \leq \lambda \leq 1/C, \varepsilon \\ h - \phi &\in W_0^{1,2}(\Omega \cap \{F > C_1\}, e^{-F/\varepsilon}\,dx),\end{aligned}$$

*satisfies for all $y \in \Sigma$*

$$(3.19)\qquad h(y) \leq C\varepsilon^{(1-d)/2}\sup_{\{F=C_1\}}\phi.$$

Moreover, if $\mu_\varepsilon > \delta\varepsilon$ for some $\delta > 0$, then there exists $C \equiv C(d, F)$ such that for all $y \in \Omega \cap \{F > C_1\}$ and all $0 \leq \lambda \leq \varepsilon$

$$(3.20) \quad h(y) \leq C \sup_{B(y,\varepsilon\delta(y))} |\nabla F|^d \varepsilon^{-C} e^{(F(y)-C_1)/(2\varepsilon) - \mathrm{dist}(y,\{F<C_1\})/C}.$$

The reason for writing the poor a priori estimate at infinity in (3.20) is that in the last section concerning the distribution function of transition times we shall need some bound on the principal eigenfunction which is uniform in volume.

PROOF OF PROPOSITION 3.3. We first assume that $\Omega$ is bounded. Since $h \leq \tilde{h} \sup_{\{F=C_1\}} \phi$, $\tilde{h} \equiv h^\lambda_{\{F \leq C_1\}, \Omega^c}$ by the weak maximum principle, it suffices to prove the assertion for $\tilde{h}$. By the boundary Hölder estimates (2.15) we may restrict ourselves to the case $\mathrm{dist}(y, \partial\Omega \cup \{F = C_1\}) > \delta(y)$, where $\delta(x)$ is defined in (2.12). Application of (3.16) and the Harnack inequality, Theorem 8.20 in [24], to $\tilde{h} \in \mathcal{C}^2(\Sigma) \cap \mathcal{C}^1(\overline{\Sigma})$ in combination with the condition on $\lambda$ imply the existence of $C(d) > 0$ such that $e^{-CF(y)} \tilde{h}(y)^2 |B(y, \varepsilon\delta(y))|/C(d)$ is bounded above by

$$(3.21) \quad Ce^{-CC_1}\left(\varepsilon|\{F=C_1\}| \sup_{\{F=C_1\}} |\nabla F| + \varepsilon \int_{\{F=C_1\}} \partial_n \tilde{h} \, d\sigma\right),$$

where $C$ is the constant appearing in Lemma 3.2. The assertion follows since by definition $\delta(y) \sup_{B(y,\varepsilon\delta(y))} |\nabla F| = 1/8$ and since the integral equals $\int_{\{F=C_1\}} \partial_n h^\lambda_{\Omega^c \cup \{F \leq C_1\}} \, d\sigma$ by Green's second formula and $\partial_n h^\lambda_{\Omega^c \cup \{F \leq C_1\}} \leq 0$ on $\partial\Omega$ by the Hopf maximum principle and the fact that $h^\lambda_{\Omega^c \cup \{F \leq C_1\}} \geq 1$ and equal to 1 on $\{F = C_1\}$.

For $\Omega$ unbounded, fix a sequence $\Omega_n \subset\subset \Omega$ of regular, open and bounded domains and denote by $h_n$ the solution to the boundary value problem with $\Omega$ replaced by $\Omega_n$. Note that $h_n \uparrow \bar{h}$ so that $\bar{h}$ is $(L_\varepsilon - \lambda)$-harmonic in $\Omega$ by (2.17) and (2.19) as the Hopf maximum principle tells us that Poisson's kernel is nonnegative. Since the solution is unique by the weak maximum principle and since the right-hand side of (3.19) on each $h_n$ does not depend on $n$, the estimate again follows.

Equation (3.20) is a consequence of (3.8). For, as already mentioned, this equation may be rewritten in terms of $L_\varepsilon^{\tilde{F}}$ as

$$(3.22) \quad \begin{aligned} &\varepsilon \int_\Sigma |\nabla e^{\varphi/\varepsilon} v|^2 e^{-\tilde{F}/\varepsilon} \, dx - \frac{1}{\varepsilon} \int_\Sigma |\nabla \varphi|^2 |e^{\varphi/\varepsilon} v|^2 e^{-\tilde{F}/\varepsilon} \, dx \\ &= \frac{\varepsilon}{2} \int_{\partial\Sigma} \partial_n |e^{-(\tilde{F}/2-\varphi)/\varepsilon} v|^2 \, d\sigma + \int_\Sigma e^{-(\tilde{F}-2\varphi)/\varepsilon} v L_\varepsilon^{\tilde{F}} v \, dx, \end{aligned}$$

where $v \equiv e^{\tilde{F}/(2\varepsilon)} u$. Again by a simple approximation argument, we may assume that $\tilde{F} = F$. Let us introduce the function $v \equiv (1 - J_\varepsilon)h$, where $J_\varepsilon$ is some smooth cut-off function equal to 1 on $\{F \leq C_1\}$, equal to zero on $\{F \leq C_1\}^\varepsilon$ and with modulus of its gradient bounded by $C/\varepsilon$. Choose $\varphi(x) \equiv \delta\varepsilon \operatorname{dist}(x, \{F < C_1\})$ and $\Sigma \equiv \{F > C_1\}$. $\varphi = 0$ on $\partial\Sigma$ and it is not difficult to see that $\varphi$ satisfies the eiconal equation $|\nabla\varphi|^2 = (\delta\varepsilon)^2$ in $\Sigma$ (see Exercise 5.7 in [23]). Since on the left-hand side there appears the quadratic form of the operator $L_\varepsilon$ and since $v$ satisfies the boundary condition zero on $\{F = C_1\}$, we obtain from (3.22)

$$(3.23) \qquad (\mu_\varepsilon - \delta^2 \varepsilon) \int_{\Sigma \setminus \{F \leq C_1\}^\varepsilon} |e^{\varphi/\varepsilon} h|^2 e^{-F/\varepsilon} \, dx \leq \varepsilon^N e^{-C_1/\varepsilon}$$

for some $N \equiv N(d)$, where we use that $h$ is bounded by $\varepsilon^N$ for some $N$ in $\{F \leq C_1\}^\varepsilon$. By the Harnack inequality in combination with the condition on $\mu_\varepsilon$, it follows for some $\delta > 0$ after possibly increasing $N$ that $h(y) \leq \delta(y)^{d/2} \varepsilon^N e^{(F(y)-C_1)/(2\varepsilon) - \delta \operatorname{dist}(y, \{F<C_1\})}$ for all $y \notin \{F \leq C_1\}^\varepsilon$ which implies the assertion if $\Omega$ is bounded. By the same approximation argument as given above we derive the estimate in the general case. $\square$

3.3. *Laplace transforms.* We now want to sharply compare in compact sets eigenfunctions to linear combinations of electrostatic equilibrium potentials $h^0_{A,B} \equiv H^0_{(A \cup B)^c} \mathbb{1}_A$ for small neighborhoods $A$ and $B$ of *relevant* local minima. More precisely, let $\phi \in W^{1,2}_0(\Omega, e^{-F/\varepsilon} dx)$ be a solution of the eigenvalue problem (2.1) for a regular domain $\{F < C_1\} \subset \Omega \subset \mathbb{R}^d$. By Lemma 2.3 for every $x \in \mathcal{M}$ we find $\tilde{x} \in B(x, \varepsilon^{1/(1+\beta)})$ such that $\phi$ does not change sign in the ball $B(\tilde{x}, \varepsilon)$ for some $\beta \equiv \beta(F) > 0$. Let $\tilde{\mathcal{M}} \equiv \tilde{\mathcal{M}}_\phi$ be a collection of such points. Let us define

$$(3.24) \quad \Omega_0 \subset\subset \Omega_1 \subset\subset \Omega_2 \subset\subset \Omega \qquad \text{where } \Sigma \subset\subset \Gamma \text{ stands for } \overline{\Sigma} \subset \Gamma$$

via $\Omega_0 \equiv B_{\tilde{I}}$ for $\tilde{I} \subset \tilde{\mathcal{M}}$, where the former set was defined before (2.32), $\Omega_1 \equiv \{F < C_1\} \cap \Omega$ where $C_1$ is given in Assumption (1.2) and $\Omega_2 \equiv \{F < C_2\} \cap \Omega$ for some large constant $C_2 > C_1$. Clearly, by, for example, (2.17)

$$(3.25) \qquad \phi = \sum_{y \in \tilde{I}} \phi^\lambda_y, \qquad \phi^\lambda_y \equiv H^\lambda_{\Omega \setminus \overline{\Omega}_0} \mathbb{1}_{\partial B_y} \phi$$

in $\Omega_2 \setminus \overline{\Omega}_0$ provided $\lambda < \lambda(\Omega \setminus \overline{\Omega}_0)$. Generally speaking, $\phi$ and $\phi^0$, where we abbreviate $\phi^\lambda \equiv \sum_{y \in \tilde{\mathcal{M}}} \phi^\lambda_y$, are not close to each other everywhere in unbounded regions $\Omega$ even if $\Omega_2 \equiv \Omega$. In fact, we allow $\Omega$ to be equal to $\mathbb{R}^d$ and in this case $\phi^0$ stays bounded while in general $\phi$ is unbounded near infinity. However, exploiting the drift of $F$ toward the local minima, we can show that $\phi$ is close to $\phi^0$ in bounded regions $\Omega_1 \setminus \overline{\Omega}_0$ independent of $\varepsilon$ and containing all relevant local minima. For similar problems in discrete space we refer the reader to [18].

Let us first generalize Lemma 3.1 to the noncompact case. Recall the definition of the maximal time scale $T_I$, $I \subset \mathcal{M}$, given in (2.31). Combining this lemma with the weighted estimates written in Proposition 3.3, we can prove:

PROPOSITION 3.4. *Assume that $F$ satisfies Assumption 1.2 and let $\Omega$ be a regular domain independent of $\varepsilon > 0$ and containing $\{F < C_1\}$, where $C_1$ is defined in Assumption 1.2. Fix $I \subset \mathcal{M}$ and let $\Omega_0$ be a union of $|I|$ balls $B_y \equiv B(\tilde{y}, \varepsilon/4) \subset B(y, \varepsilon^{1/(1+\beta)})$, $y \in I$. There are $N \equiv N(d) \geq 0$ and $\beta \equiv \beta(F) > 0$ such that for all $C_2 > C_1$ and $R > 0$ we find $C \equiv C(d, F, C_2, \{F < C_2\}, R)$ with the following property. For all $0 \leq \lambda \leq \varepsilon^N / T_I$, all nonnegative $f \in L^\infty(\mathbb{R}^d) \setminus 0$ satisfying $\sup_{\partial B_y} f \leq R \inf_{\partial B_y} f$ and all $x \in \{F < C_1\} \setminus \overline{\Omega}_0$ it follows that*

$$(3.26) \quad \begin{aligned} H^\lambda_{\Omega \setminus \overline{\Omega}_0} &\mathbb{1}_{\partial \Omega_0} f(x) \\ &\leq (1 + C\varepsilon^{-N}(\lambda T_I + e^{-(C_2-C_1)/\varepsilon})) H^0_{\{F < C_2\} \setminus \overline{\Omega}_0} \mathbb{1}_{\partial \Omega_0} f(x) \end{aligned}$$

*and for all $y \in I$*

$$(3.27) \quad w^\lambda_{B_y, B_I \cup \Omega^c}(x) \leq (1 + C\varepsilon^{-N} \lambda T_I) w^0_{B_y, B_I \cup \Omega_2^c}(x) + \varepsilon^{-N} T_I e^{-(C_2 - F(x))/\varepsilon}.$$

Before we turn to the proof of this proposition, we note the following a priori lower bound on principal eigenvalues in unbounded domains, which is an immediate consequence of (3.26).

COROLARRY 3.5. *In the situation of the previous proposition we have*

$$(3.28) \quad \lambda(\Omega \setminus \overline{B}_{\tilde{I}}) \geq \frac{\varepsilon^{-N}}{T_I}.$$

At this point we shall need the following consequence of Proposition 4.7 in [4] giving decay of the Poisson kernel in bounded sets. We identify $H^\lambda_\Sigma$ as an operator acting on functions defined on $\mathbb{R}^d$ via $H^\lambda_\Sigma = \mathbb{1}_{\Sigma \cup \partial \Sigma} H^\lambda_\Sigma \mathbb{1}_{\partial \Sigma}$ and likewise for functions a priori defined on $\Sigma$, which by definition take the value zero outside $\Sigma$.

LEMMA 3.6. *There is $C \equiv C(d, F|_{\Omega_2})$ such that for all $\alpha < \beta - C\varepsilon \log \varepsilon \leq C_2 + 1$ and all regular, open, connected sets $\Sigma \subset K \subset \Gamma$ satisfying $K \subset \{F \leq \alpha\}$ and $\{F < \beta\} \subset \Gamma$*

$$(3.29) \quad \|\mathbb{1}_{K \setminus \overline{\Sigma}} H^0_{\Gamma \setminus \overline{\Sigma}} \mathbb{1}_{\partial \Gamma}\| = \sup_{K \setminus \overline{\Sigma}} h^0_{\Gamma^c, \overline{\Sigma}} \leq C e^{-(\beta - \alpha)/\varepsilon} / \varepsilon.$$

PROOF. The existence of $C \equiv C(d, F|\Omega_2)$ follows from Propositions 4.3 and 4.7 in [4] such that for all $x \in K$

$$
\begin{aligned}
h^0_{\Gamma^c, \Sigma}(x) &\leq C \frac{\text{cap}^0_{B(x, \delta(x)\varepsilon)}(\Gamma^c)}{\text{cap}^0_{B(x, \delta(x)\varepsilon)}(\Gamma^c \cup \Sigma)} \\
&\leq C e^{-(\beta-\alpha)/\varepsilon} \bigg/ \bigg( \delta(x) \sup_{B(x, \varepsilon\delta(x))} |\nabla F| \bigg).
\end{aligned}
\tag{3.30}
$$

Equation (3.29) holds since $\delta(x) \sup_{B(x,\varepsilon\delta(x))} |\nabla F| = 1/8$. $\square$

PROOF OF PROPOSITION 3.4. Let $h^\lambda$ denote the function on the left-hand side of (3.26) and set $\Omega_2 \equiv (\Omega \setminus \overline{\Omega}_0) \cap \{F < C_2\}$. We first note that (2.33) gives the existence of constants $N \equiv N(d)$ and $C \equiv C(d, F|\{F < C_1\}, |\{F < C_1\}|)$ such that $s^0(B_{\mathcal{M}\setminus I}, B_I \cup \{F > C_1\}) < C\varepsilon^{-N}$. Therefore, for $\lambda < \varepsilon^N/(2C)$ Lemma 3.1 in combination with (3.19) for $\Sigma \equiv \Omega_2 \cap \{F > C_1\}$ and $h \equiv H^\lambda_{\Omega_2 \cap \{F>C_1\}} \mathbb{1}_{\{F=C_1\}} h^\lambda_{B_{\mathcal{M}\setminus I}, B_I \cup \Omega_2^c}$ and the weak maximum principle yields

$$
\begin{aligned}
\sup_{\Omega_2 \setminus \overline{\Omega}_0} h^\lambda &\leq \sup_{\partial \Omega_0} f \sup_{\Omega_2 \setminus \overline{\Omega}_0} h^\lambda_{B_{\mathcal{M}\setminus I}, B_I \cup \Omega_2^c} \\
&\leq C\varepsilon^{(1-d)/2} \sup_{\partial \Omega_0} f \sup_{\{F=C_1\}} h^\lambda_{B_{\mathcal{M}\setminus I}, B_I \cup \Omega_2^c} \\
&\leq C\varepsilon^{-N} \sup_{\partial \Omega_0} f,
\end{aligned}
\tag{3.31}
$$

where $C \equiv C(d, F|\Omega, \Omega_2)$. We now use the equation

$$
h^\lambda = H^\lambda_{\Omega_2 \setminus \overline{\Omega}_0} \mathbb{1}_{\partial \Omega_0} f + H^\lambda_{\Omega_2 \setminus \overline{\Omega}_0} \mathbb{1}_{\partial \Omega_2} h^\lambda
\tag{3.32}
$$

in $\Omega_1 \setminus \overline{\Omega}_0$ so that by the condition on $f$, (3.31) and (3.1) in combination with (2.33) for $A \equiv \overline{\Omega}_0$ and $B \equiv \Omega_2^c$ for some $C \equiv (d, F|\Omega, \{F < C_2\}, R)$, $N \equiv N(d)$ as above and all $x \in \overline{\Omega}_1 \setminus \overline{\Omega}_0$,

$$
\begin{aligned}
\frac{h^\lambda(x)}{H^\lambda_{\Omega_1 \setminus \overline{\Omega}_0} \mathbb{1}_{\partial \Omega_0} f(x)} &\leq 1 + C\varepsilon^{-N} \frac{h^\lambda_{\Omega_2^c, \overline{\Omega}_0}(x)}{h^\lambda_{\overline{\Omega}_0, \Omega_2^c}(x)} \\
&\leq 1 + C\varepsilon^{-N} \frac{h^0_{\Omega_2^c, \overline{\Omega}_0}(x)}{1 - h^0_{\Omega_2^c, \overline{\Omega}_0}(x)} \\
&\leq 1 + C\varepsilon^{-N} e^{-(C_2 - \tilde{C}_1)/\varepsilon},
\end{aligned}
\tag{3.33}
$$

where we have used Corollary 4.7 in [4] in the latter inequality. To derive the result, it remains to compare $H^\lambda_{\Omega_1 \setminus \overline{\Omega}_0} \mathbb{1}_{\partial \Omega_0} f$ with $H^0_{\Omega_1 \setminus \overline{\Omega}_0} \mathbb{1}_{\partial \Omega_0} f$. But (3.1) once more for $A \equiv \overline{B}_y$, $y \in I$, and $B \equiv \Omega_2^c$ shows after possibly increasing $N$ for some $C \equiv C(d, R)$, some $N \equiv N(d)$ and all $x \in \Omega_1 \setminus \overline{\Omega}_0$

$$
\frac{H^\lambda_{\Gamma_1} \mathbb{1}_{\partial \Omega_0} f(x)}{H^0_{\Gamma_1} \mathbb{1}_{\partial \Omega_0} f(x)} \leq 1 + C\varepsilon^{-N} \lambda T_I
\tag{3.34}
$$

again by the condition on $\lambda$. This proves the first equation in (3.26).

For the proof of (3.27) we note that we already have proven (3.28). Therefore, the Cauchy inequality in combination with (3.26) for $f \equiv \mathbb{1}_{B_x}$ implies for large $N$

$$(3.35) \qquad \sup_{\Omega_1 \setminus \overline{\Omega}_0} w^\lambda_{B_x, B_I \cup \Omega^c} \leq \varepsilon^{-N} T_I.$$

Similarly to the argumentation in (3.31), (2.17) in combination with (3.35) and (3.29) gives for all $y \in \Omega_1 \setminus \overline{\Omega}_0$

$$(3.36) \qquad \begin{aligned} w^\lambda_{B_x, B_I \cup \Omega^c}(y) &= w^\lambda_{B_x, B_I \cup \Omega_2^c}(y) + H^\lambda_{\Omega_1 \setminus \overline{\Omega}_0} \mathbb{1}_{\partial \Omega_2} w^\lambda_{B_x, B_I \cup \Omega^c}(y) \\ &\leq w^\lambda_{B_x, B_I \cup \Omega_2^c}(y) + \varepsilon^{-N} T_I e^{-(C_2 - F(y)/\varepsilon}/\varepsilon. \end{aligned}$$

The assertion in (3.27) now follows from (2.33) in combination with (3.2).

For unbounded $\Omega$ an argument analogous to that given at the end of the proof of Lemma 3.1 shows (3.26) and (3.27) since the constants are uniform in $\Omega$. $\square$

An immediate corollary from Proposition 3.4 is the following relation between eigenfunctions with small eigenvalue and equilibrium potentials. Recall that we have defined for $\beta > 0$

$$(3.37) \qquad A^\beta_{I,J} \equiv \{ y \in \mathbb{R}^d | \hat{F}(y, I) \leq \hat{F}(y, J \setminus I) - \beta \}.$$

COROLARRY 3.7. *Choose $N \equiv N(d)$, $\beta \equiv \beta(F) > 0$ and $C \equiv C(d, F, C_2, \{F < C_2\})$ as in Proposition 3.4 and let $\phi$ be a solution of (2.1) such that $0 \leq \lambda \leq \varepsilon^N / T_{\tilde{I}}$, where $\tilde{\mathcal{M}}$ is defined before (3.24). Let $\phi_y$, $y \in \tilde{I}$, be defined in (3.25). For large $C_2$ it follows for all $x \in \Omega_1 \setminus \overline{\Omega}_0$*

$$(3.38) \quad \phi(x) = \sum_{y \in \tilde{I}} \phi_y(x), \qquad \phi_y(x) = (1 + \mathcal{O}(1)\varepsilon^{-N}\lambda T_{\tilde{I}}) H^0_{\Omega_2 \setminus \overline{\Omega}_0} \mathbb{1}_{\partial B_y} \phi(x).$$

*Moreover, let $x \in \tilde{I}$ be such that $\sup_{\partial B_x} |\phi| = \sup_{\partial B_{\tilde{I}}} |\phi|$. Then*

$$(3.39) \qquad \phi(y) = (1 + \mathcal{O}(1)\varepsilon^{-N}(\lambda T_{\tilde{I}} + e^{-\beta/\varepsilon}))\phi(x), \qquad y \in A^\beta_{x,\tilde{I}},$$

*where the modulus of the Landau symbols is dominated by $C$.*

PROOF. We may assume that $\phi$ is normalized such that it is positive on $B_y$. The Harnack inequality ensures the existence of $C \equiv C(d, F|B(y, \sqrt{\varepsilon}))$ such that $\sup_{B_y} \phi \leq C \inf_{B_y} \phi$. Equation (3.26) for $f \equiv \mathbb{1}_{\partial B_y} \phi$ gives (3.38).

For the proof of (3.39) we first note that Corollary 4.8 in [3] in combination with (3.38) and the condition on $x$ implies for all $y \in \tilde{I} \setminus x$ and all $z \in \Sigma \equiv A^\beta_{x,\tilde{I}}$

$$(3.40) \qquad |\phi_y(z)| \leq 2 h^0_{B_y, B_{\tilde{I}}}(z) \sup_{\partial B_x} |\phi| \leq C e^{-\beta/\varepsilon} |\phi(x)|/\varepsilon,$$

where we have used the Harnack inequality to replace the supremum in the latter inequality. Furthermore, for all $z \in \Sigma$,

$$(3.41) \qquad |\phi_x(z)| \geq (1/2)(1 - h^0_{B_{\tilde{I} \setminus x}, B_x}(z)) \inf_{B_x} |\phi_x| \geq (1/3)|\phi(x)|.$$

In particular, $\phi$ does not change sign in $\Sigma$. In view of (3.38) again and (3.40) we now may choose $r \equiv C\varepsilon^{-N}(\lambda T_{\tilde{I}} + e^{-\beta/\varepsilon})$ in (2.11) applied to $\phi$ and $\Gamma \equiv \Omega$. Since $\{F < C_1\} \subset \Gamma$, (3.1) for $A \equiv \Omega_1^c$ and $B \equiv B_x$ and (3.29) for $\alpha \equiv \sup_{B_x} F$ and $\beta \equiv C_1$ show $\sup_{\Sigma \setminus B_x} h^\lambda_{\Omega^c, B_x} \leq Ce^{-(C_1 - \alpha)/\varepsilon}/\varepsilon$, and (3.39) follows from (2.11); note that we may replace $\varepsilon$ there by $\varepsilon/4$ without any harm. $\square$

**4. Small eigenvalues.** In this section we derive precise asymptotics of the exponentially small eigenvalues of $L_\varepsilon$ in a regular domain $\Omega$ containing $\{F < C_1\}$, where $C_1$ was introduced in Assumption 1.2. We first relate these eigenvalues to the capacity matrix defined in (4.1). In the last section we show that for generic $F$ they are exponentially close to certain principal eigenvalues. In the following section we therefore study principal eigenvalues in detail.

4.1. *Sharp uncertainty principle.* In the sequel we want to derive necessary conditions on small eigenvalues by relating them to a matrix which in leading order equals the capacity matrix introduced in [28]. Namely, fix an eigenvalue $0 \leq \lambda < \varepsilon$ with corresponding eigenfunction $\phi$. Recall the choice of $\tilde{\mathcal{M}} \equiv \tilde{\mathcal{M}}_\phi$ given before (3.24). The existence of $\tilde{\mathcal{M}}$ is guaranteed by Lemma 2.3. Recall the definition of $B_{\tilde{I}}$, $\tilde{I} \subset \tilde{\mathcal{M}}$, given in (1.7). We define $C_{\tilde{I}}(\lambda, \phi)$ to be the matrix with entries

$$(4.1) \quad C_{\tilde{I}}(\lambda, \phi)_{yz} \equiv \varepsilon \int_{\partial B_z} e^{-F/\varepsilon} \frac{\phi}{\phi(z)} \partial_n h^\lambda_{y, \tilde{I}} d\sigma - \delta_{zy} \lambda \int_{B_y} e^{-F/\varepsilon} \frac{\phi}{\phi(y)} dx,$$

where $y, z \in \tilde{I}$. For the sake of convenience we henceforth write shorthand for $\tilde{I}, \tilde{J} \subset \tilde{\mathcal{M}}$:

$$(4.2) \qquad \begin{aligned} h^\lambda_{\tilde{I}, \tilde{J}} &\equiv h^\lambda_{B_{\tilde{I}}, B_{\tilde{J}} \cup \Omega^c}, \qquad \mathrm{cap}^\lambda(\tilde{I}, \tilde{J}) \equiv \mathrm{cap}^\lambda_{B_{\tilde{I}}}(B_{\tilde{J}} \cup \Omega^c), \\ \lambda_{\tilde{I}} &\equiv \lambda(\Omega \setminus B_{\tilde{I}}). \end{aligned}$$

Note that the choice of $\tilde{y} \in B(y, \sqrt{\varepsilon})$, $\tilde{y} \in \tilde{I}$, $y \in \mathcal{M}$, a priori depends on $\phi$.

LEMMA 4.1. *Let $0 \leq \lambda < \lambda_{\tilde{I}}$ for some $\tilde{I} \subset \tilde{\mathcal{M}}$. Then $\lambda \in \sigma(L^\Omega_\varepsilon)$ implies $\det C_{\tilde{I}}(\lambda, \phi) = 0$. Moreover, the vector $\vec{\phi} \equiv (\phi(\tilde{y}))_{\tilde{y} \in \tilde{I}}$ solves the system $C_{\tilde{I}}(\lambda, \phi)\vec{\phi} = 0$.*

PROOF. This characterization is a consequence of Green's second formula applied to $B_{\tilde{I}}$ and $\Omega\setminus\overline{B}_{\tilde{I}}$ showing for all $y \in \tilde{I}$

$$
\begin{aligned}
0 &= \int_\Omega e^{-F/\varepsilon} h^\lambda_{y,\tilde{I}}(L-\lambda)\phi\,dx \\
&= \sum_{z \in \tilde{I}} \int_{\partial B_z} e^{-F/\varepsilon} \phi\, \partial_n h^\lambda_{y,\tilde{I}}\,d\sigma - \lambda \int_{B_y} e^{-F/\varepsilon}\phi\,dx,
\end{aligned}
\tag{4.3}
$$

where we have used that the normal derivative of $\phi$ at $\partial B_{\tilde{I}}$ taken from inside $\Omega\setminus\overline{B}_{\tilde{I}}$ equals the negative of the normal derivative taken from inside $B_{\tilde{I}}$. $\square$

Lemma 4.1 can be used to analyze principal eigenvalues leading to the sharp uncertainty principle Theorem 4.2. As in [3] or [18] one proves $\lambda_{\tilde{I}} < \lambda_{\tilde{I}\cup x}$ for $x \in \tilde{\mathcal{M}}\setminus\tilde{I}$ so that from Lemma 4.1 it follows that

$$
\varepsilon \int_{\partial B_x} e^{-F/\varepsilon} \phi_{\tilde{I}}\, \partial_n h^{\lambda_{\tilde{I}}}_{x,\tilde{I}}\,d\sigma - \lambda_{\tilde{I}} \int_{B_x} e^{-F/\varepsilon} \phi_{\tilde{I}}\,dy = 0,
\tag{4.4}
$$

where $\phi_{\tilde{I}}$ is the principal eigenfunction of $L^\Omega_\varepsilon$ such that $\phi_{\tilde{I}}(x) = 1$. This equation implies:

THEOREM 4.2. *There exist $N \equiv N(d)$ and $C \equiv C(d,F)$ with the following properties. For nonempty, properly contained $\tilde{I} \subset \tilde{\mathcal{M}}$ and $x \in \mathcal{M}\setminus\tilde{I}$ such that $T_{x,\tilde{I}} = T_{\tilde{I}} < \varepsilon^N T_{\tilde{I}\cup x}$ we have*

$$
\lambda_{\tilde{I}} = \left(1 + \mathcal{O}(1)\varepsilon^{-N}\frac{T_{\tilde{I}\cup x}}{T_{\tilde{I}}}\right)\frac{\mathrm{cap}^0(x,\tilde{I})}{\int_{A_{x,\tilde{I}}} e^{-F/\varepsilon}\,dy},
\tag{4.5}
$$

*where $A_{x,\tilde{I}} \equiv A^0_{x,\tilde{I}}$ was defined in (3.37) and where the modulus of the Landau symbol is dominated by $C$.*

PROOF. Taylor's formula shows for all $y \in B^r_x \setminus B_x$, $0 < r < \varepsilon$, and some $0 < \lambda_0 \equiv \lambda_0(y) < \lambda_{\tilde{I}}$,

$$
h^{\lambda_{\tilde{I}}}_{x,\tilde{I}}(y) = h^0_{x,\tilde{I}}(y) + \lambda_{\tilde{I}} w^0_{x,\tilde{I}}(y) + \lambda^2_{\tilde{I}} \dot{w}^{\lambda_0}_{x,\tilde{I}}(y)/2.
\tag{4.6}
$$

By the Cauchy inequality and (3.28) we may bound for some universal $C$

$$
\dot{w}^{\lambda_0}_{x,\tilde{I}}(y) \leq (C/(\lambda_1 - \lambda_0)) w^{\lambda_1}_{x,\tilde{I}}(y) \leq C\varepsilon^{-N} T_{\tilde{I}\cup x} w^{\lambda_1}_{x,\lambda_1}(y)
\tag{4.7}
$$

for some $\lambda_1 < e^N/T_{\tilde{I}\cup x}$ not depending on $y$ so that from (4.6) it follows that

$$
h^{\lambda_{\tilde{I}}}_{x,\tilde{I}}(y) = h^0_{x,\tilde{I}}(y) + \lambda_{\tilde{I}} w^0_{x,\tilde{I}}(y) + \mathcal{O}(1)\lambda_{\tilde{I}}\frac{T_{\tilde{I}\cup x}}{T_{\tilde{I}}} w^{\lambda_1}_{x,\tilde{I}}(y),
\tag{4.8}
$$

where we have used (3.28) once more. From (4.8) and (4.6) we obtain for $y \in \partial B_x$ the double side estimate

$$(4.9) \quad \frac{\lambda_{\tilde{I}}}{C\varepsilon^N} \frac{T_{\tilde{I} \cup x}}{T_{\tilde{I}}} \partial_n w_{x,\tilde{I}}^{\lambda_1}(y) \leq \partial_n h_{x,\tilde{I}}^{\lambda_{\tilde{I}}}(y) - \partial_n h_{x,\tilde{I}}^0(y) - \lambda_{\tilde{I}} \partial_n w_{x,\tilde{I}}^0(y) \leq 0,$$

where the normal derivative is taken from outside $B_x$. Denote by $\beta(F) > 0$ the optimal Hölder exponent of $F$ around $\mathcal{M}$ and fix $\beta \in (\beta(F)/2, \beta(F))$ and define $B_{x,\beta} \equiv B_x(\varepsilon^{1/(1+\beta)})$. We want to estimate for $\lambda \equiv 0, \lambda_1$

$$(4.10) \quad \begin{aligned} &\varepsilon \int_{\partial B_x} e^{-F/\varepsilon} h_{B_x, \mathbb{R}^d \setminus B_{x,\beta}}^{\lambda} \partial_n w_{x,\tilde{I}}^{\lambda} d\sigma \\ &= -\int_{B_{x,\beta} \setminus B_x} e^{-F/\varepsilon} h_{B_x, \mathbb{R}^d \setminus B_{x,\beta}}^{\lambda} dy \\ &\quad + \varepsilon \int_{\partial B_{x,\beta}} e^{-F/\varepsilon} w_{x,\tilde{I}}^{\lambda} \partial_n h_{B_x, \mathbb{R}^d \setminus B_{x,\beta}}^{\lambda} d\sigma, \end{aligned}$$

where the latter equality uses Green's second formula. By (3.27) $w_{x,\tilde{I}}^{\lambda}$ is bounded by

$$(4.11) \quad (1 + C\lambda T_{\tilde{I} \cup x}) w_{x,\tilde{I} \cup \Sigma^c}^0 + C\varepsilon^{-N} e^{-(C_2 - F)/\varepsilon} \sup_{F \leq C_1} T_{\tilde{I} \cup x}$$

for some $C \equiv C(d, F|\{F \leq C_2\}, |\{F \leq C_2\}|)$, where in slight abuse of notation $w_{x,\tilde{I} \cup \Sigma^c}^{\lambda} \equiv w_{B_x, B_{\tilde{I}} \cup \Sigma^c}^{\lambda}$ for $\Sigma \equiv \{F < C_2\}$. Combination of (4.10) and (4.11) with (3.1), (3.2) and (2.33) yields for some $C$ as before and some $N \equiv N(d)$

$$(4.12) \quad \begin{aligned} &\varepsilon \int_{\partial B_x} e^{-F/\varepsilon} \partial_n w_{x,\tilde{I}}^{\lambda} d\sigma \\ &\geq (1 - C\varepsilon^{-N} \lambda T_{\tilde{I} \cup x}) \varepsilon \int_{\partial B_x} e^{-F/\varepsilon} \partial_n w_{x,\tilde{I} \cup \Sigma^c}^0 d\sigma \\ &\quad - C\varepsilon^{-N} \mathrm{cap}_{B_x}^0(\mathbb{R}^d \setminus B_{x,\beta}) T_{\tilde{I} \cup x} e^{-(C_2 - F(x))/\varepsilon}/\varepsilon. \end{aligned}$$

Again by Green's second formula we compute for $\beta > N\varepsilon \log(1/\varepsilon)$

$$(4.13) \quad \begin{aligned} -\varepsilon &\int_{\partial B_x} e^{-F/\varepsilon} h_{x,\tilde{I}}^0 \partial_n w_{x,\tilde{I} \cup \Sigma^c}^0 d\sigma \\ &= \int_{\Sigma \setminus B_{\tilde{I} \cup x}} e^{-F/\varepsilon} h_{x,\tilde{I}}^0 dy \\ &= \sum_{y \in \tilde{I}} \int_{A_{y,x}^{\beta} \setminus B_y} e^{-F/\varepsilon} h_{x,\tilde{I}}^0 dz + \int_{A_{x,\tilde{I}}^{\beta} \setminus B_x} e^{-F/\varepsilon} h_{x,\tilde{I}}^0 dz \\ &\quad + \int_{\Sigma \setminus \bigcup_{y \in \tilde{I}} A_{y,x}^{\beta} \setminus A_{x,\tilde{I}}^{\beta}} e^{-F/\varepsilon} h_{x,\tilde{I}}^0 dz. \end{aligned}$$

Since $F$ is bounded below by $\hat{F}(x,\tilde{I}) - \beta$ on $\Sigma\backslash\bigcup_{y\in\tilde{I}} A^\beta_{y,x}\backslash A^\beta_{x,\tilde{I}}$, we simply bound the last integral in (4.13) by

$$(4.14) \qquad \left|\Sigma\backslash\bigcup_{y\in\tilde{I}} A^\beta_{y,x}\backslash A^\beta_{x,\tilde{I}}\right| e^{-(\hat{F}(x,\tilde{I})-\beta)/\varepsilon}.$$

For the integrals in the sum on the right-hand side by Corollary 4.8 in [4] we may bound $h^0_{x,\tilde{I}}$ on $B^\varepsilon_y\backslash\overline{B}_y$, $y\in\tilde{I}$, by

$$(4.15) \qquad h^0_{x,\tilde{I}} = H^0_{B^\varepsilon_y\backslash\overline{B}_y} \mathbb{1}_{\partial B^\varepsilon_y} h^0_{x,\tilde{I}} \leq Ce^{-(\hat{F}(x,\tilde{I})-F(y))/\varepsilon}/\varepsilon$$

and on $A^\beta_{y,x}\backslash B^\varepsilon_y$ by $Ce^{-(\hat{F}(x,\tilde{I}))-F)/\varepsilon}/\varepsilon$, so that

$$(4.16) \qquad \int_{A^\beta_{y,x}\backslash B_y} e^{-F/\varepsilon} h^0_{x,\tilde{I}}\, dz \leq C|A^\beta_{y,x}\backslash B_y| e^{-F(z^*(x,\tilde{I}))/\varepsilon}/\varepsilon.$$

Concerning the second term on the right-hand side of (4.13), we again use Corollary 4.8 in [4] applied to $h^0_{B_{\tilde{I}},B} \geq h^0_{\tilde{I},x}$, $B \equiv B_x\backslash(\mathbb{R}^d\backslash B_x)^{\varepsilon/100}$, on $A^\beta_{x,\tilde{I}}\backslash B_x$ in combination with Green's second formula showing that for some $N \equiv N(d)$ and some $C \equiv C(d, F|A^\beta_{x,\tilde{I}}, |A^\beta_{x,\tilde{I}}|)$

$$(4.17) \qquad \begin{aligned} &\int_{A^\beta_{x,\tilde{I}}\backslash B_x} e^{-F/\varepsilon} h^0_{x,\tilde{I}}\, dz + \int_{B_x} e^{-F/\varepsilon}\, dz \\ &\geq \int_{A^\beta_{x,\tilde{I}}\backslash B_x} e^{-F/\varepsilon}(1 - h^0_{B_{\tilde{I}},B})\, dz + \int_{B_x} e^{-F/\varepsilon}\, dz \\ &\geq (1 - C\varepsilon^{-N}/T_{\tilde{I}}) \int_{A^\beta_{x,\tilde{I}}} e^{-F/\varepsilon}\, dz. \end{aligned}$$

Inserting (4.17), (4.16) and (4.14) into (4.13) and the result into (4.12), we derive

$$(4.18) \qquad \begin{aligned} &\varepsilon\int_{\partial B_x} e^{-F/\varepsilon} \partial_n w^\lambda_{x,\tilde{I}}\, d\sigma + \int_{B_x} e^{-F/\varepsilon}\, dz \\ &\geq -\left(1 + C\lambda T_{\tilde{I}\cup x} + C\frac{\varepsilon^{-N} + e^{\beta/\varepsilon}}{T_{\tilde{I}}}\right)\int_{A^\beta_{x,\tilde{I}}} e^{-F/\varepsilon}\, dz \\ &\quad - C\varepsilon^{-N} T_{\tilde{I}\cup x} e^{-C_2/\varepsilon}, \end{aligned}$$

where we have used $\mathrm{cap}^0_{B_x}(\mathbb{R}^d\backslash B_{x,\beta}) \leq Ce^{-F(x)/\varepsilon}/\varepsilon^{d-1}$ following from Proposition 4.7 in [4]. Using (4.10) in combination with (4.18) in (4.4), we now conclude

$$(4.19) \qquad \begin{aligned} \mathrm{cap}^0(x,\tilde{I}) &= \left(1 + \mathcal{O}(1)\varepsilon^{-N}\left(\frac{1+e^{\beta/\varepsilon}}{T_{\tilde{I}}} + e^{-(C_2-F(x))/\varepsilon} + \frac{T_{\tilde{I}\cup x}}{T_{\tilde{I}}}\right)\right) \\ &\quad \times \lambda_{\tilde{I}} \int_{A^\beta_{x,\tilde{I}}} e^{-F/\varepsilon}\, dz, \end{aligned}$$

where we have replaced $\phi_{\tilde{I}}$ by $1 + \mathcal{O}(1)e^{-N}(1/T_{\tilde{I}} + e^{-(C_2-F(x))/\varepsilon})$ in view of (3.39) and (2.7). Since $C_2$ can be made arbitrarily large, the proof is completed by choosing $\beta \equiv N\varepsilon\log(1/\varepsilon)$ and that $A^{\beta}_{x,\tilde{I}}$ may be replaced by $A^0_{x,\tilde{I}}$ without harm. $\square$

Equation (4.5) implies the following intuitive sharp uncertainty principle.

COROLARRY 4.3. *In the situation of the previous theorem, we have*

$$(4.20) \qquad \lambda_{\tilde{I}} = \left(1 + \mathcal{O}(1)\varepsilon^{-N}\frac{T_{\tilde{I}\cup x}}{T_{\tilde{I}}}\right)\frac{1}{\mathbb{E}[\tau^x_{B_{\tilde{I}}\cup\Omega^c}]},$$

*where* $\mathbb{E}[\tau^x_{B_{\tilde{I}}\cup\Omega^c}]$ *is the expected time of the first visit of* $B_{\tilde{I}}\cup\Omega^c$ *of the diffusion generated by* $L_\varepsilon$ *and starting in* $x$.

PROOF. We first note that by (2.18) for $\Sigma \equiv \Omega\setminus\{F \geq C_1\}$ and $\Gamma \equiv \Omega\setminus\{F \geq C_2\}$, $C_2 > C_1$,

$$(4.21) \quad \begin{aligned} s^0_{\Sigma\setminus\overline{B}_{\tilde{I}\cup x}}(B_{\tilde{I}} \cup \Omega^c) \\ \leq s^0_{\Sigma\setminus\overline{B}_{\tilde{I}\cup x}} + \sup_{\Sigma} h^0_{\Gamma^c, B_{\tilde{I}\cup x}} s^0_{\partial\Gamma}(B_{\tilde{I}} \cup \Omega^c) + s^0_{\partial B_x}(B_{\tilde{I}} \cup \Omega^c). \end{aligned}$$

Since for large $C_2$ the first term on the right-hand side of (4.21) is bounded by $C\varepsilon^{-N}T_{\tilde{I}\cup x}$ for some $N \equiv N(d)$ and $C \equiv C(d, F|\{F < C_2\}, |\{F < C_2\}|)$ by (3.27) in combination with (2.33), after possibly increasing $C_2$ we obtain $s^0_{\Sigma\setminus\overline{B}_{\tilde{I}\cup x}}(B_{\tilde{I}} \cup \Omega^c) \leq C\varepsilon^{-N}T_{\tilde{I}\cup x} + s^0_{\partial B_x}(B_{\tilde{I}} \cup \Omega^c)$. Invoking (6.1) in [4] in combination with the fact that the nominator in this estimate is computed in (4.17), we conclude

$$(4.22) \qquad s^0_{\Sigma\setminus\overline{B}_{\tilde{I}\cup x}}(B_{\tilde{I}} \cup \Omega^c) \leq \left(1 + C\varepsilon^{-N}\frac{T_{\tilde{I}\cup x}}{T_{\tilde{I}}}\right)s^0_{\partial B_x}(B_{\tilde{I}} \cup \Omega^c).$$

It follows from (2.11), (4.22) and (6.1) in [4] that

$$(4.23) \qquad \sup_{B_x} w^0_{I\cup\Omega^c} \leq \left(1 + C\varepsilon^{-N}\frac{T_{\tilde{I}\cup x}}{T_{\tilde{I}}}\right)\inf_{B_x} w^0_{I\cup\Omega^c}.$$

It follows from Green's second formula and (4.23) that modulo the error term appearing above

$$(4.24) \quad \begin{aligned} w^0_{I\cup\Omega^c}(x)\,\mathrm{cap}^0(x, I) &= \varepsilon\int_{\partial B_x} e^{-F/\varepsilon} w^0_{I\cup\Omega^c} \partial_n h^0_{x,I}\,d\sigma \\ &= -\varepsilon\int_{\partial B_x} e^{-F/\varepsilon} \partial_{-n} w^0_{I\cup\Omega^c}\,d\sigma + \int_{\Omega\setminus\overline{B}_x} e^{-F/\varepsilon} h^0_{x,I}\,dy \\ &= \int_{\Omega} e^{-F/\varepsilon} h^0_{x,I}\,dy. \end{aligned}$$

The latter integral was computed in (4.17) so that the assertion follows from (4.5). □

We want to analyze the diagonal entries of the matrix $C_{\tilde{I}}(\lambda, \phi)$ in more detail.

LEMMA 4.4. *In the situation of the previous theorem there are $N \equiv N(d)$ and $C \equiv C(d, F)$ satisfying, for all $0 \leq \lambda \leq \lambda_0$, $\lambda_0 \equiv \varepsilon^{-N}/T_{\tilde{I} \cup x}$,*

$$
\begin{aligned}
C_{\tilde{I}}(\lambda, \phi_{\tilde{I}})_{xx} = &-(1 + \mathcal{O}(1)\varepsilon^{-N}\lambda T_{\tilde{I} \cup x}) \\
&\times (\lambda - \lambda_{\tilde{I}} - (\lambda - \lambda_{\tilde{I}})^2 \mathcal{O}(1)\varepsilon^{-N} T_{\tilde{I} \cup x}) \int_{A_{x,\tilde{I}}} e^{-F/\varepsilon}\, dy.
\end{aligned}
\tag{4.25}
$$

PROOF. Performing a Taylor expansion at $\lambda \equiv \lambda_{\tilde{I}}$ to second order of the Laplace transform on the left-hand side of (4.25), we compute similarly to (4.9) using (4.4)

$$
\begin{aligned}
C_{\tilde{I}}(\lambda, \phi_{\tilde{I}})_{xx} = (\lambda - \lambda_{\tilde{I}})\bigg(&\varepsilon \int_{\partial B_x} e^{-F/\varepsilon} \phi_{\tilde{I}}\, \partial_n w_{x,\tilde{I}}^{\lambda_{\tilde{I}}}\, d\sigma - \int_{B_x} e^{-F/\varepsilon} \phi_{\tilde{I}}\, dx\bigg) \\
&+ (\lambda - \lambda_{\tilde{I}})^2 \varepsilon \int_0^1 s\, ds \int_{\partial B_x} d\sigma e^{-F/\varepsilon} \phi_{\tilde{I}}\, \partial_n \dot{w}_{x,\tilde{I}}^{(1-s)\lambda_{\tilde{I}} + s\lambda}.
\end{aligned}
\tag{4.26}
$$

Analogously to (4.8) by the Cauchy formula and (3.28) on $\partial B_x$ we have the estimate

$$
\partial_n \dot{w}_{x,\tilde{I}}^{(1-s)\lambda_{\tilde{I}} + s\lambda} \geq (C/(\lambda_1 - \lambda))\, \partial_n w_{x,\tilde{I}}^{\lambda_1} \geq (100 C/\lambda_0)\, \partial_n w_{x,\tilde{I}}^{\lambda_1}
\tag{4.27}
$$

for some universal constant $C$, $(1/100)\lambda_0 \leq \lambda_1 < \lambda_0$ and all $0 \leq \lambda < (1/200)\lambda_0$. Moreover, (4.18) remains valid for $\lambda \equiv \lambda_I, \lambda_1$ so that analogously to (4.19) for $C_2$ large the last term on the right-hand side of (4.26) is of order

$$
(\lambda - \lambda_{\tilde{I}})^2 \varepsilon^{-N} T_{\tilde{I} \cup x} \int_{A_{x,\tilde{I}}^\beta \setminus B_x} e^{-F/\varepsilon}\, dy
\tag{4.28}
$$

while the first term can be estimated by

$$
(\lambda - \lambda_{\tilde{I}})\bigg(1 + \mathcal{O}(1)\varepsilon^{-N}\bigg(\lambda_{\tilde{I}} T_{\tilde{I} \cup x} + \frac{1 + e^{\beta/\varepsilon}}{T_{\tilde{I}}}\bigg)\bigg) \int_{A_{x,\tilde{I}}^\beta} e^{-F/\varepsilon}\, dx.
\tag{4.29}
$$

In view of (3.28) from (4.29) and (4.28) for some $N \equiv N(d)$ and $\beta \equiv N\varepsilon \log(1/\varepsilon)$ the assertion follows. □

**4.2. Small eigenvalues.** We now turn to the investigation of small eigenvalues of $L_\varepsilon$. Namely, we will show how the capacity matrix introduced in (4.1) can be used to analyze the spectrum of $L_\varepsilon^\Omega$ near zero. The proof of Theorem 4.6 proceeds close to the line of arguments of the analogous assertion in [18] or in [3]. In particular, for proofs, which are straightforward generalizations, we refer the reader to the counterparts therein.

In addition to the notation introduced in (2.31) let us define

(4.30) $\mathcal{M}(x) \equiv \{y \in \mathcal{M} | F(y) < F(x)\}, \qquad T_x \equiv T_{x,\mathcal{M}(x)}, \qquad x \in \mathcal{M},$

in case that $\mathcal{M}(x) \neq \varnothing$. We use the convention $T_x \equiv d_x \equiv \infty$ for $\mathcal{M}(x) \equiv \varnothing$.

ASSUMPTION 4.5. $F$ is generic in the sense that $\rho > N\varepsilon \log(1/\varepsilon)$ for some $N \equiv N(d) > 0$, where $\rho$ was introduced in (1.17).

From this assumption and its consequence Lemma 4.8 it follows that $\mathcal{M} \ni x \mapsto T_x$ is injective [see (4.47) for a proof]. We hence obtain an ordering of $\mathcal{M}$ via

(4.31) $\qquad\qquad x < y \quad \text{if and only if} \quad T_x > T_y.$

We also define

(4.32) $\qquad \mathcal{M}_{<x} \equiv \{y \in \mathcal{M} | y < x\}, \qquad \mathcal{M}_{\leq x} \equiv \{y \in \mathcal{M} | y \leq x\}.$

We also shall need

(4.33) $\qquad\qquad \mathcal{T}_x \equiv \min_{y \in \mathcal{M}_{<x}} \min_{z \in \mathcal{M}_{<x} \setminus y} T_{y,z}.$

Let us briefly outline the strategy of finding small eigenvalues. Starting the process in a local minimum of $F$, we believe that for exponentially long times in $1/\varepsilon$ it behaves like the process obtained by reflecting the original one at the boundary of the corresponding valley. For each $x \in \mathcal{M}$, we thus look for a solution of the equation appearing in Lemma 4.1 near the ground-state energy of the associated Dirichlet operator $L_\varepsilon^\Sigma$, where $\Sigma \equiv \Omega \setminus B_{\tilde{\mathcal{M}}_{<x}}$—recall the choice of $\tilde{\mathcal{M}}$ given before (3.24). We then show that these solutions are the only candidates for eigenvalues below $\varepsilon^N$ for some $N$ and that they have to be simple in case. Since from, for example, [39] or [8] we already know that there are $|\mathcal{M}|$ many eigenvalues, these candidates in fact constitute all eigenvalues below $\varepsilon^N$. We use the conventions $T_\varnothing \equiv \infty$, $h_{A,\varnothing}^\lambda \equiv 1$, $A \neq \varnothing$, $h_{\varnothing,B}^\lambda \equiv 0$, $B \neq \varnothing$, $\lambda_\varnothing \equiv 0$ and $\alpha/\infty \equiv 0$ for $\alpha > 0$. The result is

THEOREM 4.6. *For some $N \equiv N(d)$ there is $C \equiv C(d,F)$ dominating all moduli of the Landau symbols appearing below in case that Assumptions 1.2 and 4.5 hold. There are $|\mathcal{M}|$ simple eigenvalues $\lambda_x < \lambda_y$, $x,y \in \mathcal{M}$, $x < y$, satisfying*

(4.34) $\qquad\qquad \sigma(L_\varepsilon^\Omega) \cap [0, \varepsilon^N) = \{\lambda_x | x \in \mathcal{M}\}.$

For every $x \in \mathcal{M}$, $x \neq \min \mathcal{M}$, we have $\mathcal{T}_x \geq e^{\rho/\varepsilon} T_x$ and $T_x \geq e^{\rho/\varepsilon} \times \max_{y \in \mathcal{M} \setminus \mathcal{M}_{\leq x}} T_y$, where $\rho$ is defined in (1.17). Furthermore, there exist $\beta \equiv \beta(F) > 0$, a set $\tilde{\mathcal{M}}_{<x}$ of $|\mathcal{M}_{<x}|$ points such that $\tilde{\mathcal{M}}_{<x} \cap B(y, \varepsilon^{1/(1+\beta)})$, $y \in \mathcal{M}_{<x}$, is a singleton, and

$$(4.35) \quad \lambda_x = \left(1 + \mathcal{O}(1)\varepsilon^{-N}\left(\frac{T_x}{\mathcal{T}_x} + \frac{\max_{y \in \mathcal{M} \setminus \mathcal{M}_{\leq x}} T_y}{T_x}\right)\right)\lambda_{\tilde{\mathcal{M}}_{<x}},$$

where we use the notation introduced in (4.2) for the principal eigenvalue $\lambda_{\tilde{I}}$. Moreover, every eigenfunction $\phi_x$ corresponding to $\lambda_x$ satisfies for all $z \in \{F < C_1\} \cap \Omega$

$$(4.36) \quad \begin{aligned} \phi_x(z) &= \left(1 + \mathcal{O}(1)\varepsilon^{-N}\frac{T_{\mathcal{M}_{\leq x}}}{T_{\mathcal{M}_{<x}}}\right)h^0_{x,\tilde{\mathcal{M}}_{<x}}(z) \\ &\quad + \sum_{y \in \tilde{\mathcal{M}}_{<x}} \mathcal{O}(1)\varepsilon^{-N}\frac{T_x}{T_{y,x}}h^0_{y,\tilde{\mathcal{M}}_{\leq x}}(z). \end{aligned}$$

Combination of 4.2 with (4.35) and (1.11) yields the following.

COROLARRY 4.7. *In the situation of Theorem 4.6 we have for all $x \in \mathcal{M}$*

$$(4.37) \quad \lambda_x = (1 + \mathcal{O}(1)\varepsilon^{-N}e^{-\rho/\varepsilon})\frac{1}{T_{\mathcal{M}(x)}}.$$

*In particular, under the conditions of Theorem 1.1 it follows that*

$$(4.38) \quad \lambda_x = (1 + \mathcal{O}(1)\varepsilon\log(1/\varepsilon))|\lambda^*|\sqrt{\frac{\det \operatorname{Hess} F(x)}{|\det \operatorname{Hess} F(z^*)|}}e^{-(\hat{F}(x,I)-F(x))/\varepsilon}.$$

PROOF. It only remains to show that we may replace $\tilde{\mathcal{M}}(x)$ by $\mathcal{M}(x)$ within the error estimates in $T_{\tilde{\mathcal{M}}(x)} = T_{x,\tilde{I}}$, $\tilde{I} \equiv \tilde{\mathcal{M}}(x)$, where we use (4.43). This is obvious for the nominator of the latter time scale. For the denominator we have by the Hopf maximum principle on $\partial B_x$ for small $\delta > 0$

$$(4.39) \quad \begin{aligned} \partial_n h^0_{x,\tilde{I}} &= -\partial_n h^0_{\tilde{I},x} = -\partial_n H^0_{(B_x \cup B_{\tilde{I}}(\delta))^c}\mathbb{1}_{\partial B_{\tilde{I}(\delta)}}(1 - h^0_{x,\tilde{I}}) \\ &= -(1 - \mathcal{O}(1)e^{-(\hat{F}(x,\tilde{I})-\max F(B_{\tilde{I}}(\delta)))/\varepsilon}/\varepsilon)\partial_n h^0_{B_{\tilde{I}}(\delta),B_x} \\ &= (1 - \mathcal{O}(1)e^{-(\hat{F}(x,I)-\max F(B_{\tilde{I}}(\delta)))/\varepsilon}/\varepsilon)\partial_n h^0_{B_x,B_{\tilde{I}}(\delta)}, \end{aligned}$$

where the second equality again is a consequence of Corollary 4.8 in [4]. We thus obtain

$$(4.40) \quad \operatorname{cap}^0_{B_x}(B_{\tilde{I}}) = (1 - \mathcal{O}(1)e^{-(\hat{F}(x,I)-\max F(B_{\tilde{I}}(\delta)))/\varepsilon}/\varepsilon)\operatorname{cap}^0_{B_x}(B_{\tilde{I}}(\delta)).$$

Applying the same arguments to $I \equiv \mathcal{M}(x)$, we obtain

$$(4.41) \quad T_{x,\tilde{I}} = (1 - \mathcal{O}(1)\varepsilon^{-N}/T_{x,I})T_{x,I}$$

for $\delta \equiv \varepsilon^{1/(1+\beta)}$, where $\beta > 0$ is the Hölder exponent of $F$ locally at $\mathcal{M}$. The assertion then follows from (4.43). □

We have to introduce some more notation. Fix a set $\tilde{\mathcal{M}}$ of cardinality $|\mathcal{M}|$ such that $B(x, \varepsilon^{1/(1+\beta)}) \cap \tilde{\mathcal{M}}$ is a singleton, where $\beta \equiv \beta(F) > 0$ is the constant appearing in Proposition 3.4. $\tilde{\mathcal{M}}$ inherits the ordering of $\mathcal{M}$ in an obvious way. Let us define $\tilde{\mathcal{E}}_x \equiv \infty$ for $x \equiv \min \tilde{\mathcal{M}}$ and for $x \in \tilde{\mathcal{M}} \setminus \min \tilde{\mathcal{M}}$ set

$$\tilde{\mathcal{E}}_x \equiv \min_{y \in \tilde{\mathcal{M}}_{<x}} T_{y, \tilde{\mathcal{M}}_{\leq x} \setminus y}. \tag{4.42}$$

The first lemma actually is a special case of Lemma 4.5 in [18] or Lemma 5.3 in [3].

LEMMA 4.8. *For all $x \in \tilde{\mathcal{M}}$ it follows for some $N \equiv N(d)$, $C \equiv C(d, F) > 0$ if $\rho > C\varepsilon \log \varepsilon$*

$$(1 + \mathcal{O}(1)\varepsilon^{-N} e^{-\rho/\varepsilon}) T_x = T_{\tilde{\mathcal{M}}_{<x}} = T_{x, \tilde{\mathcal{M}}_{<x}}. \tag{4.43}$$

*For $x \in \tilde{\mathcal{M}} \setminus \min \tilde{\mathcal{M}}$ we have*

$$\tilde{\mathcal{T}}_x \geq \tilde{\mathcal{E}}_x \geq e^{\rho/\varepsilon} T_x. \tag{4.44}$$

*Moreover, for $x, y \in \tilde{\mathcal{M}}$, $y < x$,*

$$\max_{z \in \tilde{\mathcal{M}} \setminus \tilde{\mathcal{M}}_{\leq x}} T_{z, \tilde{\mathcal{M}}_{\leq x} \setminus y} \leq e^{-\rho/\varepsilon} T_{y, \tilde{\mathcal{M}}_{\leq x} \setminus y}. \tag{4.45}$$

*In particular, for $y \in \tilde{\mathcal{M}}_{<x}$,*

$$T_{y, \tilde{\mathcal{M}}_{\leq x} \setminus y} = T_{\tilde{\mathcal{M}}_{\leq x} \setminus y}. \tag{4.46}$$

For the proof that the map $x \mapsto T_x$ is one-to-one, we first note that (4.41) also holds for some $\beta \equiv \beta(F) > 0$ and arbitrary $I \subset \mathcal{M}$ and $\tilde{I}$ such that $\tilde{I} \cap B_x(\varepsilon^{1/(1+\beta)})$, $x \in \mathcal{M}$, is a singleton. We then have for $x < y$ by definition and Assumption 4.5

$$\begin{aligned} T_x &\sim T_{\tilde{\mathcal{M}}_{<x}} = \max_{z \in \tilde{\mathcal{M}}_{<x}} T_{z, \tilde{\mathcal{M}}_{<x}} = T_{x, \tilde{\mathcal{M}}_{<x}} \\ &> e^{\rho/\varepsilon} T_{y, \tilde{\mathcal{M}}_{<x}} \geq e^{\rho/\varepsilon} T_{y, \tilde{\mathcal{M}}_{<y}} \sim e^{\rho/\varepsilon} T_y. \end{aligned} \tag{4.47}$$

We would like to explain the geometrical background of the previous crucial lemma. The exponential rate of the time scale $T_{x,I}$ is given by

$$\hat{F}(x, I) - F(x) = \varepsilon \log T_{x,I} + \mathcal{O}(1) \log(1/\varepsilon), \tag{4.48}$$

where $\hat{F}$ is the communication height defined in (1.6). The latter equality is a consequence of Proposition 4.7 in [4], where the Landau symbol denotes a quantity with modulus bounded by a constant $C \equiv C(d, F|\{F \leq C_1\})$. The first observation is that the restriction $\hat{E}$ of the communication height $\hat{F}$ to singletons in $\mathcal{M}$ satisfies the ultrametric triangle inequality, that is, $\hat{E}(x,y) \leq \max(\hat{E}(x,z), \hat{E}(z,y))$ for all $x, y, z \in \mathcal{M}$. Using the convention $\hat{E}(x,x) \equiv 0$, it is also positive definite and symmetric and therefore it is an ultrametric by definition. It is not difficult to see that the ultrametric triangle inequality is equivalent to the assertion that an ultrametric ball is centered at each of its interior points, that is, for $I \subset \mathcal{M}$ and $x \in \mathcal{M} \setminus I$ and all $y \in \mathcal{M} \setminus I$ such that $\hat{E}(y,x) < r \equiv \hat{E}(x,I) \equiv \max_{y \in I} \hat{E}(x,y)$ we have $\hat{E}(y,I) = r$. We would like to point out that the time scales still exhibit this ultrametric structure under very general conditions without knowing (4.48). If, for example, $F \equiv F_\varepsilon$ also depends on $\varepsilon$ with degenerate growing level sets at local minimal values, the process behaves like a Brownian motion when started there. It therefore might be that the process stays in such regions rather because it takes much time for Brownian motion to leave a large set. This feature is taken care of in the definition of the time scale $T_{x,I}$ in (1.9), which then is large not because the capacity $\mathrm{cap}^0(x,I)$ is small but because the invariant measure $\int_{A_{x,I}} e^{-F_\varepsilon/\varepsilon} dx$ of a basin $A_{x,I}$ is fairly large. Under the same genericity assumption as Assumption 4.5 one can still prove ultrametricity and we refer the interested reader to [18], where under minimal conditions this point is made rigorous.

Lemma 4.8 is a special case of Lemma 5.3 in [3]. Indeed, within the notation used therein a glance at the proof of Lemma 5.3 in [3] shows that it actually holds for any set of times scales $T_{x,J}$, $x \in \mathcal{M}_N$, $J \subset \mathcal{M}_N \setminus x$ (depending on a parameter $N \equiv 1/\varepsilon$), on some finite set $\mathcal{M}_N$ such that $(1/N) \log T_{x,J} = e(x,J) - f(x) + \mathcal{O}(1) \log(N)/N$ for some function $f \equiv f_N$ and some ultrametric $e \equiv e_N$ on $\mathcal{M}_N$. The quantity $\delta_N$ appearing there (see Definition 1.2 in [3]) corresponds to $e^{-\rho/\varepsilon}$ in our case choosing $\mathcal{M}_N \equiv \mathcal{M}, \tilde{\mathcal{M}}$. With our favorite choice of the set of time scales we then have $e_N = \hat{E} \equiv \hat{F}|\mathcal{M} = \hat{F}|\tilde{\mathcal{M}}$ and $f_N = F$ in our context.

We now analyze the possible solutions of the capacity matrix introduced in (4.1).

*Henceforth, we assume for some $x \in \mathcal{M} \setminus \min \mathcal{M}$ that $\lambda_x$ is an eigenvalue of $L_\varepsilon^\Omega$ satisfying for some $\alpha > 0$*

$$(4.49) \qquad e^{\alpha/\varepsilon}/\mathcal{E}_x < \lambda_x < e^{-\alpha/\varepsilon} T_{\mathcal{M}_{\leq x}}.$$

*Furthermore, $\phi_x$ denotes a corresponding eigenfunction.*

For $y \in \mathcal{M}_{\leq x}$ according to Lemma 2.3 we choose a set $\tilde{\mathcal{M}}_{\leq x}$ of $|\mathcal{M}_{\leq x}|$ points each of them lying in one ball $B(y, \varepsilon^{1/(1+\beta)})$, $y \in \mathcal{M}_{<x}$, for $\tilde{\beta} \equiv \beta(F)$ appearing in the choice of $\tilde{I}$ in Corollary 3.7 such that $\phi_x$ does not change sign in $B(\tilde{y}, \varepsilon)$. Equation (3.39) implies that $\phi_x$ does not change sign in $B(x, \varepsilon/4)$. Indeed, assume that $|\phi_x|$ attains its supremum in some ball

$B(y, \varepsilon/4)$, where $y \in \tilde{I} \equiv \tilde{\mathcal{M}}_{<x}$, $y \notin B(x, \varepsilon^{1/(1+\beta)})$ in Corollary 3.7. It follows from (3.39) in combination with the upper bound in (4.49) that $\phi_x$ does not change sign in $A^\beta_{y,\tilde{\mathcal{M}}_{\leq x}}$ for arbitrary but fixed $\beta > 0$. Choosing $u \equiv |\phi_x|$ restricted to $\overline{A}^\beta_{y,\tilde{\mathcal{M}}_{\leq x}}$ in (2.6) for $\Sigma \equiv A^\beta_{y,\tilde{\mathcal{M}}_{\leq x}}$, it follows $\lambda(\Sigma) \geq \lambda_x$. On the other hand, we obtain from (2.7) in combination with Proposition 4.7 in [4] that $\lambda(\Sigma) \leq \varepsilon^{-N} \varepsilon^{\beta/\varepsilon} T_{y,\tilde{\mathcal{M}}_{\leq x}\backslash y} \leq \varepsilon^{-N} \varepsilon^{\beta/\varepsilon}/\tilde{\mathcal{E}}_x$. Since $\tilde{\mathcal{E}}_x \leq \varepsilon^N \mathcal{E}_x$, we thus have derived a contradiction to the lower bound in (4.49) for $\beta < \alpha$. In view of (3.39) we now may assume that $\tilde{\mathcal{M}}_{\leq x} \cap B(x, \varepsilon^{1/(1+\beta)}) = x$. Note that by the obvious generalization of (4.41) the various time scales in the error estimates appearing in Theorem 4.6 may be replaced without harm by those defined with respect to $\tilde{\mathcal{M}}_{<x}$ instead. In addition, the replacement of $\mathcal{M}_{\leq x}$ by $\tilde{\mathcal{M}}_{\leq x}$ changes $\alpha$ only by an amount of order $N\varepsilon \log(1/\varepsilon)$. Therefore,

*In the sequel we identify $\mathcal{M}_{<x}$ with the set $\tilde{\mathcal{M}}_{<x}$.*

Let $C_x \equiv \text{diag}(e^{F(y)/\varepsilon})_{y \in \mathcal{M}_{\leq x}} C_{\mathcal{M}_{\leq x}}(\lambda, \phi_x)$, where $C_{\mathcal{M}_{\leq x}}(\lambda, \phi_x)$ is the matrix defined in (4.1) for $\tilde{I} \equiv \mathcal{M}_{\leq x}$. Note that by Green's second formula and $h^\lambda_{z,\mathcal{M}_{\leq x}}|\partial B_z = 1$ it follows that the latter matrix is symmetric. Hence

$$
\begin{aligned}
&\begin{pmatrix} \mathcal{K}_x(\lambda) & -\vec{g}_x(\lambda) \\ -(\text{diag}(e^{(F(x)-F(y))/\varepsilon})_{y \leq x} \vec{g}_x(\lambda))^t & e^{F(x)/\varepsilon} C_{\mathcal{M}_{\leq x}}(\lambda, \phi_x)_{xx} \end{pmatrix} \\
&\equiv C_x(\lambda) \\
(4.50) \quad &= \left( \varepsilon \int_{\partial B_z} e^{-(F-F(z))/\varepsilon} \frac{\phi_x}{\phi_x(z)} \partial_n h^\lambda_y \, d\sigma \right. \\
&\qquad \left. - \delta_{zy} \lambda \int_{B_y} e^{-(F-F(y))/\varepsilon} \frac{\phi_x}{\phi_x(y)} \, du \right)_{zy \in \mathcal{M}_{\leq x}}.
\end{aligned}
$$

During the rest of the section we write shorthand $h^\lambda_y \equiv h^\lambda_{y,\mathcal{M}_{\leq x}}$. Let us furthermore define

$$(4.51) \qquad \mathcal{N}_x \equiv \mathcal{D}_x - \mathcal{K}_x,$$

where

$$
(4.52) \quad \mathcal{D}_x \equiv \text{diag}\left( \varepsilon \int_{\partial B_y} e^{-(F-F(y))/\varepsilon} \frac{\phi_x}{\phi_x(y)} \partial_n h^\lambda_y \, d\sigma - \lambda \int_{B_y} e^{-(F-F(y))/\varepsilon} \frac{\phi_x}{\phi_x(y)} \, du \right)_{y<x}.
$$

Equipped with the ultrametric structure in the form written in Lemma 4.8 and the control of Laplace transforms and of eigenfunctions obtained in the previous section, one simply can write a Neumann series [see (4.58); recall that a matrix $A$ is invertible if the series $\sum_{k \geq 0} \|\mathbb{1} - A\|^k$ converges in one multiplicative norm $\|\cdot\|$ in which case $A^{-1} = \sum_{k \geq 0} (\mathbb{1} - A)^k$] for $\mathbb{1} - \mathcal{D}_x(\lambda)^{-1} \mathcal{N}_x(\lambda)$ for $\lambda$ near $\lambda_{\mathcal{M}_{<x}}$ proving invertibility of $\mathcal{K}_x(\lambda)$. We then

compute

$$(4.53) \quad \det C_x = \det \begin{pmatrix} \mathcal{K}_x & 0 \\ -(\mathrm{diag}(e^{(F(x)-F(y))/\varepsilon})_{y \leq x} \vec{g}_x)^t & G_x \end{pmatrix} = G_x \det \mathcal{K}_x$$

where

$$(4.54) \qquad G_x \equiv C_x(\cdot)_{xx} - (\mathrm{diag}(e^{(F(x)-F(y))/\varepsilon})_{y \leq x} \vec{g}_x)^t \cdot \mathcal{K}_x^{-1} \vec{g}_x.$$

This follows by simply adding the column vector

$$\begin{pmatrix} \mathcal{K}_x \\ -(\mathrm{diag}(e^{(F(x)-F(y))/\varepsilon})_{y \leq x} \vec{g}_x)^t \end{pmatrix} \mathcal{K}_x^{-1} \vec{g}_x$$

(which clearly is a linear combination of the first columns of $C_x$) to the last column in $C_x$. From this representation we obtain that $\lambda_x$ is very close to $\lambda_{\mathcal{M}_{<x}}$. We begin with:

LEMMA 4.9. *There is $N \equiv N(d)$ such that for all $\alpha > 0$ and some $C \equiv C(d, F, \alpha)$ dominating the supremum norms of the Landau symbols appearing below and all*

$$(4.55) \qquad\qquad e^{\alpha/\varepsilon}/\mathcal{E}_x < \lambda < e^{-\alpha/\varepsilon}/T_{\mathcal{M}_{\leq x}},$$

*the inverse of $\mathcal{K}_x(\lambda)$ exists. More precisely, uniformly in $\lambda$*

$$(4.56) \qquad (\mathcal{K}_x(\lambda)^{-1} \vec{g}_x(\lambda))_y = \mathcal{O}(1) \frac{1}{\varepsilon^N \lambda T_{y,x}}, \qquad y \leq x.$$

*Moreover, we obtain*

$$(4.57) \qquad\qquad \lambda \in \sigma(L_\varepsilon) \iff G_x(\lambda) = 0,$$

*where $G_x(\lambda)$ is defined in (4.54).*

PROOF. Formally the inverse $\mathcal{K}_x(\lambda)^{-1}$ is given by the Neumann series, that is,

$$(4.58) \quad \begin{aligned} \mathcal{K}_x(\lambda)^{-1} &= (\mathbb{1} - \mathcal{D}_x(\lambda)^{-1} \mathcal{N}_x(\lambda))^{-1} \mathcal{D}_x(\lambda)^{-1} \\ &= \sum_{s=0}^{\infty} (\mathcal{D}_x(\lambda)^{-1} \mathcal{N}_x(\lambda))^s \mathcal{D}_x(\lambda)^{-1}. \end{aligned}$$

In order to make sense out of this calculation and to extract the exponential decay estimate written in (4.56) out of the sum, a straightforward computation for $s \in \mathbb{N} \setminus 0$ gives the random walk representation

$$(4.59) \quad \begin{aligned} &(\mathcal{D}_x(\lambda)^{-1} \mathcal{N}_x(\lambda))^s \mathcal{D}_x(\lambda)^{-1} \\ &= \left( \sum_{\substack{\omega : y \to z \\ |\omega| = s}} \prod_{t=1}^{|\omega|} \frac{C_x(\lambda)_{\omega_{t-1} \omega_t}}{C_x(\lambda)_{\omega_{t-1} \omega_{t-1}}} \frac{1}{C_x(\lambda)_{zz}} \right)_{yz, y, z < x}, \end{aligned}$$

where for $J \subset \mathcal{M}_{<x}$ we write shorthand $\omega\colon y \to J$ for a sequence $\omega = (\omega_0, \ldots, \omega_T)$ such that $\omega_0 = y$, $\omega_T \in J$, $\omega_t \in \mathcal{M}_{\leq x}\setminus J$ and $\omega_{t-1} \neq \omega_t$ for all $t = 1, \ldots, T$. $|\omega|$ denotes the length $T$ of the sequence. By means of (4.46) we may apply (4.25) for $x \equiv y < x$—in slight abuse of notation—and $\tilde{I} \equiv \mathcal{M}_{\leq x}\setminus y$ and conclude using (4.5) in combination with Proposition 4.7 in [4] that for some $C \equiv C(d, F)$ and $N \equiv N(d)$ and all $\lambda$ satisfying (4.55)

$$
\begin{aligned}
(4.60)\quad C_x(\lambda)_{yy} &\geq (1/C) \int_{A_{x,\mathcal{M}_{<x}}} e^{-(F - F(y))/\varepsilon}\, dy \\
&\quad \times (\lambda - \lambda_{\mathcal{M}_{\leq x}\setminus y})(1 + (\lambda - \lambda_{\mathcal{M}_{\leq x}\setminus y})\mathcal{O}(1)T_{\mathcal{M}_{\leq x}}) \\
&= (1/C)\lambda \varepsilon^N.
\end{aligned}
$$

In addition, (3.26) for $I \equiv \mathcal{M}_{<x}$ in combination with the upper bound in (4.55) proves for some $C \equiv C(d, F)$ and all $y, z \in \mathcal{M}_{\leq x}$, $y \neq z$,

$$
(4.61)\qquad \partial_n h_z^\lambda \geq C\,\partial_n h_z^0 \geq C\,\partial_n h_{z,y}^0 = -C\,\partial_n h_{y,z}^0
$$

on $\partial B_y$. Harnack's inequality applied to $\phi_x|B(y,\varepsilon)$ and Corollary 4.8 in [4] show for some $N \equiv N(d)$

$$
\begin{aligned}
(4.62)\quad -C_x(\lambda)_{yz} &= -\varepsilon \int_{\partial B_y} e^{-(F - F(y))/\varepsilon} \frac{\phi_x}{\phi_x(y)}\, \partial_n h_z^\lambda\, d\sigma \\
&\leq C\varepsilon \int_{\partial B_y} e^{-(F - F(y))/\varepsilon}\, \partial_n h_{y,z}^0\, d\sigma \\
&\leq C\varepsilon^{-N}/T_{y,z}.
\end{aligned}
$$

Now fix $y \in \mathcal{M}_{\leq x}$, $z \in \mathcal{M}_{\leq x}\setminus y$ and a sequence $\omega \equiv (\omega_0, \ldots, \omega_T)\colon z \to y$ such that $\omega_{t+1} \neq \omega_t$ for $t = 0, \ldots, T-1$, $\omega_t \in \mathcal{M}_{<x}$ for $t = 1, \ldots, T-1$. We observe that there is $1 \leq t_0 \leq T$ such that

$$
(4.63)\qquad \hat{F}(\omega_{t_0-1}, \omega_{t_0}) \geq \hat{F}(z, y).
$$

For, if we assume that the contrary is true for all $t_0 = 1, \ldots, T-1$, it follows from ultrametricity of $\mathcal{M}_{<x} \times \mathcal{M}_{<x} \ni (y, z) \mapsto \hat{F}(y, z)$ that $\hat{F}(\omega_1, y) = \hat{F}(z, y)$. By the same argument this implies $\hat{F}(\omega_2, y) = \hat{F}(z, y)$ and so forth. We conclude $\hat{F}(\omega_{T-1}, \omega_T) = \hat{F}(\omega_{T-1}, y) = \hat{F}(z, y)$ so that $t_0 \equiv T$ does the job. Choose $t_0$ satisfying (4.63). Combining this with the triviality $T_{y,z} \geq \mathcal{E}_x$ for $y, z \in \mathcal{M}_{\leq x}$, $y \neq z$, $y < x$, where $\mathcal{E}_x$ is defined in (4.42), we obtain

$$
\begin{aligned}
(4.64)\quad \prod_{t=1}^{|\omega|} T_{\omega_{t-1}, \omega_t} &= \prod_{t=1}^{t_0 - 1} T_{\omega_t, \omega_{t-1}} e^{\hat{F}(\omega_{t_0-1}, \omega_{t_0}) - F(z)/\varepsilon} \prod_{t=t_0+1}^{|\omega|} T_{\omega_{t-1}, \omega_t} \\
&\geq T_{z,y}\mathcal{E}_x^{|\omega|-1}.
\end{aligned}
$$

Combination of (4.64), (4.62) and (4.60) tells us that the computations in (4.58) are justified and that for $y < x$

$$
\begin{aligned}
(\mathcal{K}_x(\lambda)^{-1}\vec{g}_x(\lambda))_y &= \sum_{\omega:y\to x} \prod_{t=1}^{|\omega|} \frac{C_x(\lambda)_{\omega_t\omega_{t-1}}}{C_x(\lambda)_{\omega_{t-1}\omega_{t-1}}} \\
&= \sum_{\omega:y\to x} \mathcal{O}(1)\frac{1}{\varepsilon^N \lambda T_{y,x}}\left(\mathcal{O}(1)\frac{1}{\varepsilon^N \lambda \mathcal{E}_x}\right)^{|\omega|-1} \\
&= \mathcal{O}(1)\frac{1}{\varepsilon^N \lambda T_{y,x}} \sum_{t=1}^{\infty} (|\mathcal{M}_{<x}|-1)^{t-1} \varepsilon^{-N(t-1)} e^{-\rho(t-1)/\varepsilon} \\
&= \mathcal{O}(1)\frac{1}{\varepsilon^N \lambda T_{y,z}}.
\end{aligned}
\tag{4.65}
$$

We thus obtain (4.56).

Equation (4.57) then is a direct consequence of (4.53) and Lemma 4.1 for $\tilde{I} \equiv \mathcal{M}_{\leq x}$. $\square$

We are searching for solutions $\lambda$ near $\lambda_{\mathcal{M}_{<x}}$ of the equation appearing in (4.57). We want to apply Lagrange's theorem to this equation (see [42]) which tells us the following: Fix a point $a \in \mathbb{C}$ and an analytic function $\Psi$ defined on a domain containing the point $a$. Assume that there is a contour in the domain surrounding $a$ such that on this contour the estimate $|\Psi(\zeta)| < |\zeta - a|$ holds. Then the equation

$$
\zeta = a + \Psi(\zeta) \tag{4.66}
$$

has a unique solution in the interior of the contour. Furthermore, the solution can be expanded in the form

$$
\zeta = a + \sum_{n=1}^{\infty} (n!)^{-1} \partial_\zeta^{n-1} \Psi(a)^n. \tag{4.67}
$$

We are in a position to prove Theorem 4.6.

PROOF OF THEOREM 4.6. Equation (4.57) can be written as

$$
-C_x(\lambda)_{xx} + \Phi_j(\zeta) = 0, \tag{4.68}
$$

where we have set $\zeta \equiv \lambda \int_{A^\beta_{x,\mathcal{M}_{<x}}} e^{-F/\varepsilon} du / \operatorname{cap}^0(x, \mathcal{M}_{<x})$ and

$$
\Phi_j(\zeta) \equiv \sum_{y<x} C_x(\lambda)_{xy} (\mathcal{K}_x(\lambda)^{-1} \vec{g}_x(\lambda))_y. \tag{4.69}
$$

Fix constant arbitrary $\alpha > 0$ and let us denote by $U_x$ the interval of all $\zeta$ such that

$$
e^{\alpha/\varepsilon} T_x / \mathcal{E}_x < \zeta < e^{-\alpha/\varepsilon} T_x / T_{\mathcal{M}_{\leq x}}. \tag{4.70}
$$

Defining $\zeta_{\mathcal{M}_{<x}} \equiv \lambda_{\mathcal{M}_{<x}} \operatorname{cap}^0(x, \mathcal{M}_{<x})/\int_{A^\beta_{x,\mathcal{M}_{<x}}} e^{-F/\varepsilon} \, du$, $\beta > 0$ small, it follows $\zeta_{\mathcal{M}_{<x}} = e^{\mathcal{O}(1)}$ from Theorem 4.2 and we may apply (4.25) for all $\zeta \in U_x$ to obtain for some $N \equiv N(d)$

$$(4.71) \quad -C_x(\lambda)_{xx} = \frac{\int_{A^\beta_{x,\mathcal{M}_{<x}}} e^{-F/\varepsilon} \, du}{\operatorname{cap}^0(x, \mathcal{M}_{<x})} \times e^{\mathcal{O}(1)} \varepsilon^{-N} (\zeta - \zeta_{\mathcal{M}_{<x}} + (\zeta - \zeta_{\mathcal{M}_{<x}})^2 R_x(\zeta)),$$

where

$$(4.72) \quad R_x(\zeta) = \mathcal{O}(1) e^{\rho/\varepsilon} T_{\mathcal{M}_{\leq x}} \frac{\operatorname{cap}^0(x, \mathcal{M}_{<x})}{\int_{A^\beta_{x,\mathcal{M}_{<x}}} e^{-F/\varepsilon} \, du} = \mathcal{O}(1) \varepsilon^{-N} e^{\rho/\varepsilon} \frac{T_{\mathcal{M}_{\leq x}}}{T_x}$$

by Proposition 4.7 in [4]. In view of (4.71) it follows that (4.68) is equivalent to

$$(4.73) \quad \zeta = \zeta_{\mathcal{M}_{<x}} + \Psi_x(\zeta)$$

for some function $\Psi_x$ satisfying

$$(4.74) \quad \Psi_x(\zeta) = \frac{\operatorname{cap}^0(x, \mathcal{M}_{<x})}{\int_{A^\beta_{x,\mathcal{M}_{<x}}} e^{-F/\varepsilon} \, du} \varepsilon^{-N} e^{\mathcal{O}(1)} \Phi_x(\zeta) + (\zeta - \zeta_{\mathcal{M}_{<x}})^2 R_x(\zeta).$$

Furthermore, (4.62) shows for all $\zeta \in U_x$ and all $y < x$

$$(4.75) \quad -C_x(\lambda)_{xy} = \mathcal{O}(1) \frac{1}{T_{x,y}}.$$

Thus for all $|\zeta - \zeta_{\mathcal{M}_{<x}}| = 1$ we deduce from (4.5) and (4.56) that for all $\zeta \in U_x$

$$(4.76) \quad \frac{\operatorname{cap}^0(x, \mathcal{M}_{<x})}{\int_{A^\beta_{x,\mathcal{M}_{<x}}} e^{-F/\varepsilon} \, du} \varepsilon^{-N} |\Phi_x(\zeta)| \leq \sum_{y<x} \frac{\varepsilon^{-N} e^{\rho/\varepsilon} T_x^2}{T_{x,y} T_{y,x}} \leq \frac{\varepsilon^{-N} e^{\rho/\varepsilon} T_x}{\mathcal{T}_x}.$$

By means of (4.72) and (4.76) it follows for $|\zeta - \zeta_{\Sigma_{j-1}}| = 1$

$$(4.77) \quad |\Psi_x(\zeta)| \leq \frac{\varepsilon^{-N} e^{\rho/\varepsilon} T_x}{\mathcal{T}_x} + \frac{\varepsilon^{-N} e^{\rho/\varepsilon} T_{\mathcal{M}_{\leq x}}}{T_x}.$$

Since $\mathcal{T}_x \geq \mathcal{E}_x$, in view of (4.44) we may apply Lagrange's theorem to (4.73) giving the existence of a unique solution $\zeta_x = \lambda_x \int_{A^\beta_{x,\mathcal{M}_{<x}}} e^{-F/\varepsilon} \, du / \operatorname{cap}^0(x, \mathcal{M}_{<x})$ of (4.68) satisfying $|\tilde{\zeta}_j - \zeta_{\mathcal{M}_{<x}}| < 1$. We rewrite (4.73) in the form

$$(4.78) \quad \zeta_x = \zeta_{\mathcal{M}_{<x}} + \mathcal{O}(1) \left( \frac{\varepsilon^{-N} e^{\rho/\varepsilon} T_x}{\mathcal{T}_x} + \frac{\varepsilon^{-N} e^{\rho/\varepsilon} T_{\mathcal{M}_{\leq x}}}{T_x} \right).$$

Since from invertibility of $\mathcal{K}_x(\lambda_x)$ it follows that the kernel of $C_x(\lambda_x)$ is at most one-dimensional, Lemma 4.1 implies that $\lambda_x$ is simple. Using (4.5) for $\tilde{I} \equiv \mathcal{M}_{<x}$ from (4.78) we derive that (4.35) holds. Moreover, using $\lambda_x < \lambda_{\mathcal{M}_{\leq x}}$, which follows from (4.5) and (2.5) in combination with (2.33) from Lemma 4.1 we conclude that

$$(4.79) \qquad (\phi_x(y))_{y<x} = \phi_x(x) \mathcal{K}_x(\lambda_x)^{-1} \vec{g}_x(\lambda_x).$$

Hence from (4.56) and $\lambda_x = e^{\mathcal{O}(1)} \lambda_{\mathcal{M}_{<x}}$ we obtain from (3.38) and (3.39) that (4.36) is satisfied. Now it is very easy to finish the theorem. In view of Lemma 4.8 and Assumption 4.5 choosing $\alpha < \liminf_{\varepsilon \downarrow 0} \rho$ the union of the intervals described in (4.55) contains an interval of the form $[0, \varepsilon^N)$, $N \equiv N(d)$. Noting that [39] after possibly increasing $N$ gives the existence of $|\mathcal{M}|$ eigenvalues in this interval, we obtain (4.34). Actually, similarly to (2.7) one can also use the variational formula, Theorems 4.5.1 and 4.5.2 in [12], by (4.36) obvious choice $\tilde{\phi}_x \equiv \mathbb{1}^\beta_{A_{x,\mathcal{M}\setminus x}}$, $x \in \mathcal{M}$, for some small $\beta > 0$ as a trial function for $\phi_x$ to obtain the existence of these eigenvalues—and already rather precise upper bounds. $\square$

**5. Distribution of metastable transition times.** In the sequel we show how the structure of the low-lying part of the spectrum developed in the previous section determines in a precise way the asymptotic behavior of the distribution of metastable transition times $\tau(x)$ defined in (1.16).

We first would like to point out that Theorem 4.6 holds in more generality with only obvious changes in the notation. Namely, if we define $\Omega \setminus B_I$ instead of $\Omega$, for some nonempty, properly contained subset $I \subset \mathcal{M}$, then Theorem 4.6 still holds for the Dirichlet realization $L_\varepsilon^I \equiv L_\varepsilon^{\Omega \setminus \overline{B}_I}$ in $L^2(\Omega \setminus \overline{B}_I, e^{-F/\varepsilon} dx)$ when $\mathcal{M}(x)$ is replaced by $\mathcal{M}(x) \cup I$ for all $x \in \mathcal{M} \setminus I$. Moreover, we could have looked only for the principal eigenvalue and its eigenfunctions of $L_\varepsilon^I$ and the same procedure then leads to:

THEOREM 5.1. *Fix a nonempty, properly contained subset $I \subset \mathcal{M}$. Then Theorem 4.6 still holds for the operator $L_\varepsilon^I$ with the modification that $\mathcal{M}$ has to be replaced by $\mathcal{M} \setminus I$, $\mathcal{M}(y)$, $y \in \mathcal{M} \setminus I$, by its union with $I$, $\mathcal{M}_{<y}$, $\mathcal{M}_{\leq y}$ are defined with respect to the time scales $T_y \equiv T_{y, \mathcal{M}(y) \cup I}$ and $\tilde{\mathcal{M}}_{<y}$, $\tilde{\mathcal{M}}_{\leq y}$ are chosen depending on the exchange of $\mathcal{M}_{<y}$, $\mathcal{M}_{\leq x}$.*

*Let $x \in \mathcal{M} \setminus I$ be the unique solution to the equation $T_{x,I} = T_I$, where the latter was defined in (2.31). In addition to the equivalent of (4.36) in this situation we have on $\{F < C_1\}$*

$$(5.1) \qquad \phi_x = \left(1 + \mathcal{O}(1) \varepsilon^{-N} \frac{T_{I \cup x}}{T_I}\right) \phi_x(x) h^0_{x,I}$$

*with the usual dependence of the constant $N$ and the constant dominating the Landau symbol and where we use the convention $T_\mathcal{M} \equiv 1$.*

For $I \equiv \mathcal{M}(x) \neq \varnothing$ modulo factors of order $1 + \mathcal{O}(1)\varepsilon^{-N}(T_{\mathcal{M}_{\leq y}}/T_{\mathcal{M}_{<x}})$, respectively, the small eigenvalues of $L_\varepsilon^I$ equal $\lambda_y$, $y \in \mathcal{M}\setminus\mathcal{M}_{<x}$, where $\lambda_x$ are given by (4.34) for $I \equiv \varnothing$.

As we shall see, it is not difficult to obtain the leading part in the following result from Theorem 5.1. But for reasonable control of the remainder term we have to prove additional a priori large deviation type estimates to which most of this section is devoted. Recall the definition (1.2) of the diffusion generated by $L_\varepsilon$ starting in $x$ and that of the hitting time $\tau_B^x$ in (1.14) and set $\tau_{I\cup\Omega^c}^x \equiv \tau_{B_I\cup\Omega^c}^x$ in slight abuse of notation.

THEOREM 5.2. *Suppose the assumptions of Theorem 4.6 are met and assume that either $\Omega$ is bounded or $\int_{\{F>C_1\}} |\nabla F|^d e^{-(F-C_1)/\gamma}\, dy < \infty$ for some $\gamma > 0$. Assume, moreover, that $\mu_\varepsilon$ defined in Assumption 1.2 satisfies $\mu_\varepsilon \geq \delta\varepsilon$ for some constant $\delta > 0$ independent in $\varepsilon > 0$. Let $x \in \mathcal{M}\setminus I$ be the unique local minimum such that $T_{x,I} = T_I$. Let $\lambda_x$ be the principal eigenvalue of $L_\varepsilon^I$. Then for all $t \geq 0$*

$$(5.2) \qquad \mathbb{P}[\tau_{I\cup\Omega^c}^x > t] = e^{-\lambda_x t}\left(1 + \mathcal{O}(1)\varepsilon^{-N}\frac{T_{I\cup x}}{T_I}\right),$$

*where the modulus of each Landau symbol is bounded by a constant $C \equiv C(d, F)$ and $N \equiv N(d) > 0$.*

For $I \equiv \mathcal{M}(x)$, $x \in \mathcal{M}$, we obtain Theorem 1.4 from (4.20), the generalization of (4.35) and (4.41).

REMARK. We note that the a priori bound on $\mu_\varepsilon$ is the natural choice. More precisely, for $F$ of sufficient regularity this is the case as can be proven by a semiclassical approximation similarly to Theorem 11.1 in [8]. Moreover, this property is trivially fulfilled if $\liminf_{|x|\to\infty} |\Delta F(x)|/|\nabla F(x)|^2 < \infty$.

REMARK. As the reader might observe while reading the proof of (5.2), the methods are fully sufficient to obtain the following result from the generalization of Theorem 4.6 described in the first part of Theorem 5.1. In the situation of the previous theorem let $\sigma(L_\varepsilon^I) = \{\lambda_y | y \in \mathcal{M}\setminus I\}$ be the low-lying spectrum of $L_\varepsilon^I$. Then for some $\kappa \equiv \kappa(d, F) > 0$ and all $t \geq 0$

$$(5.3) \qquad \begin{aligned} \mathbb{P}[\tau_{I\cup\Omega^c}^x > t] &= e^{-\lambda_x t}\left(1 + \mathcal{O}(1)\varepsilon^{-N}\frac{T_{I\cup x}}{T_I}\right) \\ &\quad + \sum_{y\in\mathcal{M}\setminus I\setminus x} e^{-\lambda_y t}\mathcal{O}(1)\varepsilon^{-N}\frac{T_{\mathcal{M}_{<y}}}{T_I} \\ &\quad + e^{-\varepsilon^\kappa t}\mathcal{O}(1)\varepsilon^{-N}\left(\frac{1}{T_I} + e^{-(C_1-F(x))/\varepsilon}\right) \end{aligned}$$

with the usual dependence of the errors. Since the computations [starting with (5.4)] necessary for this result are a bit tedious, to keep this work at a reasonable length we omit its proof. Instead we refer the interested reader to [3] for a proof in discrete space, which unfortunately does not generalize directly to the continuous state space setting. The full strength of this expansion would be achieved if one proves lower bounds on eigenfunctions, that is, in view of (3.38) on transition probabilities $h_{A,B}^0(z)$, in regions where they are small by, for example, applying large deviation principles. This would lead to a replacement of $\mathcal{O}(1)$ in the sum above by $\exp(-(\text{rate} + o(1))/\varepsilon)$, where the rate depends on the properties of the flow of $\nabla F$.

PROOF OF THEOREM 5.2. Assumption 4.5 assures that $x$ is unique. Let $\phi_x$ be a corresponding eigenfunction to the principal eigenvalue $\lambda_x$ of $L_\varepsilon^I$ normalized such that $\phi_x(x) = 1$. Fix $y \in B_x$ and write

$$(5.4) \qquad \mathbb{P}[\tau_{I \cup \Omega^c}^y > t] = e^{-tL_\varepsilon^I}(\mathbb{1}_{\Omega \setminus B_I})(y).$$

Using the spectral decomposition corresponding to the principal eigenvalue of $L_\varepsilon^I$, we compute

$$(5.5) \quad \begin{aligned} & e^{(F(y)-F(x))/\varepsilon} e^{-tL_\varepsilon^I}(\mathbb{1}_{\Omega \setminus \overline{B}_I})(y) \\ & = e^{-\lambda_x t} \frac{(\phi_x, \mathbb{1}_{\Omega \setminus \overline{B}_I})_{F-F(x)}}{\|\phi_x\|_{F-F(x)}^2} \phi_x(y) + e^{-tL_\varepsilon^I} \Pi_x(\mathbb{1}_{\Omega \setminus \overline{B}_I})(y), \end{aligned}$$

where $\|\cdot\|_F$ denotes the norm induced by the inner product $(\cdot,\cdot)_F$ on $L^2(\Omega \setminus B_I, e^{-F/\varepsilon} dz)$ and where $\Pi_x$ is the orthogonal projection onto $(\phi_x)^\perp$. To estimate the second term on the right-hand side of (5.5) we introduce $\Omega_1 \equiv \{F < C_1\}$ and $\Gamma_k \equiv \{-\varepsilon \leq F - C_1 - k\varepsilon < 0\}$ for $k \in \mathbb{N}$ and write

$$(5.6) \quad \begin{aligned} & e^{-tL_\varepsilon^I} \Pi_x(\mathbb{1}_{\Omega \setminus \overline{B}_I})(y) \\ & = \mathbb{E}[\Pi_x(\mathbb{1}_{\Omega \setminus \overline{B}_I})(X_t^y), \tau_{I \cup \Omega^c}^y > t] \\ & = \mathbb{E}[\Pi_x(\mathbb{1}_{\Omega \setminus \overline{B}_I})(X_t^y), X_t^y \in A_{x,I}(\beta), \tau_{I \cup \Omega^c}^y > t] \\ & \quad + \mathbb{E}[\Pi_x(\mathbb{1}_{\Omega \setminus \overline{B}_I})(X_t^y), X_t^y \in \Omega_1 \setminus A_{x,I}(\beta), \tau_{I \cup \Omega^c}^y > t] \\ & \quad + \sum_{k \in \mathbb{N}} \mathbb{E}[\Pi_x(\mathbb{1}_{\Omega \setminus \overline{B}_I})(X_t^y), X_t^y \in \Gamma_k, \tau_{I \cup \Omega^c}^y > t], \end{aligned}$$

where we recall definition (3.37) of $A_{x,I}^\beta$. To estimate the right-hand side we start with the claim that

$$(5.7) \qquad \mathbb{1}_{\Omega_1 \setminus \overline{B}_I} \Pi_x(\mathbb{1}_{\Omega \setminus \overline{B}_I}) = \mathcal{O}(1) \varepsilon^{-N} \frac{T_{I \cup x}}{T_I} \mathbb{1}_{A_{x,I}^\beta} + \mathcal{O}(1) \mathbb{1}_{\Omega_1 \setminus \overline{A_{x,I}^\beta} \setminus \overline{B}_I}.$$

Combination of (5.1) with Corollary 4.7 in [3] gives on $A_{x,I}^\beta$, $\beta > C\varepsilon \log(1/\varepsilon)$,

$$(5.8) \quad \begin{aligned} \phi_x &= \sum_{y \in I} \mathbb{1}_{A_{y,x}^\beta \setminus \overline{B}_I} \mathcal{O}(1) e^{-(\hat{F}(y,x)-F)/\varepsilon}/\varepsilon + \mathcal{O}(1) \mathbb{1}_{\Omega_1 \setminus \overline{A_{I,x}^\beta} \setminus \overline{A_{x,I}^\beta}} \\ & \quad + \left(1 + \mathcal{O}(1)\varepsilon^{-N}\left(\frac{T_{I \cup x}}{T_I} + e^{-(\hat{F}(x,I)-F)/\varepsilon}/\varepsilon\right)\right) \mathbb{1}_{A_{x,I}^\beta}, \end{aligned}$$

where $A^\beta_{I,x} \equiv \bigcup_{y\in I} A^\beta_{y,x}$ and thus

$$\|\phi_x\|^2_{F-F(x)} = \left(1 + \mathcal{O}(1)\varepsilon^{-N}\left(\frac{T_{I\cup x}}{T_I} + e^{-2(\hat{F}(x,I)-F(x)-\beta)/\varepsilon}\right)\right) \int_{A^\beta_{x,I}} e^{-(F-F(x))/\varepsilon} \tag{5.9}$$

whereas

$$(\phi_x, \mathbb{1}_{\Omega\setminus\overline{B}_I})_{F-F(x)} = \left(1 + \mathcal{O}(1)\varepsilon^{-N}\left(\frac{T_{I\cup x}}{T_I} + e^{-(\hat{F}(x,I)-F(x)-\beta)/\varepsilon}\right)\right) \\ \times \int_{A^\beta_{x,I}} e^{-(F-F(x))/\varepsilon}. \tag{5.10}$$

Equations (5.9) and (5.10) give the claim and thus the first term on the right-hand side of (5.6) is of order $\varepsilon^{-N}(T_{I\cup x}/T_I)\mathbb{P}[\tau^y_{I\cup\Omega^c} > t]$ for some $N \equiv N(d, F)$. In addition, together with (5.8) they imply that the first term on the right-hand side of (5.5) equals the leading part in (5.2).

The remaining part of the proof is devoted to the estimate of the second and the third terms on the right-hand side of (5.6). Next, we generalize Proposition 5.9 in [2] to our setting. This generalization is—besides technical details—straightforward. We have:

LEMMA 5.3. *For some $N \equiv N(d, F)$ uniformly in $\varepsilon > 0$ and $t > 0$,*

$$\mathbb{P}[\tau^x_I > t] \geq e^{-\varepsilon^{-N}t/T_I}\varepsilon^N. \tag{5.11}$$

For the proof of this lemma it will be useful to introduce the following renewal structure. Let $W$, $W^{(k)}$, $k \in \mathbb{N} \cup 0$, be a sequence of independent Brownian motions on $\mathbb{R}^d$ starting in zero defined on a common probability space; denote by $X^z$, $X^{z,(l)}$, the strong solution to (1.2) defined with respect to $W$, $W \equiv W^{(l)}$, respectively, starting in $x \equiv z$. Fix two regular domains $A \subset \Sigma \subset B^c$ and define the stopping times

$$\sigma^z_l \equiv \inf\{t \geq 0 | \exists_{0<s<t} X^{z,(l)}_s \notin \Sigma, X^{z,(l)}_t \in A\}, \\ \rho^z_l \equiv \inf\{t \geq 0 | X^{z,(l)}_t \in B\}. \tag{5.12}$$

Moreover, for $y \in A$ set $S_0 \equiv 0$, $z_0 \equiv y$, $S_l \equiv S_{l-1} + \sigma^{z_{l-1}}_l$, $z_l \equiv X^{z_{l-1},(l)}_t$ for $t \equiv S_l$.

PROOF OF LEMMA 5.3. Now choose $B \equiv B_I \cup \Omega^c$, $A \equiv B_x$ and $\Sigma \equiv B_x(\varepsilon/2)$. We may write, using the strong Markov property and independence, for every $k \in \mathbb{N}$

$$\begin{aligned}\mathbb{P}[\tau^x_I \geq t] &\geq \mathbb{P}[S_k > t, \forall_{1\leq l\leq k}\sigma^{z_{l-1}}_l < \rho^{z_{l-1}}_l] \\ &= \mathbb{P}[\forall_{1\leq l\leq k}\sigma^{z_{l-1}}_l < \rho^{z_{l-1}}_l]\mathbb{P}[S_k > t | \forall_{1\leq l\leq k}\sigma^{z_{l-1}}_l < \rho^{z_{l-1}}_l] \\ &\geq \left(\inf_{z\in\partial B_x}\mathbb{P}[\sigma^z_1 < \rho^z_1]\right)^k \mathbb{P}[S_k > t | \forall_{1\leq l\leq k}\sigma^{z_{l-1}}_l < \rho^{z_{l-1}}_l] \\ &\geq (1 - Ce^{-(\hat{F}(z,I)-F(x))/\varepsilon}/\varepsilon)^k \mathbb{P}[S_k > t | \forall_{1\leq l\leq k}\sigma^{z_{l-1}}_l < \rho^{z_{l-1}}_l],\end{aligned} \tag{5.13}$$

where the last inequality follows by the strong Markov property and Corollary 4.8 in [4] yielding

$$(5.14) \quad \mathbb{P}[\sigma_1^z < \rho_1^z] \geq \inf_{y \in \partial B(x,\varepsilon/2)} \mathbb{P}[\tau_x^y < \tau_{I \cup \Omega^c}^y] \geq 1 - Ce^{-(\hat{F}(z,I) - F(x))/\varepsilon}/\varepsilon.$$

The claim is proven once we show that for some $N$ and all $t > 0$ and $k \equiv \varepsilon^{-N} t$ the second term in the product of (5.13) is bounded below. For, we recall an inequality going back to Paley and Zygmund—also referred to as the second moment inequality—saying that

$$(5.15) \quad P[X > (1-\delta)\mathbb{E}[X]] \geq \delta^2 E[X]^2 / E[X^2], \qquad \delta \in (0,1),$$

for any random variable $X$ with finite expectation. We want to apply this inequality to the variable $X \equiv S_k/k$, where we choose $\delta \equiv 1 - l/(Rk/\inf_{z \in \partial B_x} \mathbb{E}[\sigma_1^z | \sigma_1^z < \rho_1^z])$, where $R$ is specified later on, and the probability measure $P \equiv \mathbb{P}[\cdot | \forall_{1 \leq l \leq k} \sigma_l^{z_{l-1}} < \rho_l^{z_{l-1}}]$. Therefore, we notice that by independence

$$(5.16) \quad \begin{aligned} &\mathbb{E}[S_k | \forall_{1 \leq l \leq k} \sigma_l^{z_{l-1}} < \rho_l^{z_{l-1}}] \\ &= \sum_{n \leq k} \mathbb{E}[\sigma_n^{z_{n-1}} | \forall_{1 \leq l \leq k} \sigma_l^{z_{l-1}} < \rho_l^{z_{l-1}}] \\ &\geq \sum_{n \leq k} \mathbb{E}[\mathbb{E}[\sigma_1^y, \sigma_1^y < \rho_1^y]|_{y \equiv z_{n-1}}, \forall_{1 \leq l < n} \sigma_l^{z_{l-1}} < \rho_l^{z_{l-1}}] \\ &\qquad \times \frac{\inf_{z_0 \in \partial B_x} \mathbb{P}[\forall_{1 \leq l \leq k-n} \sigma_l^{z_{l-1}} < \rho_l^{z_{l-1}}]}{\mathbb{P}[\forall_{1 \leq l \leq k} \sigma_l^{z_{l-1}} < \rho_l^{z_{l-1}}]} \\ &\geq (k/C) \inf_{z \in \partial B_x} \mathbb{E}[\sigma_1^z | \sigma_1^z < \rho_1^z]. \end{aligned}$$

In the last inequality we have used independence once more in combination with the Harnack inequality for harmonic measures (see, e.g., Theorem 4.3 in [36] which is applicable after a scaling argument to get rid of the dependence on $\varepsilon$) saying that for $Z^{y,(1)} \equiv X_t^{y,(1)}$ at $t \equiv \tau_{B_x(\varepsilon/2)^c}^y$ and $Z^y \equiv X_t^y$ at $t \equiv \tau_x^y$ again by independence for some $C \equiv C(d, F)$ and all $y, \tilde{y} \in \partial B_x$,

$$(5.17) \quad \begin{aligned} &\mathbb{P}[\sigma_1^y < \rho_1^y, \forall_{2 \leq l \leq n} \sigma_l^{z_{l-1}} < \rho_l^{z_{l-1}}] \\ &= \int_{\partial B_x(\varepsilon/2)} \mathbb{P}[Z^{y,(1)} \in dz] \mathbb{P}[\tau_x^z < \tau_{I \cup \Omega^c}^z, Z^z \in d\tilde{z}] \\ &\qquad \times \mathbb{P}[\forall_{1 \leq l < n} \sigma_l^{z_{l-1}} < \rho_l^{z_{l-1}}] \\ &= e^{\mathcal{O}(1)} \mathbb{P}[\sigma_1^{\tilde{y}} < \rho_1^{\tilde{y}}, \forall_{2 \leq l \leq n} \sigma_l^{z_{l-1}} < \rho_l^{z_{l-1}}]. \end{aligned}$$

Similarly, one proves

$$
\begin{aligned}
\mathbb{E}[(S_k)^2|\forall_{1\leq l\leq k}\sigma_l^{z_{l-1}} &< \rho_l^{z_{l-1}}] \\
&= \sum_{n\leq k} \mathbb{E}[(\sigma_n^{z_{n-1}})^2|\forall_{1\leq l\leq k}\sigma_l^{z_{l-1}} < \rho_l^{z_{l-1}}] \\
&\quad + \sum_{\substack{n,m\leq k \\ m\neq n}} \mathbb{E}[\sigma_n^{z_{n-1}}\sigma_m^{z_{m-1}}|\forall_{1\leq l\leq k}\sigma_l^{z_{l-1}} < \rho_l^{z_{l-1}}] \\
&= \mathcal{O}(1)k \sup_{z\in\partial B_x} \mathbb{E}[(\sigma_1^z)^2|\sigma_1^z < \rho_1^z] \\
&\quad + \mathcal{O}(1)k(k-1) \sup_{z\in\partial B_x} \mathbb{E}[\sigma_1^z|\sigma_1^z < \rho_1^z]^2.
\end{aligned}
\tag{5.18}
$$

We thus obtain for all $z \in \partial B_x$

$$
\begin{aligned}
\mathbb{P}[S_k > t|\forall_{1\leq l\leq k}\rho_l^{z_{l-1}} &< \sigma_l^{z_{l-1}}] \\
&\geq (1/C)\left(1 - t\Big/\left(Rk \inf_{z\in\partial B_x} \mathbb{E}[\sigma_1^z|\sigma_1^z < \rho_1^z]\right)\right)^2 \\
&\quad \times k^2 \inf_{z\in\partial B_x} \mathbb{E}[\sigma_1^z|\sigma_1^z < \rho_1^z]^2 \\
&\quad \times \left(k(k-1)\sup_{z\in\partial B_x}\mathbb{E}[\sigma_1^z|\sigma_1^z<\rho_1^z]^2 + k\sup_{z\in\partial B_x}\mathbb{E}[(\sigma_1^z)^2|\sigma_1^z<\rho_1^z]\right)^{-1}.
\end{aligned}
\tag{5.19}
$$

It remains to estimate the right-hand side of the previous inequality from below. By the strong Markov property we compute

$$
\begin{aligned}
\mathbb{E}[(\sigma_1^z)^2, \sigma_1^z < \rho_1^z] &= \mathbb{E}[(\tau_{B_x(\varepsilon/2)^c}^z)^2 \mathbb{P}[\tau_x^y < \tau_{I\cup\Omega^c}^y y]_{|y\equiv X_t^z, t\equiv \tau_{B_x(\varepsilon/2)^c}^z}] \\
&\quad + 2\mathbb{E}[\tau_{B_x(\varepsilon/2)^c}^z \mathbb{E}[\tau_x^y, \tau_x^y < \tau_{I\cup\Omega^c}^y]_{|y\equiv X_t^z, t\equiv \tau_{B_x(\varepsilon/2)^c}^z}] \\
&\quad + \mathbb{E}[\mathbb{E}[(\tau_x^y)^2, \tau_x^y < \tau_{I\cup\Omega^c}^y]_{|y\equiv X_t^z, t\equiv \tau_{B_x(\varepsilon/2)^c}^z}].
\end{aligned}
\tag{5.20}
$$

By the Cauchy inequality in combination with (2.21) and (2.18) for $h \equiv w_{B_x(\varepsilon/2)^c}^\lambda$, $K \equiv \Sigma \equiv \Gamma \equiv B_x(\varepsilon/2)$ and for $h \equiv h_{B_x(\varepsilon/2)^c}^\lambda$ the first term on the right-hand side is bounded by $C\mathbb{E}[\tau_{B_x(\varepsilon/2)^c}^z]/\lambda(B_x(\varepsilon/2)) \leq C$. For the second term we compute by the strong Markov property and Assumption 4.5 for some $N$ and all $y \in \partial B_x(\varepsilon/2)$ and $z \in \mathcal{M}\setminus I\setminus x$ satisfying $T_{z,I\cup x} = T_{I\cup x}$

$$
\begin{aligned}
\mathbb{E}[\tau_x^y, &\tau_x^y < \tau_{I\cup\Omega^c}^y] \\
&\leq \mathbb{E}[\tau_x^y, \tau_x^y = \tau_{\mathcal{M}\cup\Omega^c}^y] \\
&\quad + \mathbb{E}[\tau_z^y|\tau_z^y = \tau_{\mathcal{M}\cup\Omega^c}^y]\mathbb{P}[\tau_z^y = \tau_{\mathcal{M}\cup\Omega^c}^y]\sup_{u\in\partial B_z}\mathbb{P}[\tau_x^z < \tau_{I\cup\Omega^c}^z] \\
&\quad + \mathbb{P}[\tau_z^y = \tau_{\mathcal{M}\cup\Omega^c}^y]\sup_{u\in\partial B_z}\mathbb{P}[\tau_x^z < \tau_{I\cup\Omega^c}^z]\mathbb{E}[\tau_x^z|\tau_x^z < \tau_{I\cup\Omega^c}^z] \\
&\leq \varepsilon^{-N} + (\varepsilon^{-N} + C\varepsilon^{-N}T_{I\cup x})e^{-(\hat F(x,z)-F(x)+\hat F(z,x)-F(z))/\varepsilon}\big/\varepsilon^2 \\
&\leq \varepsilon^{-N}(1 + Ce^{-(\hat F(x,z)-F(x))/\varepsilon}),
\end{aligned}
\tag{5.21}
$$

where we have used (2.33) in the second inequality. For the third term on the right-hand side of (5.20) we obtain similarly to (5.21) by the strong Markov

property

$$
\begin{aligned}
&\mathbb{E}[(\tau_x^y)^2, \tau_x^y < \tau_{I \cup \Omega^c}^y] \\
&\leq \mathbb{E}[(\tau_x^y)^2, \tau_x^y = \tau_{\mathcal{M} \cup \Omega^c}^y] \\
&\quad + \mathbb{P}[\tau_z^y = \tau_{\mathcal{M} \cup \Omega^c}^y]\mathbb{E}[(\tau_z^y)^2 | \tau_z^y = \tau_{\mathcal{M} \cup \Omega^c}^y] \sup_{u \in \partial B_z} \mathbb{P}[\tau_x^u < \tau_{I \cup \Omega^c}^u] \\
&\quad + \mathbb{P}[\tau_z^y = \tau_{\mathcal{M} \cup \Omega^c}^y] \sup_{u \in \partial B_z} \mathbb{P}[\tau_x^u < \tau_{I \cup \Omega^c}^u]\mathbb{E}[(\tau_x^u)^2 | \tau_x^u < \tau_{I \cup \Omega^c}^u] \\
&\quad + 2\mathbb{P}[\tau_z^y = \tau_{\mathcal{M} \cup \Omega^c}^y]\mathbb{E}[\tau_z^y | \tau_z^y = \tau_{\mathcal{M} \cup \Omega^c}^y] \\
&\quad\quad \times \sup_{u \in \partial B_z} \mathbb{P}[\tau_x^u < \tau_{I \cup \Omega^c}^u]\mathbb{E}[\tau_x^u | \tau_x^u < \tau_{I \cup \Omega^c}^u],
\end{aligned}
\tag{5.22}
$$

so that in combination with the Cauchy inequality for $\lambda < \varepsilon^N / T_{I \cup x} \wedge (\lambda(\Omega \backslash B_{I \cup x})/C)$ and some $\mu < \varepsilon^N \wedge (\lambda(\Omega \backslash B_{\mathcal{M}})/C)$, Proposition 3.4 and (2.33)

$$
\begin{aligned}
&\mathbb{E}[(\tau_x^y)^2, \tau_x^y < \tau_{I \cup \Omega^c}^y] \\
&\leq \varepsilon^{-N}/\mu + (\varepsilon^{-N}/\mu + \varepsilon^{-N}T_{I \cup x}/\lambda)e^{-(\hat{F}(x,z) - F(x) + \hat{F}(z,x) - F(z))/\varepsilon}/\varepsilon^2 \\
&\leq \varepsilon^{-N} + \varepsilon^{-N}e^{-(\hat{F}(x,z) - F(x))/\varepsilon} + \varepsilon^{-N}e^{-\rho/\varepsilon}.
\end{aligned}
\tag{5.23}
$$

In case that $\mathcal{M} = I \cup x$ the bounds in (5.21) and (5.23) remain valid if the terms involving exponentials are replaced by zero by an even simpler argumentation and an obvious generalization of (2.33) to the case $I \cup J = \mathcal{M}$. Furthermore, for some $C \equiv C(d)$ and $\varepsilon > 0$ small enough, using (2.21) again,

$$
\begin{aligned}
\mathbb{E}[\sigma_1^z, \sigma_1^z < \rho_1^z] &\geq \mathbb{E}[\sigma_{B(x,\varepsilon/2)^c}^z] \inf_{y \in \partial B(x,\varepsilon/2)} \mathbb{P}[\tau_x^y < \tau_{I \cup \Omega^c}^y] \\
&\quad + \inf_{y \in \partial B(x,\varepsilon/2)} \mathbb{E}[\tau_x^y, \tau_x^y < \tau_{I \cup \Omega^c}^y] \\
&\geq \mathbb{E}[\sigma_{B(x,\varepsilon/2)^c}^z]/2 \geq \varepsilon/C
\end{aligned}
\tag{5.24}
$$

since the latter probability converges exponentially fast in $1/\varepsilon$ to 1. On the other hand, combination of (2.21) and (5.21) gives

$$
\begin{aligned}
\mathbb{E}[\sigma_1^z, \sigma_1^z < \rho_1^z] &\leq \mathbb{E}[\sigma_{B(x,\varepsilon/2)^c}^z] + \sup_{y \in \partial B(x,\varepsilon/2)} \mathbb{E}[\tau_x^y, \tau_x^y < \tau_{I \cup \Omega^c}^y] \\
&\leq C\varepsilon + \varepsilon^{-N}(1 + Ce^{-(\hat{F}(x,z) - F(x))/\varepsilon}).
\end{aligned}
\tag{5.25}
$$

Combination of (5.20) with the remark following, (5.21) and (5.23) and the resulting bound with (5.24), (5.25) and (5.14) shows that the right-hand side of (5.19) is bounded below by $\varepsilon^{-N}$ for some $N \equiv N(d, F)$ and all $k \geq \varepsilon^{-N}t$. □

With the uniform a priori estimate (5.11) we can proceed with the generalization of Proposition 6.1 in [2] to our setting. For $\beta \in (F(x), \infty)$ we denote by $\mathcal{C}_x(\beta)$ the connected component of $x$ in $\{F < \beta\}$.

LEMMA 5.4. *There is $C \equiv C(d, F) > 0$ such that for some $N \equiv N(d, F)$, all $\beta > F(x) + C\varepsilon \log(1/\varepsilon)$, all $t > 0$ and all $y \in B_x$,*

$$
\mathbb{P}[X_t^y \notin \mathcal{C}_x(\beta) | \tau_{I \cup \Omega^c}^y > t] \leq \varepsilon^{-N}T_{I \cup x}e^{-(\beta - F(x))/\varepsilon}.
\tag{5.26}
$$

PROOF. Let $T \equiv \inf_{z \in \partial B_x(\varepsilon/2)} \mathbb{E}[\sigma_1^z|\sigma_1^z < \rho_1^z]/2$, where we have chosen $A \equiv B_x(\varepsilon/2)^c$, $\Sigma \equiv B_x^c$ and $B \equiv B_I \cup \Omega^c$ in the definition (5.12). Decomposing the event $\{X_t^y \notin \mathcal{C}_x(\beta), \tau_I^y > t\}$ according to the number of returns to $B_x$ from $B_x(\varepsilon/2)^c$ before time $t$, we have for $K \equiv \min\{k \in \mathbb{N}|kT \geq t\}$

$$\begin{aligned}
&\mathbb{P}[X_t^y \notin \mathcal{C}_x(\beta), \tau_{I \cup \Omega^c}^y > t] \\
&\quad = \sum_{n \geq 0} \sum_{1 \leq k \leq K} \mathbb{P}[X_t^y \notin \mathcal{C}_x(\beta), T_n^y < t < T_{n+1}^y \wedge \tau_{I \cup \Omega^c}^y, \\
&\qquad\qquad\qquad\qquad\qquad\qquad\qquad \eta_n^y \in [(k-1)T, kT \wedge t)],
\end{aligned}$$
(5.27)

where $\eta_n^y$ denotes the first instant of reaching $B_x(\varepsilon/2)^c$ after the moment $T_n^y$ of the $n$th return to $B_x$ from $B_x(\varepsilon/2)^c$ before time $t$ when starting in $y$. For the $k$th term in the $n$th inner sum on the right-hand side of the last equation we may write using the strong Markov property

$$\begin{aligned}
&\mathbb{P}[X_t^y \notin \mathcal{C}_x(\beta), T_n^y < t < T_{n+1}^y \wedge \tau_{I \cup \Omega^c}^y, \eta_n^y \in [(k-1)T, kT \wedge t)] \\
&\quad = \mathbb{E}[\mathbb{P}[X_{t-r}^z \notin \mathcal{C}_x(\beta), \tau_{I \cup x \cup \Omega^c}^z > t - r]_{|z \equiv X_r^y, r \equiv \eta_n^y}, \\
&\qquad\qquad \eta_n^y \in [(k-1)T, kT \wedge t \wedge \tau_{I \cup \Omega^c}^y)],
\end{aligned}$$
(5.28)

where for some $C \equiv C(d, F)$, some $N \equiv N(d)$, all $\lambda < \varepsilon^N/T_{I \cup x}$, all $r \in [(k-1)T, kT)$ and all $z \in \partial B_x(\varepsilon/2)$ by the exponential Chebyshev inequality, the strong Markov property, (3.26) and Corollary 4.8 in [4],

$$\begin{aligned}
&\mathbb{P}[t - r < \tau_{I \cup x \cup \Omega^c}^z, X_{t-r}^z \notin \mathcal{C}_x(\beta)] \\
&\quad \leq e^{-\lambda(t-kT)} \mathbb{E}[e^{\lambda \tau_{I \cup x \cup \Omega^c}^z}, \tau_\beta^z \leq \tau_{I \cup x \cup \Omega^c}^z] \\
&\quad \leq e^{-\lambda(t-kT)} \mathbb{E}[e^{\lambda \tau_\beta^z}, \tau_\beta^z \leq \tau_{I \cup x \cup \Omega^c}^z] \sup_{y \in \partial \mathcal{C}_x(\beta)} \mathbb{E}[e^{\lambda \tau_{I \cup x \cup \Omega^c}^y}] \\
&\quad \leq e^{-\lambda(t-kT)} C \mathbb{P}[\tau_\beta^z \leq \tau_{I \cup x \cup \Omega^c}^z] \\
&\quad \leq e^{-\lambda(t-kT)} C e^{-(\beta - F(x))/\varepsilon}/\varepsilon,
\end{aligned}$$
(5.29)

where we have introduced $\tau_\beta^z \equiv \inf\{t > 0|F(x_t^z) \geq \beta\}$. Combination of this estimate with (5.28) and (5.27) leads to

$$\begin{aligned}
&\mathbb{P}[X_t^y \notin \mathcal{C}_x(\beta), \tau_{I \cup \Omega^c}^y > t] \\
&\quad \leq \frac{C}{\varepsilon} e^{-(\beta - F(x))/\varepsilon} \\
&\qquad \times \sum_{1 \leq k \leq K} e^{-\lambda(t-kT)} \sum_{n \geq 0} \mathbb{P}[\eta_n^y \in [(k-1)T, kT \wedge t \wedge \tau_{I \cup \Omega^c}^y)].
\end{aligned}$$
(5.30)

From the definition of $T$ in the beginning of the proof and the second moment inequality (5.15) it follows for some $N \equiv N(d, F)$ and all $z \in \partial B_x(\varepsilon/2)$ and all $r \in [(k-1)T, kT \wedge t \wedge \tau_{I \cup \Omega^c}^y)$

$$\begin{aligned}
\mathbb{P}[kT \wedge t - r < \sigma_1^z < \rho_1^z] &\geq \mathbb{P}[T < \sigma_1^z < \rho_1^z] \\
&\geq \mathbb{P}[\sigma_1^z < \rho_1^z] \frac{\mathbb{E}[\sigma_1^z|\sigma_1^z < \rho_1^z]^2}{4\mathbb{E}[(\sigma_1^z)^2|\sigma_1^z < \rho_1^z]} \geq \varepsilon^N,
\end{aligned}$$
(5.31)

where the last line involves a computation almost the same as in (5.14) and in (5.20) to (5.23) which we leave to the reader. Moreover, similarly to (5.24) the reader may convince himself that

$$\mathbb{E}[\sigma_1^z | \sigma_1^z < \rho_1^z] \geq \varepsilon^N \qquad (z \in \partial B_x(\varepsilon/2)) \tag{5.32}$$

after possibly increasing $N$. Using (5.31), we compute by the strong Markov property

$$\begin{aligned}
&\mathbb{P}[(k-1)T \leq \eta_n^y < kT \wedge t < \eta_{n+1}^y < \tau_{I \cup \Omega^c}^y] \\
&= \mathbb{E}[\mathbb{P}[kT \wedge t - r < \sigma_1^z < \rho_1^z]|_{z \equiv X_r^y, r \equiv \eta_n^y}, \\
&\qquad \eta_n^y \in [(k-1)T, kT \wedge t \wedge \tau_{I \cup \Omega^c}^y)] \\
&\geq \varepsilon^N \mathbb{P}[\eta_n^y \in [(k-1)T, kT \wedge t \wedge \tau_{I \cup \Omega^c}^y)].
\end{aligned} \tag{5.33}$$

Let $\eta_x^y$ be the first time $s$ after $(k-1)T$ such that $X_s^y$ reaches $B_x(\varepsilon/2)$. Combination of (5.32) and (5.33) with (5.30) shows for some $N$

$$\begin{aligned}
&\mathbb{P}[X_t^y \notin \mathcal{C}_x(\beta), \tau_{I \cup \Omega^c}^y > t] \\
&\leq \varepsilon^{-N} e^{-(\beta - F(x))/\varepsilon} \\
&\quad \times \sum_{1 \leq k \leq K} e^{-\lambda(t-kT)} \mathbb{P}[(k-1)T \leq \eta_x^y < kT \wedge t \wedge \tau_{I \cup \Omega^c}^y].
\end{aligned} \tag{5.34}$$

On the other hand, we have by the strong Markov property again in combination with the a priori lower bound (5.11)

$$\begin{aligned}
\mathbb{P}[\tau_{I \cup \Omega^c}^y > t] &\geq \mathbb{P}[(k-1)T \leq \eta_x^y < kT \wedge t < \tau_{I \cup \Omega^c}^y] \\
&= \mathbb{E}[\mathbb{P}[\tau_{I \cup \Omega^c}^y > t - r]|_{z \equiv X_r^y, r \equiv \eta_x^y}, \\
&\qquad (k-1)T \leq \eta_x^y < kT \wedge t \wedge \tau_{I \cup \Omega^c}^y] \\
&\geq e^{-\varepsilon^{-N}(t-(k-1)T)/T_I} \varepsilon^N \mathbb{P}[(k-1)T \leq \eta_x^y < kT \wedge t \wedge \tau_{I \cup \Omega^c}^y].
\end{aligned} \tag{5.35}$$

The last two estimates, the choice of $\lambda$ and $T$ in combination with (5.32) prove (5.26). $\square$

Now we are in a position to estimate the second and third terms on the right-hand side of (5.6). For $y \in B_x$ combination of (5.26) and (3.20) leads to

$$\begin{aligned}
&\mathbb{E}[\Pi_x(\mathbb{1}_{\Omega \setminus \overline{B}_I})(X_t^y), X_t^y \in \Omega_1 \setminus A_{x,I}(\beta), \tau_{I \cup \Omega^c}^y > t] \\
&\quad + \sum_{k \in \mathbb{N}} \mathbb{E}[\Pi_x(\mathbb{1}_{\Omega \setminus \overline{B}_I})(X_t^y), X_t^y \in \Gamma_k, \tau_{I \cup \Omega^c}^y > t] \\
&= \mathcal{O}(1) \varepsilon^{-N} \mathbb{P}[\tau_{I \cup \Omega^c}^y > t] \\
&\quad \times \Bigg( \frac{T_{I \cup x}}{T_I} + T_{I \cup x} e^{-(\hat{F}(x,I) - F(x) - \beta)/\varepsilon} \\
&\qquad + T_{I \cup x} e^{-(C_1 - F(x))/\varepsilon} \\
&\qquad \times \int_{\{F > C_1\}} \sup_{B(y, \varepsilon \delta(y))} |\nabla F|^d e^{-(F(y) - C_1)/(2\varepsilon) - \operatorname{dist}(y, \Omega_1)/C} \, dy \Bigg)
\end{aligned} \tag{5.36}$$

for sufficiently large $C$ and $\beta > C\varepsilon \log(1/\varepsilon)$. Here we also have used the trivial bound $(\mathbb{1}_{\{\phi_x > 1\}}, \phi_x)_F / \|\phi_x\|_F^2 \leq 1$. It is not difficult to see that by the integrability condition on $F$ the latter integral on the right-hand side of the last display is bounded uniformly in $\varepsilon > 0$ small enough. The theorem is proven since in the bounded case we do not need (3.20). $\square$

MATHEMATISCHES INSTITUT ZÜRICH
UNIVERSITÄT ZÜRICH
WINTERTHURERSTRASSE 190
CH-8057 ZÜRICH
SWITZERLAND
E-MAIL: eckhoff@amath.unizh.ch
URL: www.math.unizh.ch/assistenten/eckhoff